\documentclass[10pt]{article}
\usepackage[all,2cell]{xy} \UseAllTwocells \SilentMatrices
\usepackage{latexsym,amsfonts,amssymb}
\usepackage{amsmath,amsthm,amscd}
\usepackage{hyperref,psfrag}
\usepackage{color}
\usepackage{etoolbox}
\usepackage[dvips]{epsfig}
\usepackage{psfrag}
\usepackage{pinlabel}

\usepackage{graphicx}

\usepackage{a4wide}
\usepackage{geometry}

\usepackage{epigraph,wrapfig}

\normalfont\upshape

\usepackage{fancyheadings}
\theoremstyle{definition}
\newtheorem{thm}{Theorem}[section]
\newtheorem{cor}[thm]{Corollary}
\newtheorem{conj}[thm]{Conjecture}
\newtheorem{lem}[thm]{Lemma}
\newtheorem{rem}[thm]{Remark}
\newtheorem{prop}[thm]{Proposition}
\newtheorem{defn}[thm]{Definition}
\newtheorem{example}[thm]{Example}

\numberwithin{equation}{section}

\newcommand{\ang}[1]{{\langle{#1}\rangle}}

\usepackage{bbm}
\def\ot{\otimes}
\def\C{{\mathbbm C}}
\def\N{{\mathbbm N}}

\def\Z{{\mathbbm Z}}
\def\Q{{\mathbbm Q}}

\def\F{{\mathbbm F}}

\def\1{{\mathbbm{1}}}

\newcommand{\Hom}{{\rm Hom}}

\renewcommand{\to}{\rightarrow}

\newcommand{\co}{\colon}

\newcommand{\id}{{\rm id}}

\newcommand{\pa}{\partial}

\newcommand{\End}{{\rm End}}

\def\sltwo{\mathfrak{sl}_2}

\def\spmod{{\mathrm{\mbox{-}spmod}}}

\def\shuffle{\,\raise 1pt\hbox{$\scriptscriptstyle\cup{\mskip
               -4mu}\cup$}\,}

\newcommand{\refequal}[1]{\xy {\ar@{=}^{#1}
(-1,0)*{};(1,0)*{}};
\endxy}

\def\BM{\mathbbm{B}}

\newcommand{\ig}[2]{\vcenter{\xy (0,0)*{\includegraphics[scale=#1]{fig/#2}} \endxy}}

\newcommand{\dup}{\mathbf{d}}
\newcommand{\ddown}{\mathbf{z}}
\renewcommand{\dh}{\mathbf{h}}
\newcommand{\dx}{\mathbf{x}}
\newcommand{\dy}{\mathbf{y}}
\newcommand{\dz}{\mathbf{z}}
\newcommand{\bardup}{\overline{\mathbf{d}}}

\DeclareMathOperator{\DL}{DL}
\DeclareMathOperator{\LL}{LL}
\DeclareMathOperator{\TL}{LT}

\newcommand{\uw}{\underline{w}}
\newcommand{\ux}{\underline{x}}
\newcommand{\uy}{\underline{y}}
\newcommand{\uv}{\underline{v}}
\newcommand{\uz}{\underline{z}}
\newcommand{\eb}{\mathbf{e}}
\newcommand{\fb}{\mathbf{f}}

\newcommand{\OC}{\mathcal{O}}
\newcommand{\HC}{\mathcal{H}}
\newcommand{\CC}{\mathcal{C}}

\newcommand{\EC}{\mathcal{E}}
\DeclareMathOperator{\NH}{NH}

\DeclareMathOperator{\BS}{BS}

\DeclareMathOperator{\Span}{Span}
\DeclareMathOperator{\Mat}{Mat}

\DeclareMathOperator{\core}{Core}
\DeclareMathOperator{\oops}{oops}
\newcommand{\sqot}{\boxtimes}
\newcommand{\UC}{\mathcal{U}}
\newcommand{\FC}{\mathcal{F}}
\DeclareMathOperator{\Ker}{Ker}

\newcommand{\gln}{\mathfrak{gl}_n}
\DeclareMathOperator{\Jac}{Jac}
\DeclareMathOperator{\semis}{ss}

\newcommand{\un}{\underline}
\renewcommand{\gg}{\mathfrak{g}}
\newcommand{\bb}{\mathfrak{b}}

\newcommand{\smsum}{{\scriptstyle \Sigma}}
\newcommand{\ui}{\underline{i}}
\newcommand{\uj}{\underline{j}}
\newcommand{\GG}{\Gamma\Gamma}

\renewcommand{\Im}{\mathrm{Im}}

\newcommand{\poly}[1]{{
\labellist
\small\hair 2pt
 \pinlabel {$#1$} [ ] at 7 15
\endlabellist
\centering
\ig{1}{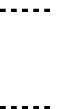}
}}

\DeclareMathOperator{\enrich}{enrich}
\newcommand{\sumset}{\stackrel{\oplus}{\subset}}

\title{Actions of $\sltwo$ on algebras appearing in categorification}
\author{Ben Elias, You Qi}

\date{\today}

\begin{document}
%

\maketitle

\begin{abstract} We prove that many of the recently-constructed algebras and categories which appear in categorification can be equipped with an action of $\sltwo$ by
derivations. The $\sltwo$ representations which appear are filtered by tensor products of coverma modules. In a future paper, we will address the implications of the $\sltwo$ structure for categorification. \end{abstract}

\setcounter{tocdepth}{2} \tableofcontents


\section{Introduction}

\subsection{$\sltwo$ structures}

Here are several $2$-categories with graded $2$-morphism spaces which play fundamental roles in categorical representation theory in type $A$.
\begin{itemize} \item The Khovanov-Lauda-Rouquier category $\UC^+(\gln)$, which categorifies the positive half of the quantum group of $\gln$ \cite{KhoLau09, Rouq2KM-pp}.
	\item The Lauda's category $\UC(\sltwo)$, which categorifies the entire quantum group of $\sltwo$ \cite{LauSL2}.
	\item The thickened category $\dot{\UC}(\sltwo)$ of Khovanov-Lauda-Mackaay-Stosic, which also categorifies the entire quantum group of $\sltwo$ \cite{KLMS}.
	\item The diagrammatic Hecke category $\HC(S_n)$, which categorifies the Iwahori-Hecke algebra of $S_n$ \cite{EKho}.
\end{itemize}
They are each defined by generators and relations over some base ring, using the technology of planar diagrammatics. Since the relations only have integral coefficients, we can and will assume the base ring is the integers.

The first main theorem of this paper says that these categories have a surprising new structure.

\begin{thm} \label{thm:sl2action} Each of the categories listed above admits an action of $\sltwo$ by derivations, compatible with the graded $2$-category structure. Moreover, divided powers of the raising and lowering operators act on the integral form. \end{thm}

More precisely, we will define three operators $\{\dup, \dh, \ddown\}$ on the morphism spaces in these categories, and each operator $\dx \in \{\dup, \dh, \ddown\}$ will satisfy the Leibniz rule
\begin{equation} \dx(f \circ g) = \dx(f) \circ g + f \circ \dx(g), \quad \dx(f \ot g) = \dx(f) \ot g + f \ot \dx(g) \end{equation}
with respect to both vertical ($\circ$) and horizontal ($\otimes$) composition. The triple $(\dup, \dh, -\ddown)$ will act as an $\sltwo$ triple on each morphism space. A graded preadditive category equipped with an action of $\sltwo$ by derivations will be called an \emph{$\sltwo$-category}. An $\sltwo$-category is analogous to a dg-category, but for an unusual kind of homological algebra; morphism spaces in both categories are modules over a Hopf algebra. The concept of an $\sltwo$-category does not fit precisely into the framework of Hopfological algebra \cite{KhoHopf, QiHopf}, but we hope to address this in a follow-up paper.

Let us discuss these three operators in turn.

The degree $+2$ operator $\dup$ has been the central topic of study in the recent programme which attempts to categorify key objects in representation theory (e.g. quantum groups
and Hecke algebras) at a root of unity. For each of the categories above, the operator $\dup$ has been defined in previous works \cite{KQ, EQpDGsmall, EQpDGbig, EQHecke}. There is a
large family of degree $+2$ derivations one could place on each of these categories, but $\dup$ is unique (up to duality, see Remark \ref{rmk:duality}) in satisfying certain key
properties important for categorification. There is still no geometric understanding for the existence and importance of $\dup$, though see \cite{BeliakovaCooper} for more on the
connection to Steenrod operations. We will not discuss categorification at a root of unity any further in this paper, though we will have much to say in the next paper.

The degree $0$ operator $\dh$ is the degree operator. It multiplies any homogeneous morphism by a scalar, equal to its degree. Thus the weight grading for $\sltwo$ matches the ordinary grading in these categories.

The degree $-2$ operator $\ddown$ is new in this paper. It is also remarkably easy to define, once you know that it exists. Most of the generating morphisms of these categories live
in the minimal degree in their respective morphism spaces, and $\ddown$ must send them to zero for degree reasons. From these considerations, it is rather simple to verify Theorem
\ref{thm:sl2action}, and even to prove that $\ddown$ is the unique derivation of degree $-2$ (up to scalar). Note that $-\ddown$ is the lowering operator in the $\sltwo$ triple,
while we prefer to discuss $\ddown$ because it eliminates many signs from the formulas.

To reiterate, proving Theorem \ref{thm:sl2action} is quite easy. However, as far as we are aware there was no expectation at all that these categories should admit $\sltwo$ actions,
and no one had bothered to look for a degree $-2$ derivation. We discovered it by accident, as part of an effort to explain certain ``hard Lefschetz style'' phenomena which appeared
in the study of $\dup$.

\begin{rem} \label{rmk:duality} Each of these categories admits a duality functor, a contravariant automorphism $f \mapsto \bar{f}$ which flips diagrams upside-down. The operators
$\dh$ and $\ddown$ both intertwine with duality, but $\dup$ does not. Instead, duality intertwines $\dup$ with another derivation $\bardup$, where
\begin{equation}\bardup(f) := \overline{(\dup(\bar{f}))}.\end{equation}
Then $(\bardup, \dh, -\ddown)$ is another $\sltwo$ triple acting on the category. Readers familiar with the Jacobson-Morozov theorem might be surprised by the existence of two different $\sltwo$ triples which share the same lowering and degree operators, but we should reiterate that these representations of $\sltwo$ are all infinite-dimensional. \end{rem}

\begin{rem} In \S\ref{sec:KLR}, we actually equip the Khovanov-Lauda-Rouquier category $\UC^+(\gg)$ with an $\sltwo$ action for any (oriented) simply laced root datum. Meanwhile, it
was proven in \cite[Proposition 6.8]{EQHecke} that there is no ``good'' degree $+2$ derivation on the Hecke category in simply laced types outside of finite and affine type $A$, so
we do not expect an $\sltwo$ action either. We have not yet tried to place an $\sltwo$ action on other categorifications of interest though we expect the phenomenon to be common in
certain categories related to equivariant constructible sheaves, see Remark \ref{rem:equivariant}. We certainly expect there to be an $\sltwo$ action on $\UC(\mathfrak{gl}_n)$ and
$\UC(\mathfrak{sl}_n)$ in general, and it is reasonable to expect that one exists on quantum Heisenberg 2-categories (see e.g. \cite[\S 4]{LicataSavageHeisenberg}). \end{rem}

\subsection{The polynomial ring}

If $A$ and $B$ are algebras equipped with an action of $\sltwo$ by derivations, and $M$ and $N$ are bimodules equipped with a compatible action of $\sltwo$, then the space of
bimodule morphisms $\Hom(M,N)$ is naturally equipped with an $\sltwo$ action as well. This is analogous to the internal Hom between two chain complexes, which is itself another
chain complex.

There is a well-known action of $\sltwo$ on the ring $R_n = \Z[x_1, \ldots, x_n]$, where $\deg x_i = 2$ for all $i$ (which determines the $\dh$ action), and where \begin{equation}
\label{actiononRn} \dup = \sum_i x_i^2 \frac{\partial}{\partial {x_i}}, \qquad \ddown = \sum_i \frac{\partial}{\partial {x_i}}. \end{equation} Here, $\sltwo$ appears as a subalgebra
of the Witt Lie algebra acting on polynomials\footnote{The Witt Lie algebra is generated by differential operators $L_k = x^{k+1} \frac{\partial}{\partial x}$ for all $k \in \Z$,
and it acts on the space of Laurent polynomials $\Z[x,x^{-1}]$. The subalgebra generated by $L_k$ for $k \ge -1$ preserves the subring of ordinary polynomials $\Z[x]$. The action on
a polynomial ring in $n$ variables is just the $n$-fold tensor product of the action on the polynomial ring in one variable. Note that the operators $L_k$ are quite different from the divided powers $\dup^{(k)}$ or $\ddown^{(-k)}$.}. Note that $\ddown$ acts trivially on the subring
generated by the roots $(x_i - x_j)$, and $\dup$ does not preserve this subring. For all the categories above, morphism spaces have polynomial subalgebras which play an important
role. In all cases, these polynomial subalgebras are preserved by the $\sltwo$ action, and the two different $\sltwo$ triples $(\dup, \dh, -\ddown)$ and $(\bardup, \dh, -\ddown)$
restrict to the same standard $\sltwo$ triple on the polynomial ring.


Most of the categories we are discussing have full faithful embeddings into the category of bimodules over polynomial rings. For example, the Hecke category $\HC$ is equivalent
(after base change) to the category of Soergel bimodules, certain $(R_n, R_n)$ bimodules. To give another example, the nilHecke algebra $\NH_n$ is isomorphic to
$\End_{R_n^{S_n}}(R_n)$. However, we wish to emphasize that this realization in terms of bimodules does not equip $\NH_n$ or $\HC$ with an $\sltwo$ action! One must still choose an
action of $\sltwo$ on the bimodules in question; this choice is not unique and is rather subtle in practice. Even a free module of rank $1$ over $R_n$ admits many compatible
$\sltwo$ actions.

\begin{example} In \cite[equations (65) and (66) on page 44]{KQ}, a one-parameter family of degree $+2$ derivations $\dup_a$ is defined on $\NH_n$. We define the lowering operator $\ddown$ on $\NH_n$ below in \eqref{loweringNHintro}. For any scalar $a$, $(\dup_a, \dh, -\ddown)$ is an $\sltwo$ triple. A major point in \cite{KQ} is that only two of these raising operators, $\dup_{\pm 1}$, have the desired properties for categorification. The subtle properties of $\dup_a$ are invisible from the perspective of the polynomial ring $R_n$, since $(\dup_a, \dh, -\ddown)$ restricts to the standard $\sltwo$ triple on $R_n$ for all $a$.

\end{example}

\begin{rem} In \cite{KRWitt}, Khovanov and Rozansky use the action of the positive half\footnote{This positive half includes all operators $L_k$ for $k \ge 0$. Their action
extends to include $L_{-1}$ as well, though they did not note this.} of the Witt Lie algebra on $R_n$ to place an action of this same algebra on triply graded knot homology, which
is built using Hochschild homology of Soergel bimodules. This (positive half of the) Witt action is an important precursor to our $\sltwo$ action, though just as in the previous remark, an action on the polynomial ring does not determine a unique action on the Hecke category itself.

\end{rem}

\begin{rem} \label{rem:equivariant} The $\sltwo$ action (even on the polynomial ring) currently lacks a geometric motivation. The raising operator is related to general homological
operations (e.g. Steenrod squares), though the connection is subtle, see \cite{BeliakovaCooper} (and also see \cite{KitchlooSteenrod} for more on Steenrod operations and Soergel
bimodules). However, the lowering operator does not seem to arise from a general construction. Though a good explanation is missing, one thing is clear: $\sltwo$ acts on the
$\C^*$-equivariant cohomology of a point. Geometric constructions of these categories involve perverse sheaves which are equivariant over an algebraic group. We suspect that there
is a relationship between the $\sltwo$ action and the existence of a copy of $\C^*$ (i.e. $\mathbb{G}_m$) inside the algebraic group which is acting trivially. Insisting upon
equivariance for a trivial action often leads to extra ``homological'' operations of higher degree, such as the ``log of monodromy'' maps from \cite{BezYunMonodromy}.

For example, the Hecke category studies $B \times B$-equivariant sheaves on $G$, and any element of $Z(G) \cap B$ will act the same way on both sides, so its antidiagonal copy in $B
\times B$ will act trivially. When $G = GL_n(\C)$, there is a central copy of $\C^*$ in the torus, and there is also an $\sltwo$ action on the Hecke category. When $G = SL_n(\C)$, there is no center, and $\sltwo$ does not act (c.f. \cite[Proposition 6.9]{EQHecke})! As a shadow of this fact, the reader can already verify that $\sltwo$ acts on $\Z[x_1, \ldots, x_n]$, but it does not have a (nontrivial) $S_n$-invariant action on the subring generated by $(x_i - x_j)$ for $i < j$. This contrasts the $B$-equivariant cohomology of a point for the Borel subgroup in $GL_n(\C)$ versus $SL_n(\C)$. \end{rem}

\begin{rem} There are a number of (typically non-monoidal) categories which also play major roles in categorical representation theory in type $A$, and for which the operator $\dup$ has already been studied \cite{KQS, QiSussan3}.
	\begin{itemize}
		\item The cyclotomic quotients of Khovanov-Lauda-Rouquier categories, which categorify irreducible representations of the quantum group.
		\item Webster's categories, which categorify tensor products of irreducible representations.
	\end{itemize}
However, neither the cyclotomic quotients nor Webster's categories admit actions of $\sltwo$. For example, cyclotomic quotients are quotients by an ideal inside $\UC^+(\gln)$ which is preserved by $\dup$ but not by $\ddown$. This suggests some nuance in how one should interpret modules over $\sltwo$-categories. \end{rem}

\begin{rem} Other important categories are the cell subquotients of $\HC(S_n)$, which categorify the irreducible modules over the Hecke algebra, and the Schur quotients of quantum group
categorifications \cite{MSV}. They arise from monoidal ideals generated by identity maps of various objects, and identity maps are killed by the $\sltwo$ action. Thus the $\sltwo$
action on $\HC(S_n)$ descends to its cell subquotients, and the $\sltwo$ action on $\UC(\mathfrak{gl}_2)$ descends to the Schur quotient. \end{rem}

\subsection{Contrasting algebra and representation theory} \label{subsec:NHintro}

Given an $\sltwo$-category, we can forget some structure and study it in two ways:
\begin{itemize} \item Forgetting the $\sltwo$ action, we can study the category algebraically. We can study the splitting of objects into direct summands, the Jacobson radical, and so forth. \item Forgetting the algebra structure, we can study Hom spaces as representations of $\sltwo$. We can ask about their characters, their finite-dimensional subrepresentations, and so forth. \end{itemize}
There seems to be an incredible connection between the structure of these categories as algebras and as $\sltwo$ representations. It is so astounding to the authors that it needs to be showcased immediately. We hope this example will whet the appetite, and drum up excitement for our $\sltwo$ action.

First we present a toy example. For ease
of discussion\footnote{We could make most of the same statements over $\Z$, but our use of terms like the Jacobson radical would be inappropriate.} let us work over a field $\Bbbk$ of characteristic zero. The ring $\Bbbk[x]$ has graded Jacobson radical $(x)$, and the quotient by this ideal is $\Bbbk$. So we have a short exact sequence
\begin{equation} \label{introseq1} 0 \to (x) \to \Bbbk[x] \to \Bbbk \to 0 \end{equation}
of $\Bbbk[x]$-modules, and $\Bbbk$ is the graded semisimplification of $\Bbbk[x]$. The submodule $(x)$ is also preserved by the raising operator $\dup$, so this is a short exact sequence of $U(\bb^+)$-modules, where $\bb^+$ is the Lie algebra inside $\sltwo$ generated by $\dup$ and $\dh$. However, $(x)$ is not preserved by $\ddown$. On the other hand, $\Bbbk \subset \Bbbk[x]$ is a subalgebra, and is also preserved by the $\sltwo$ action. So we have a short exact sequence
\begin{equation} \label{introseq2} 0 \to \Bbbk \to \Bbbk[x] \to Q \to 0 \end{equation}
of $\sltwo$-modules, where $Q$ is the quotient module. Note that $Q$ is simple, so it has no finite-dimensional submodules.

The short exact sequences \eqref{introseq1} and \eqref{introseq2} live in different categories, but they are both sequences of vector spaces. They split each other, in that the first map of \eqref{introseq2} will give a section for the quotient map of \eqref{introseq1}. Consequently, we can identify $Q$ with $(x)$ as a vector space and as a complementary direct summand to $\Bbbk$. In this fashion, the maximal finite-dimensional $\sltwo$-submodule $\Bbbk$ forms a semisimple subalgebra which maps isomorphically to the semisimplification (i.e. the quotient by the Jacobson radical).

Let us reproduce the same behavior in a more interesting example. The nilHecke algebra $\NH_n$ is the endomorphism algebra of the object $E^n$ inside $\UC^+(\sltwo)$. Put together, the nilHecke algebras form a monoidal category with morphism algebra $\NH = \bigoplus_{n
\ge 0} \NH_n$, which is monoidally generated by morphisms depicted as a dot and a crossing. The $\sltwo$-module structure is defined on the generators below, and is extended to the whole category using the Leibniz rule.
\begin{subequations}
\begin{equation} \dup \left(~\ig{1}{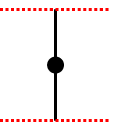} ~\right) = {
\labellist
\small\hair 2pt
 \pinlabel {$2$} [ ] at 7 19
\endlabellist
\centering
\ig{1}{dot}
}, \qquad \dup \left(~\ig{1}{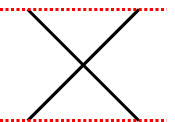}~ \right) = - \ig{1}{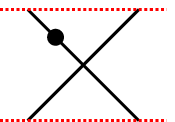} - \ig{1}{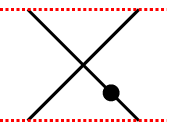}, \end{equation}
\begin{equation} \label{loweringNHintro} \ddown \left(~\ig{1}{dot} ~\right) = \ig{1}{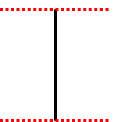}, \qquad \ddown \left(~\ig{1}{Xii} ~\right) = 0. \end{equation}
\end{subequations}

It is well-known that $\NH_n$ is the endomorphism algebra of the polynomial ring 
$$R = R_n = \Bbbk[x_1, x_2, \ldots, x_n]$$
 over its subalgebra $R^{S_n}$ of invariant polynomials. The $\sltwo$-action on $R_n$ is $S_n$-equivariant and thus descends to an $\sltwo$ action on $R^{S_n}$. By the Chevalley theorem, $R$ is free over $R^{S_n}$ of rank $n!$. Choosing a basis, one obtains an isomorphism between $\NH_n$ and a matrix algebra
\begin{equation} \NH_n \cong \Mat_{n!}(R^{S_n}). \end{equation}
On \cite[page 56]{KQ} one can find a basis of $\NH_n$ which corresponds to the basis of matrix entries, and results about the action of $\dup$ on this basis.

The (graded) Jacobson radical of $R^{S_n}$ is the ideal $R^{S_n}_+$ spanned by positive degree elements. The Jacobson radical of $\NH_n$ is therefore
\begin{equation} \Jac(\NH_n) \cong \Mat_{n!}(R^{S_n}_+). \end{equation}
Consequently there is a short exact sequence
\begin{equation} \label{ssofNHnMod} 0 \to \Jac(\NH_n) \to \NH_n \to \semis(\NH_n) \to 0, \end{equation}
where $\semis(\NH_n) \cong \Mat_{n!}(\Bbbk)$ is the semisimplification of $\NH_n$, viewed as a quotient. This is a short exact sequence of $\NH_n$-modules.

Now let us examine the $\sltwo$ structure on $\NH_n$. The ideal $R^{S_n}_+ \subset R^{S_n}$ is preserved by $\dup$ and $\dh$ but not by $\ddown$, since $\ddown(e_1) =
n \cdot \id$. Consequently, $\Jac(\NH_n)$ is preserved by $\dup$ and $\dh$ but not by $\ddown$, so there is no induced $\sltwo$ structure on $\semis(\NH_n)$. The short exact sequence \eqref{ssofNHnMod} is not a short exact sequence of $\sltwo$ representations.  However, $\Mat_{n!}(\Bbbk)$ is not just a quotient of $\NH_n \cong \Mat_{n!}(R^{S_n})$, it is also a subring. See \S\ref{subsec-NH} for the proof of the following theorem.

\begin{thm} \label{thm:NHnintro} The maximal finite-dimensional $\sltwo$-subrepresentation of $\NH_n$ is a subalgebra, which maps isomorphically onto the semisimplification of $\NH_n$. Moreover, under the identification of $\NH_n$ with $\Mat_{n!}(R^{S_n})$ used in \cite[Proposition 3.24]{KQ}, this subalgebra is precisely $\Mat_{n!}(\Bbbk)$. \end{thm}


\begin{example} \label{ex:NH2} When $n=2$, the following morphisms correspond to the matrix entries in $\Mat_2(\Bbbk)$.
\begin{equation} \label{NH2matrix} \left( \begin{array}{cc} \ig{1}{XiiNW} & -\ig{1}{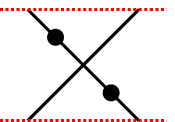} \\ \\ \ig{1}{Xii} & -\ig{1}{XiiSE} \end{array} \right) \end{equation}	
The reader should confirm that this four-dimensional subspace
of $\NH_2$ is preserved by $\ddown$ and $\dup$, and is isomorphic as an $\sltwo$-representation to $V \ot V^*$, where $V$ is the standard representation of $\sltwo$. \end{example}

Thus one has a short exact sequence of $\sltwo$-representations \begin{equation} \label{ssofNHnsl2} 0 \to \Mat_{n!}(\Bbbk) \to \NH_n \to Q \to 0 \end{equation} where $Q$ is defined
as this quotient. This is not a short exact sequence of $\NH_n$-modules, since $\Mat_{n!}(\Bbbk)$ is a subring but not an ideal. However, both sequences \eqref{ssofNHnMod} and
\eqref{ssofNHnsl2} are sequences of $\Bbbk$-modules, and split each other. Thus we can identify $Q$ with the Jacobson radical as a vector space complementary to $\Mat_{n!}(\Bbbk)$.
Said another way, the finite-dimensional part of the $\sltwo$-representation $\NH_n$ is precisely a complement to the Jacobson radical!

To summarize, there is a splitting of $U(\bb^+)$-modules
\begin{equation} \label{eq:splittingofNH} \NH_n = \Mat_{n!}(R^{S_n}) \cong \Mat_{n!}(R^{S_n}_+) \oplus \Mat_{n!}(\Bbbk). \end{equation}
When viewed as modules over $\NH_n$ or its subring $R$, this is not a splitting but a filtration, with $\Mat_{n!}(R^{S_n}_+)$ being the submodule. When viewed as modules over $\sltwo$, this is not a splitting but a filtration, with $\Mat_{n!}(\Bbbk)$ being the submodule.

\begin{rem} We come across this situation often in this paper. A vector space (resp. $\Z$-module) $X$ is isomorphic to a direct sum $A \oplus B \oplus C$ as vector spaces. However,
when $X$ is equipped with additional structure (e.g. an action of $\sltwo$, or of $R_n$), there is no longer a direct sum decomposition. Instead, the subspaces $A$ and $A \oplus B$
are preserved by the action. In this context we use the slightly imprecise phrasing that the \emph{splitting becomes/is a filtration} in the presence of the extra structure. We
call $A$ the subobject and $C$ the quotient object. The objects $A$, $B$, and $C$ were not equipped a priori with any additional structure, and they inherit this structure only by
the canonical identification of vector spaces with subquotients in the filtration $0 \subset A \subset A \oplus B \subset X$. \end{rem}

Note that $\NH_n$ is infinite-dimensional and even infinitely-generated as an $\sltwo$ module (when $n > 1$), so this kind of representation does not conform to most familiar regimes
(though they are direct limits of representations in category $\OC$). Thankfully, $\NH_n$ has finite-dimensional weight spaces, with weights bounded below. It is easy to prove that such
an $\sltwo$-representation contains a unique maximal finite-dimensional subrepresentation, which we call its \emph{core}. Equivalently, the core consists of all vectors on which $\dup$
acts nilpotently. From the Leibniz rule, one can see that the core must be closed under multiplication. Basic facts about the core are proven in \S\ref{sec-basics}. The example of
$\NH_n$ is supposed to demonstrate that the core of an $\sltwo$-category has dramatic significance to the algebraic structure of the category.

One nice feature of the core is that it can be found using basic linear algebra. One need only compute the kernel of $\dup$, and then use the lowering operator $\ddown$ to produce
the rest. This can be done in each Hom space independently, and does not involve the composition of morphisms. This simplicity is in contrast to other attempts to find complements
to the Jacobson radical, by computing inclusion and projection maps to indecomposable summands. This involves much more complicated linear algebra and in-depth knowledge of the category.

We do not wish the reader to expect that the core of an $\sltwo$-algebra always projects isomorphically to the semisimplification, as this is false in more complicated examples, see \S\ref{subsec:iteratedcoreexamples}. However, we conjecture that the map from the core to the semisimplification is injective, at least in geometrically-motivated $2$-categories.

\begin{conj} \label{conj:coreissplit} In the $2$-categories discussed in this paper, the core of each morphism space intersects the Jacobson radical trivially (every morphism in the
core is split), so the map from the core to the semisimplification is injective. \end{conj}

For example, one implication of this conjecture is that any degree zero morphism in the kernel of $\dup$ must split. If $p$ is the projection map to an indecomposable direct summand and $\dup$ acts nilpotently on $p$, then every nonzero $\dup^k(p)$ is also a projection map.

\begin{rem} This conjecture is already quite deep, indicating the non-degeneracy of certain ``local intersection forms'' upon restriction to finite-dimensional $\sltwo$
representations. One expects the eventual proof of this conjecture to come from a new kind of Hodge-theoretic argument. In all the examples computed to date, local intersection
forms have satisfied (an analogue of) the Hodge-Riemann bilinear relations upon restriction to the core. Note that the raising operator $\dup$ differs from the Lefschetz operator of
left multiplication studied in \cite{EWHodge}, and $\dup$ is not even a Lefschetz operator in the traditional sense. \end{rem}

\subsection{Filtrations on morphism spaces} \label{subsec:happyintro}

Hopefully, we have convinced the reader that the study of these categories as $\sltwo$-modules, and in particular the study of their finite-dimensional submodules, is of great
interest. Morphism spaces in these categories are free modules over a polynomial ring, as well as being modules over $\sltwo$. This is a useful tool in our effort to
understand the $\sltwo$-module structure on these morphism spaces.

\begin{defn} \label{defn:polySLTDGA} The \emph{polynomial $\sltwo$-algebra} $(R_n,\sltwo)$ is the polynomial ring $R_n = \Z[x_1, \ldots, x_n]$ equipped with the
$\sltwo$ action given in \eqref{actiononRn}. We also let $R_n$ denote the base change of $R_n$ to any commutative base ring $\Bbbk$.

An \emph{$(R_n, \sltwo)$-module} is an $R_n$-module $M$ which is also an $\sltwo$-module, satisfying a Leibniz rule. For $\dx \in \{\dup, \dh, \ddown\}$, if we write $\dx_M$ for the
action on $M$ and $\dx_R$ for the action on $R_n$, the Leibniz rule states that \begin{equation} \dx_M(r \cdot m) = \dx_R(r) \cdot m + r \cdot \dx_M(m). \end{equation} \end{defn}

Now we ask: what kinds of $(R_n, \sltwo)$-modules appear as morphism spaces in categories of interest? We know that morphism spaces will be free over $R_n$, meaning that they are a
direct sum of rank 1 free modules, but one should not expect such a splitting to be respected by the $\sltwo$ structure. The punchline will be: morphism spaces have $(R_n,
\sltwo)$-filtrations which are split over $U(\bb^-)$, whose subquotients are rank 1 free modules over $R_n$. Before stating the result, let us investigate rank $1$ free modules, which are easy to classify.

\begin{defn} Let $p = \sum a_i x_i$ be a linear polynomial in $R_n$, and let $\smsum(p) \in \Z$ be an integer whose image in $\Bbbk$ agrees with $\sum a_i$. There is a free rank one graded $R_n$-module $R_n\ang{p}$ with
generator $1_p$ living in degree $\smsum(p)$. We define an $(R_n,\sltwo)$-module structure on $R_n\ang{p}$ by setting \begin{equation} \label{actiononshift} \dup(g \cdot 1_p) =
\dup(g) \cdot 1_p + gp \cdot 1_p, \qquad \ddown(g \cdot 1_p) = \ddown(g) \cdot 1_p \end{equation} for any (homogeneous) $g \in R$. \end{defn}

The formulas \eqref{actiononshift} are determined by the Leibniz rule from the action on the generator $1_p$:
\begin{equation} \label{actionon1pintro} \dup(1_p) = p \cdot 1_p, \qquad \ddown(1_p) = 0. \end{equation}
Note that $\smsum(p) = \ddown_R(p)$, though $\smsum$ is more descriptive notation. In Proposition \ref{prop:rank1} we prove that every $(R_n,\sltwo)$-module which is free of rank one as a graded $R_n$-module is isomorphic to $R_n\ang{p}$ for some $p$. If $p \ne p'$ then $R_n\ang{p}$ and $R_n\ang{p'}$ are non-isomorphic.

%

\begin{defn} Let $M$ be an $(R_n,\sltwo)$-module which is free and finitely generated as a graded $R_n$-module. Thus there is a finite set $I$ such that $M = \bigoplus_{i \in I}
M_i$ as graded $R_n$-modules, and each $M_i$ is free of rank $1$ over $R_n$. A \emph{downfree filtration} on $M$ is a splitting into free rank one $R_n$-modules as above, where
\begin{enumerate} \item Each $M_i$ is preserved by $\ddown$. \item The indexing set $I$ can be equipped with a partial order, so that $\dup(M_i) \subset \bigoplus_{j \le i} M_j$
for all $i \in I$. \end{enumerate} A homogeneous basis of $M$ as an $R_n$-module is called \emph{downfree} if it induces a downfree
filtration. \end{defn}

\begin{defn} Let $M$ be an $(R_n, \sltwo)$-module equipped with a downfree filtration. In particular, $\bigoplus_{j \le i} M_j$ is an $(R_n, \sltwo)$-submodule, and $M$ is equipped
with an $I$-indexed filtration by $(R_n, \sltwo)$-submodules, where the subquotients are free of rank 1 over $R_n$. Each subquotient must be isomorphic as an $(R_n, \sltwo)$-module to $R_n\ang{p_i}$ for some
unique $p_i \in R_n$. The multiset of linear polynomials $\{p_i\}_{i=1}^n$ will be called the \emph{downfree character} of $M$, with respect to this filtration. \end{defn}

The second main theorem of this paper says that well-known bases of morphism spaces in the categories of interest are actually downfree, and computes their downfree characters.

\begin{example} The nilHecke algebra $\NH_2$ on two strands is free of rank $(1+q^{-2})$ as a left $R_2$-module, spanned by the identity and the crossing. Since $\id$ is killed by $\dup$ and
$\ddown$, it generates an $\sltwo$-submodule $R \cdot \id \cong R\ang{0}$. Now
\begin{equation} \label{dXaltintro} \dup \left( \ig{1}{Xii} \right) \quad = \quad \ig{1}{line} \ig{1}{line} \; - 2 \; \ig{1}{XiiNW}. \end{equation}
In the quotient by $R \cdot \id$, $\dup$ will send the crossing $X$ to $-2 x_1 X$. So the basis $\{\id, X\}$ is downfree, and the downfree character of $\NH_2$ is $\{0, -2x_1\}$. Note that
$\{\smsum(0), \smsum(-2x_1)\} = \{0,-2\}$ which matches the degrees of this graded basis. If instead we had chosen the raising operator $\bardup$,  we would have gotten downfree character $\{0, -2x_2\}$.

Note that $\NH_2$ is also free as a right $R_2$ module, with the same basis. Because
\begin{equation} \label{dXaltaltintro} \dup \left( \ig{1}{Xii} \right) \quad = \quad - \; \ig{1}{line} \ig{1}{line} \; - 2 \; \ig{1}{XiiSE}, \end{equation}
the basis is downfree with character $\{0, -2x_2\}$. The right module character for $\dup$ matches the left module character for $\bardup$ because they are related by duality,
which also swaps the left and right action of $R_n$. \end{example}

\begin{rem} For the Hecke category, morphism spaces will also be $R_n$-bimodules, but the duality functor will not interchange these actions. A given basis may induce four
different characters, based on whether one selects the left or right action of $R_n$, and whether one chooses $\dup$ or $\bardup$. \end{rem}

More generally, associated to any element $w \in S_n$ and any reduced expression of $w$, one can construct the corresponding diagram in $\NH_n$ built from crossings, and this
element in $\NH_n$ is independent of the choice of reduced expression. We denote it $\psi_w \in \NH_n$. Then the elements $\{\psi_w\}_{w \in S_n}$ form a basis of $\NH_n$ as a left
or right $R_n$-module, which we call the \emph{nilCoxeter basis}. It is not hard to compute that
\begin{equation} \dup(R_n \cdot \psi_w) \subset \bigoplus_{y \le w} R_n \cdot \psi_y, \qquad \ddown(R_n \cdot \psi_w) \subset R_n \cdot \psi_w. \end{equation}

\begin{thm} \label{thmintro:HAPPYNHn} The nilCoxeter basis of $\NH_n$ over $R_n$ is a downfree basis, with partial order given by the Bruhat order on $S_n$. \end{thm}

This is proven in Theorem \ref{thm:KLRassgr}, which also contains an explicit formula for the downfree character, and the generalization to all simply-laced KLR algebras.

Similarly we can study the Hecke category $\HC(S_n)$. Now $R_n$ is the endomorphism ring of the monoidal identity, so all morphism spaces are naturally $R_n$-bimodules. Following ideas of Libedinsky \cite{LibLL}, Elias-Williamson in \cite{EWGr4sb} define the \emph{double leaves basis}, a basis of morphism spaces as left $R_n$-modules, which is indexed by coterminal Bruhat strolls. Let us summarize Theorem \ref{thm:Heckeassgr}, which contains an explicit formula for the downfree character.

\begin{thm} \label{thmintro:HAPPYHC} The double leaves basis of morphism spaces in $\HC(S_n)$ over $R_n$ is downfree, with partial order given by the lexicoBruhat order on
coterminal Bruhat strolls. \end{thm}

The natural bases of morphism spaces of $\UC(\sltwo)$ and $\dot{\UC}(\sltwo)$ are also downfree (conjecturally, since we do not prove it here), though over a different base ring.
Any given (nonzero) morphism space in $\UC(\sltwo)$ has $2n$ points on the boundary ($n$ oriented in and $n$ oriented out). This morphism space is free of rank $n!$
over $R_n \ot \Lambda$, where $\Lambda$ is the ring of symmetric functions acting by bubbles, and $R_n$ acts by dots on the inward-oriented strands. We place an $\sltwo$ structure
on $\Lambda$ in \S\ref{subsec-symfunc}; the lowering operator $\ddown$ depends on the choice of ambient weight. The combinatorics involved in describing the downfree character have
not been developed.

\subsection{What kinds of $\sltwo$-modules appear?} \label{subsec:sl2modules}

Having just described the kinds of $(R_n, \sltwo)$-modules which appear in practice, we can ask about what these modules look like as $\sltwo$-representations, with an eye towards
understanding their cores. 

In this paper, all $\sltwo$ representations have weights which are bounded below rather than above, so verma modules $\Delta(k)$ are defined by inducing from $U(\bb^-)$ rather than
$U(\bb^+)$. Let $L_k$ denote the irreducible weight representation of \emph{lowest weight} $k$ for all $k$, a quotient of $\Delta(k)$ and submodule of $\nabla(k)$. Note that $L_k$
is finite-dimensional if and only if $k \le 0$.

\begin{example} Consider $\Z[x]$ as an $\sltwo$-module. With its usual basis of monomials, the module looks like this.
\begin{equation}
	\xymatrix{
	\underset{1}{\bullet} \ar@/^/[rr]^{\dup=0} && \underset{x}{\bullet} \ar@/^/[ll]^{\ddown=1} \ar@/^/[rr]^{\dup=1}  && \cdots \ar@/^/[ll]^{\ddown=2}  \ar@/^/[rr]^{\dup=m-2} && \underset{x^{m-1}}{\bullet} \ar@/^/[ll]^{\ddown=m-1}  \ar@/^/[rr]^{\dup=m-1} && \underset{x^m}{\bullet} \ar@/^/[ll]^{\ddown=m}  \ar@/^/[rr]^{\dup=m} && \cdots \ar@/^/[ll]^{\ddown=m+1} }
\end{equation}
Hence $\Z[x] \cong \nabla(0)$. This is a
coverma module, with the trivial module (spanned by the identity element) as a submodule. The quotient by the trivial submodule is isomorphic to $\Delta(2)$. Note that $\Delta(2)
\not\cong \nabla(2)$ when we work over $\Z$ or in finite characteristic. \end{example}

As both a ring and as an $\sltwo$-module, we have $\Z[x_1, \ldots, x_n] \cong \Z[x]^{\otimes n}$. The following proposition is very easy to prove.

\begin{prop} (See Proposition \ref{prop:rank1}) For $p = \sum a_i x_i \in R_n$, there is an isomorphism of $\sltwo$ modules
\begin{equation} R_n\ang{p} \cong \nabla(a_1) \ot \cdots \ot \nabla(a_n). \end{equation}
In particular, $\core(R_n\ang{p})$ is nonzero if and only if $a_i \in \Z_{\le 0}$ for all $i$, in which case
\begin{equation} \core(R_n\ang{p}) \cong L_{a_1} \ot \cdots \ot L_{a_n}. \end{equation} \end{prop}

\begin{example} Suppose $M = R_n\ang{p}$ is generated in degree $-2$. If $p = -2 x_1$ then $\core(M)$ is three dimensional, if $p = -x_1 - x_2$ then $\core(M)$ is four dimensional,
and if $p = -3x_1 + x_2$ then $\core(M) = 0$. This illustrates why the character of an $(R_n, \sltwo)$-module is more useful than the graded degree. \end{example}

Suppose one has an $(R_n, \sltwo)$-module with a downfree filtration, and one knows the downfree character. By the proposition above, one knows the core of the associated graded
module. A priori, this does not make it any easier to determine the core of $M$, because a finite-dimensional submodule of a subquotient of $M$ need not lead to a finite-dimensional
submodule of $M$ itself. However, extensions between $(R_n, \sltwo)$-modules are even more limited than extensions between their underlying $\sltwo$-modules, and sometimes the
downfree character of $M$ will determine the core of $M$! Let us illustrate this with the following result.

\begin{thm} \label{thm:downfreeintrononsplit} Let $n=1$, so that $R_n = \Z[x]$. Let $M = \Z[x]\ang{ax} \oplus \Z[x]\ang{bx}$ be an $(\Z[x],\sltwo)$-module with a downfree filtration, where $\Z[x]\ang{bx}$ is the submodule. If the downfree filtration does not split then $b = a+2$. Unless $a=0$ and $b=2$, we have $\core(M) \cong \core(\Z[x]\ang{ax}) \oplus \core(\Z[x]\ang{bx})$. \end{thm}

\begin{proof} (Sketch) Let $1_a$ denote the generator of $\Z[x]\ang{ax}$, living in degree $a$. Then
\begin{equation} \label{dup1aintro} \ddown(1_a) = 0, \qquad \dup(1_a) = a x 1_a + m \end{equation} for some $m \in \Z[x]\ang{bx}$. For the $\sltwo$ relations to hold we need $\ddown(m) = 0$. But the kernel of $\ddown$ inside $\Z[x]\ang{b}$ is just the span of the generator $1_b$. Thus $m$ is a scalar multiple of $1_b$, and for degree reasons $b = a+2$. 
	
Even the associated graded has no core unless $a \le 0$, so assume $a \le 0$. There are no extensions between $\nabla(a)$ and $\nabla(a+2)$ unless $a = 0$, by the usual theory of
central characters. \end{proof}

\begin{rem} The lack of extensions between $\nabla(a)$ and $\nabla(a+2)$ does not mean that the downfree filtration splits. The splitting as an $\sltwo$-module and the splitting as an $R_1$-module are not compatible. \end{rem}
	
\begin{rem} When $a = 0$ and $b=2$, one can find an extension of $\Z[x]\ang{0}$ by $\Z[x]\ang{2x}$ with zero core, whereas the core of the associated graded would be one-dimensional. \end{rem}

This theorem implies that, for a two-step downfree filtration in one variable, the difference between the core of the original module and its associated graded is at most one copy
of the trivial module. Similarly, one can prove that a three-step filtration can remove a copy of $L_0$ or $L_{-1}$, but not $L_{-k}$ for $k \ge 2$.

\begin{rem} \label{rem:rootsapain} The situation is more complex in more than one variable, because polynomials in the roots $(x_i - x_j)$ are killed by $\ddown$, and this allows for more extensions (such
polynomials times $1_b$ are valid choices for $m$ in \eqref{dup1aintro}). See \S\ref{subsec:rexmoves} for an example. Interestingly, many of these extensions do not admit integrally-defined divided powers! Keeping track of divided powers and integrality properties does seem to rigidify the possible extensions. \end{rem}

The $(R_n, \sltwo)$-modules admitting a downfree filtration form a reasonably nice category which we feel is important to study. We hope to provide a methodical study of $(R_n,
\sltwo)$-modules in future work, and provide only the basics in this paper.

\subsection{Conclusion} \label{subsec:conclusion}

In \S\ref{sec-basics} we provide some basic results and definitions related to $\sltwo$-categories and their modules, and the special case of polynomial rings. In the subsequent sections we examine the categories $\UC^+(\gg)$ (in simply laced type), $\HC$, $\UC(\sltwo)$, and $\dot{\UC}(\sltwo)$ in turn, constructing the derivation $\ddown$, establishing the $\sltwo$ action, and verifying the claims made in this introduction about downfree filtrations. In \S\ref{subsec-sympoly} and \S\ref{subsec-symfunc} we discuss the $\sltwo$ action on symmetric polynomials and symmetric functions, which may be of independent interest.

In \S\ref{sec:hecke}, in order to prove results about the downfree filtration on the Hecke category, we need to establish some basic properties of light leaves and double leaves. Aside from this, the proofs in this paper are all relatively straightforward computations.

We find this new $\sltwo$ structure to be extremely tantalizing. In a follow-up paper we will introduce an $\sltwo$-enrichment of the categories of interest, where multiplicity spaces
are naturally finite-dimensional $\sltwo$ representations. We also introduce a conjectural approach to defining a triangulated category from an $\sltwo$-category, analogous to the
derived categories associated to $p$-dg algebras in \cite{KhoHopf, QiHopf}, whose Grothendieck group is naturally a module over $\Z[q+q^{-1}]$, the Grothendieck ring of finite-dimensional $\sltwo$
representations. We will state a number of conjectures about cores and Jacobson radicals, which translate into the existence of nice bases for the triangulated Grothendieck group. If
true, these conjectures would produce a new kind of categorification and a new kind of canonical basis, where structure coefficients are naturally unimodal, being the graded dimensions
of $\sltwo$ representations. For example, the Grothendieck group of the diagrammatic Hecke category would be a variant of the Hecke algebra with base ring $\Z[q+q^{-1}]$. Our
conjectures would also imply several open conjectures about categorification at a root of unity, e.g. that the $p$-dg Grothendieck group of the diagrammatic Hecke category is the Hecke
algebra at a root of unity.

While some aspects of the material developed in the follow-up paper are quite technical and lengthy to explain, other aspects are accessible and directly computational. We continue our introduction in \S\ref{sec:blathering} by explaining the basics of the theory of $\sltwo$ decompositions. In particular, this discussion further motivates the study of the core, and of \emph{iterated cores} that one obtains by combining the algebra structure and the $\sltwo$ structure on an $\sltwo$-category. We discuss examples in the diagrammatic Hecke category in \S\ref{subsec:iteratedcoreexamples}, and state our new \emph{$\sltwo$-adapted Soergel categorification conjecture} in \S\ref{subsec:Soergelconjecture}.

\begin{rem} This final remark is for those readers familiar with $p$-dg algebras and categorification at a root of unity. By forgetting the lowering operator and restricting from
$U(\sltwo)$ to $U(\bb^+)$, one obtains a $p$-dg structure on these categories. The biggest problem in computing the $p$-dg Grothenieck group is to prove that any object has a fantastic
filtration whose subquotients are certain $p$-dg-indecomposable objects. Historically this has been done by computing idempotent decompositions explicitly, but this method becomes
intractable quickly.

Our $\sltwo$-adapted Soergel categorification conjecture essentially states that Bott-Samelson objects have $\sltwo$-fantastic filtrations, whose subquotients are ``indecomposable
objects'' tensored with a multiplicity space which is a finite-dimensional $\sltwo$ representation. Because a finite-dimensional $\sltwo$ representation has a filtration by
one-dimensional $U(\bb^+)$-modules, this will yield a filtration by $p$-dg indecomposable objects (with multiplicity one); an $\sltwo$-fantastic filtration will be a $p$-dg fantastic
filtration. However, $\sltwo$-fantastic filtrations are more restrictive and include more structure, which ironically makes them easier to find. Once you compute the highest weight
vector, you can apply $\ddown$ to find a basis for the rest of the representation; this was a tool which was not previously available. \end{rem}

\paragraph{Acknowledgments.} B.~E. was partially supported by NSF CAREER grant DMS-1553032 and NSF FRG
grant DMS-1800498. This paper was completed while B.~E. was visiting the Institute of Advanced Study, where he was supported by NSF grant DMS-1926686. Y.Q. was partially supported by the NSF grant DMS-1947532. The authors would like to thank Mikhail Khovanov for his interest, and his suggestion in \S\ref{subsec-symfunc}.

\section{Representations of polynomial $\sltwo$-algebras} \label{sec-basics}

Most of the results in this chapter are relatively straightforward, but because we are not aware of any literature on the topic, we provide some details.

\subsection{$\gg$-algebras and Leibniz exercises} \label{subsec-leibniz}

\begin{defn} Let $\Bbbk$ be a commutative domain, and let $\gg$ be a Lie algebra over $\Bbbk$. A \emph{$\gg$-algebra} is a $\Bbbk$-algebra $A$ equipped with an action of $\gg$ by derivations. We sometimes write $(A,\gg)$ for this structure. One can define a \emph{$\gg$-category} similarly, and $\gg$ will act on each Hom space. A \emph{monoidal $\gg$-category} is a $\gg$-category with the additional requirement for each $\dx \in \gg$ that
\begin{subequations}
\begin{equation} \dx(f \ot g) = \dx(f) \ot g + f \ot \dx(g). \end{equation}
By the interchange law, it is equivalent to require that
\begin{equation} \dx(f \ot 1) = \dx(f) \ot 1, \quad \dx(1 \ot f) = 1 \ot \dx(f). \end{equation}
\end{subequations}
\end{defn}

By default in this paper $\Bbbk = \Z$. Before continuing, let us address the practical question of what it takes to place a $\gg$-algebra structure on $A$. Here are two simplifying
lemmas. The first reduces the data required to define a single derivation.

\begin{lem} \label{lem:defderivongen} Let $A$ be a $\Bbbk$-algebra given by generators and relations. To define a derivation $\dx$ on $A$, it suffices to specify $\dx(a)$ for each generator of $A$, and to check that
$\dx$ preserves all the relations of $A$ (what this means precisely is stated in the proof). \end{lem}

\begin{proof} Let $S$ be the generating set of $A$, and $F\ang{S}$ the free algebra on these generators. Any assignment $\dx \co S \to F\ang{S}$ will extend to a unique derivation on $F\ang{S}$, using the Leibniz rule to define the action of $\dx$ on a word in $S$. If $I$ is the ideal in $F\ang{S}$ generated by the relations, and $\dx$ sends the generating relations to elements of $I$, then $\dx$ preserves $I$ by the Leibniz rule. Thus $\dx$ descends to a derivation on $A$. \end{proof} 

The second lemma reduces the work required to check that a collection of derivations gives an action of a particular Lie algebra $\gg$.

\begin{lem} \label{lem:checksl2ongen} Let $\dx, \dy, \dz$ be three derivations on an algebra $A$, and suppose that 
\begin{equation}\label{commutator} [\dx, \dy](a) = \dz(a)\end{equation}
holds for a generating set of elements $a \in A$. Then it holds for all $a \in A$. \end{lem}

\begin{proof} The equation \eqref{commutator} is clearly linear in $a$, so it remains to check that if \eqref{commutator} holds for $a$ and $b$, then it holds for the product $ab$. We compute
\begin{align} [\dx, \dy](ab) &= \dx(\dy(a) b + a \dy(b)) - \dy(\dx(a) b + a \dx(b)) \\ \nonumber &= \dx \dy(a) b + \dy(a) \dx(b) + \dx(a) \dy(b) + a \dx \dy(b) - \dy \dx(a) b - \dx(a) \dy(b) - \dy(a) \dx(b) - a \dy \dx(b) \\ \nonumber &= [\dx, \dy](a) b + a [\dx, \dy](b) = \dz(a) b + a \dz(b) = \dz(ab). \end{align}
\end{proof}

Here is one more practical consideration.

\begin{defn} Let $\dx$ be a derivation on a $\Z$-algebra $A$. The \emph{divided powers} of $\dx$ are the operators
\begin{equation} \dx^{(k)} := \frac{\dx^k}{k!}, \end{equation}
which a priori send $A$ to $A \ot_\Z \Q$.  We say that the \emph{divided powers} $\dx^{(k)}$ are \emph{defined over $\Bbbk$} or just \emph{defined integrally}\footnote{We use the words ``defined integrally'' in general, but for most of our applications in this paper $\Bbbk = \Z$.} if there are operators $\dx^{(k)}$ defined in $A$ (without any base change) such that
\begin{equation} k! \cdot \dx^{(k)} = \dx^k, \qquad x^{(k)} x^{(\ell)} = \binom{k+\ell}{k} x^{(k+\ell)}.\end{equation} \end{defn}

\begin{lem} \label{lem:checkdivpowerongen} Let $\dx$ be a derivation on a $\Z$-algebra $A$, and suppose that $\dx^{(k)}(a)$ is defined integrally for a generating set of elements in $A$. Then $\dx^{(k)}$ is defined integrally (on all of $A$). \end{lem}

\begin{proof} Again, $\dx^{(k)}$ is a linear operator. It is an easy exercise in the Leibniz rule (and the binomial theorem) that
\begin{equation} \label{eq:dividedpowerleibniz} \dx^{(k)}(ab) = \sum_{i + j = k} \dx^{(i)}(a) \dx^{(j)}(b). \end{equation}
In particular, if $\dx^{(k)}$ is defined on $a$ and $b$ for all $k$, then it is defined on $ab$. \end{proof}

We write $\ot$ for $\ot_{\Bbbk}$.

\begin{cor} \label{cor:tensorproduct} Let $A$ and $B$ be $\gg$-algebras. Then the tensor product algebra $A \ot B$, with its tensor product $\gg$-action, is a $\gg$-algebra. If the
divided powers of $\dx \in \gg$ are defined integrally on $A$ and on $B$, then they are defined integrally on $A \ot B$. \end{cor}

\begin{proof} Since $A \ot B$ is generated by elements of $A$ and elements of $B$, Lemmas \ref{lem:checksl2ongen} and \ref{lem:checkdivpowerongen} imply the result. \end{proof}

\begin{rem}\label{rem:Hmodulealgebra}
The $\gg$-action on an algebra $A$ can be more conceptually understood in the language of a module-algebra over a Hopf algebra $H$ (see, e.g., \cite{Montgomerybook} for more details). Let $H$ be a Hopf algebra over $\Bbbk$. We will always assume that $H$ is a projective $\Bbbk$-module.  The comultiplication on $H$ is denoted $\Delta$ in Sweedler's notation:
\begin{equation}
\Delta (h)= \sum h_1\otimes h_2,
\end{equation}
for any $h\in H$, and the counit map is denoted $\epsilon: H \to \Bbbk$.  A $\Bbbk$-algebra $A$ is called an \emph{$H$-module algebra} if $A$ is an $H$-module, and the multiplication map of $A$ is compatible with the $H$-action: for any $a,b\in A$ and $h\in H$
\begin{equation}
h\cdot(ab)=\sum (h_1\cdot a)(h_2\cdot b), \qquad  h\cdot 1_A = \epsilon(h) 1_A,
\end{equation}
where $1_A$ is the identity element of $A$. When $H$ and $A$ are graded (super)algebras, the notions should be adapted so that the $H$ action on $A$ respects the graded (super)algebra structures.

As a particular case, for a $\gg$-algebra $A$, one may take $H$ to be the universal enveloping algebra of a Lie algebra $\gg$ over $\Bbbk$, with the $H$-action on $A$ induced by derivations of $\gg$ on $A$. Similarly, taking $H=\Bbbk[d]/(d^2)$ to be the graded superalgebra of dual numbers, where the degree of $d$ is set to be $1$ and 
\begin{equation}
\Delta(d)=d\otimes 1+1\otimes d, \quad \quad \epsilon(d)=0,
\end{equation}
one recovers the usual notion of a \emph{differential graded algebra} as an $H$-module algebra.
\end{rem}

\subsection{Modules over $\gg$-algebras}

\begin{defn} Let $(A,\gg)$ be a $\gg$-algebra. An \emph{$(A,\gg)$-module} is an $A$-module $M$ which is equipped with a $\gg$ action. For $\dx \in \gg$ write $\dx_M$ for the action of $\dx$ on $M$. We require the compatibility
\begin{equation} \label{eq:compat} \dx_M(a \cdot m) = \dx_A(a) \cdot m + a \cdot \dx_M(m) \end{equation} for any $a \in A$ and $m \in M$.
A morphism of $(A,\gg)$-modules is an $A$-module morphism $M \to N$ which intertwines the action of $\gg$. These form a space denoted $\Hom_{(A,\gg)}(M,N)$.
\end{defn}

This definition is completely analogous to the definition of a dg-module over a dg-algebra. Instead of just one differential, we keep track of a $\gg$ action. Just as for dg-modules, we can consider the internal Hom space.

\begin{defn} Let $M$ and $N$ be two $(A,\gg)$-modules. Then the space of $A$-module maps $\Hom_A(M,N)$ can be equipped with the structure of a $\gg$-module, where 
\begin{equation} \dx_{\Hom}(\phi)(m) := \dx_N(\phi(m)) - \phi(\dx_M(m)). \end{equation}
\end{defn}

\begin{prop} If $A$ is commutative, then $\Hom_A(M,N)$ is also an $A$-module. The actions of $A$ and $\gg$ are compatible, making $\Hom_A(M,N)$ into an $(A,\gg)$-module. \end{prop}

\begin{proof} We only need to check the Leibniz rule. We compute
\begin{align}
 \dx_{\Hom}(a \cdot \phi)(m) & = \dx_N(a \phi(m)) - a \phi(\dx_M(m)) = \dx_A(a) \phi(m) + a \dx_N(\phi(m)) - a \phi(\dx_M(m)) \nonumber \\
 & = (\dx_A(a) \cdot \phi)(m) + (a \cdot \dx_{\Hom}(\phi))(m).
 \end{align}
 The result follows.
\end{proof}

\begin{lem} Let $M$ be an $A$-module with a presentation. To give $M$ the structure of an $(A,\gg)$-module, it suffices to define $\dx_M(m)$ for all generators $m \in M$ and all
$\dx \in \gg$, and to check the relations. For a given $\dx \in \gg$, the divided powers of $\dx_M$ are defined integrally if and only if they are defined integrally on the
generators of $M$. \end{lem}

\begin{proof} This is entirely analogous to Lemmas \ref{lem:defderivongen}, \ref{lem:checksl2ongen}, and \ref{lem:checkdivpowerongen}. We leave the proof to the reader. \end{proof}

\begin{rem}\label{rem:modoverggalgebraasHmod}
Let $A$ be an $H$-module algebra as in Remark \ref{rem:Hmodulealgebra}, one may form the \emph{smash product ring} $A\# H$, which, as an algebra, is isomorphic to the tensor product algebra $A\otimes_{\Bbbk} H$, whose multiplication is given by
\begin{equation}
(a\otimes h) \cdot (b\otimes k)=\sum a (h_1\cdot b)\otimes h_2k,
\end{equation}
for any $a,b\in A$ and $h,k\in H$. The abelian category of $A\# H$-modules is acted upon by the abelian monoidal category of $\Bbbk$-projective $H$-representations. Let us denote this action by $\boxtimes$. Given an $A\# H$-module $M$ and a $\Bbbk$-projective $H$-representation $V$, $M\boxtimes V$ is isomorphic to $M\otimes_{\Bbbk} V$, with any element $a\otimes h\in A\# H$ acting on $x\otimes v\in M\boxtimes V$ by
\begin{equation}
(a\otimes h)\cdot (x\otimes v):=\sum \big((a\otimes h_1)\cdot x\big)\boxtimes h_2 v.
\end{equation}
By construction, $(\mbox{-})\boxtimes V$ is an exact functor on $A\# H$-modules.

In this language, an $(A,\gg)$-module is no other than a module over $A\# H$, where $H$ is the universal enveloping algebra of $\gg$ over $\Bbbk$. It follows readily that the category of $(A,\gg)$-modules constitute an abelian category, with a categorical module structure over the monoidal category of projective $H$-representations.
\end{rem}

\subsection{Weights}

When $\gg = \sltwo$ we change notation to impose one additional condition: that $\dh$ acts semisimply.

\begin{defn} A \emph{(weight) $\sltwo$-algebra (over $\Bbbk$)} is a $\Bbbk$-algebra $A$ equipped with an action of $\sltwo$ by derivations, on which the Cartan element $\dh \in
\sltwo$ acts diagonalizably (with eigenvalues in $\Bbbk$). A \emph{(weight) $(A,\sltwo)$-module} is an $A$-module $M$ with a compatible $\sltwo$ action, on which the Cartan element
$\dh \in \sltwo$ acts diagonalizably. A \emph{divided powers $\sltwo$-algebra} is an $\sltwo$-algebra where the divided powers of $\dup$ and $\ddown$ are
integrally defined, and similarly for a divided powers $(A,\sltwo)$-module. \end{defn}

We let $(\dup, \dh, -\ddown)$ be the standard $\sltwo$ triple, and we write $\dx_A$ or $\dx_M$ for the action of $\dx$ on $A$ or $M$, where $\dx \in \{\dup, \dh, \ddown\}$.

\begin{lem} Equip $A$ (resp. $M$) with a $\Bbbk$-grading by the eigenvalues of $\dh$, the usual weight grading. Then $A$ is a graded algebra (resp. and $M$ is a graded module).
\end{lem}

\begin{proof} For homogeneous elements $a_1$ and $m_2$ of weights $d_1$ and $d_2$ respectively, the Leibniz rule for $\dh$ implies that $a_1 m_2$ has weight $d_1 + d_2$. \end{proof}

Most often in this paper the eigenvalues of $\dh$ will be integers, and our rings and modules will be $\Z$-graded. Only in this chapter on generalities, and to some extent in
\S\ref{subsec-symfunc}, will we care about more general weights.

\begin{rem} For all the examples we study in this paper, the eigenvalues live in the image of the map $\Z \to \Bbbk$. Even when this map is not injective, all our examples can be 
compatibly $\Z$-graded. In this context a few statements need to be modified in the obvious way. For example, for an element $a \in \Bbbk$ we will often write ``$a \in \Z$'' to indicate that $a$ is in the image of $\Z \to \Bbbk$, but in the $\Z$-graded context one should instead choose a preimage $a \in \Z$. Without further ado, we assume that our examples are $\Z$-graded even in finite characteristic. \end{rem}

Let us make one remark about divided powers.

\begin{defn} Let $U = U_{\Z}(\sltwo)$ denote the idempotented divided powers form of the enveloping algebra of $\sltwo$. A definition can be found in, for instance, \cite[Section 3.1]{KLMS} (by specializing the $q$ to be $1$ in the quantum setting). All modules over $U$
will be assumed to be weight modules, so that the idempotent $1_n$ acts by multiplication by $n$ for $n \in \Z$. \end{defn}

Note that $U$ acts on any $\sltwo$-module where divided powers exist and the weights are integers. An action of $U$ on $A$ and $B$ extends to an action on $A \ot B$
(Corollary \ref{cor:tensorproduct}) since $U$ is a Hopf algebra. The proof that $U$ is a Hopf algebra is essentially the same as the proof of Lemma
\ref{lem:checkdivpowerongen}. A divided powers $\sltwo$-algebra over $\Z$ is the same thing as an algebra in the category of (weight) $U$-modules.

\subsection{On representations of $\sltwo$ over the integers} \label{subsec:integers}

The reader may be familiar with the properties of category $\OC$, but many things are different and slightly unfamiliar when working over the integers. The goal of this section is
to make precise what we mean by Verma and coVerma modules, and to warn the reader of some pitfalls. Because we study $\ddown$ instead of $-\ddown$, and bounded-below modules rather than bounded-above modules, some signs may differ from expectations.

\begin{defn} For $k \in \Bbbk$ we will define modules $\Delta(k)$ and $\nabla(k)$ over $\sltwo$ as follows. 
\begin{subequations} \label{vermacovermapictures} 
\begin{equation}
\Delta(k): \quad
\xymatrix{
\underset{0}{\bullet} \ar@/^/[rr]^{\dup=1} && \underset{1}{\bullet} \ar@/^/[ll]^{\ddown=k} \ar@/^/[rr]^{\dup=2}  && \cdots \ar@/^/[ll]^{\ddown=k+1}  \ar@/^/[rr]^{\dup=m-1} && \underset{m-1}{\bullet} \ar@/^/[ll]^{\ddown=k+m-2}  \ar@/^/[rr]^{\dup=m} && \underset{m}{\bullet} \ar@/^/[ll]^{\ddown=k+m-1}  \ar@/^/[rr]^{\dup=m+1} && \cdots \ar@/^/[ll]^{\ddown=k+m}  
}
\end{equation}
\begin{equation}
\nabla(k): \quad
\xymatrix{
\underset{0}{\bullet} \ar@/^/[rr]^{\dup=k} && \underset{1}{\bullet} \ar@/^/[ll]^{\ddown=1} \ar@/^/[rr]^{\dup=k+1}  && \cdots \ar@/^/[ll]^{\ddown=2}  \ar@/^/[rr]^{\dup=k+m-2} && \underset{m-1}{\bullet} \ar@/^/[ll]^{\ddown=m-1}  \ar@/^/[rr]^{\dup=k+m-1} && \underset{m}{\bullet} \ar@/^/[ll]^{\ddown=m}  \ar@/^/[rr]^{\dup=k+m} && \cdots \ar@/^/[ll]^{\ddown=m+1}  
}
\end{equation}
\end{subequations}
 Let us rephrase these pictures in formulas. The module $\Delta(k)$ has a free $\Z$-basis $\{v_{k,m}\}_{m \ge 0}$, corresponding to the dots labeled by $m$ in the picture above, where $v_{k,m}$ has weight $k + 2m$. Set $v_{k,m}:= 0$ for $m < 0$. We have
\begin{equation} \label{vermadef} 
\dup(v_{k,m}) = (m+1) v_{k,m+1}, \qquad \ddown(v_{k,m+1}) = (m+k) v_{k,m}. 
\end{equation}
Similarly, $\nabla(k)$ has a free $\Z$-basis $\{w_{k,m}\}$, where $w_{k,m}$ has weight $k + 2m$. Set $w_{k,m} := 0$ for $m < 0$. We have
\begin{equation} \label{covermadef} \dup(w_{k,m}) = (m+k) w_{k,m+1}, \qquad \ddown(w_{k,m+1}) = (m+1) w_{k,m}. \end{equation}
\end{defn}	

\begin{lem} When $a \in \Z$, divided powers are defined integrally on $\Delta(a)$ and $\nabla(a)$. More generally, divided powers are defined in $\Bbbk$ whenever the \emph{$a$-binomial coefficients} are defined in $\Bbbk$ for all $l \in \Z$ and $m \in \Z_{\ge 0}$: these are the elements
\begin{equation} \binom{a + l}{m} := \frac{(a+l)(a+l-1) \cdots (a+l+1-m)}{m!}. \end{equation}  \end{lem}

\begin{proof} It is a simple computation that
\begin{subequations}\label{actiononcoverma}
\begin{equation} \dup^{(l)} v_{k,m} = \binom{m+l}{l} v_{k,m+l}, \qquad \ddown^{(l)} w_{k,m} = \binom{m}{l} w_{k,m-l}, \end{equation}
\begin{equation} \ddown^{(l)} v_{k,m} = \binom{m+k-1}{l} v_{k,m-l}, \qquad \dup^{(l)} w_{k,m} = \binom{m+k + l - 1}{l} w_{k,m+l}. \end{equation}
\end{subequations}
\end{proof} 

For example, suppose $\Bbbk = \Z[y]$ and $a = y$. Then \[ \binom{a+7}{3} = \frac{(y+7)(y+6)(y+5)}{3!}, \] is not an element of $\Bbbk$, so divided powers are not defined integrally on $\Delta(a)$ or $\nabla(a)$.

For the rest of this section, $\Bbbk = \Z$, so we can think of Verma and coVerma modules as $U$-modules.

Let $v_k := v_{k,0}$ and $w_k := w_{k,0}$. Note that \begin{equation} \dup^{(m)}(v_k) = v_{k,m}, \qquad \ddown^{(m)}(w_{k,m}) = w_k. \end{equation} The module $\Delta(k)$ is
generated by $v_k$ over $U$, and the divided powers $\dup^{(m)}$ applied to this generator give the basis. Meanwhile, $\nabla(k)$ is infinitely-generated over $U$, but it is
``co-generated'' by $w_k$, as the divided powers $\ddown^{(m)}$ bring every basis element to $w_k$.

Given any $U$-module $M$ which is free as a $\Z$-module, $2M \subset M$ will be a proper submodule. Simplicity, in its naive sense as when working over a field, is not as useful a concept.

\begin{defn} A morphism of $U$-modules is called \emph{$\dh$-split} if it is a split morphism of $\Z$-modules for each weight space. Let $U\spmod$ denote the category whose objects
are $U$-modules where weight spaces are free over $\Z$, and whose morphisms are $\dh$-split maps. If $M$ is a module in $U\spmod$, a $U$-submodule is called \emph{$\dh$-split} if
the inclusion map is $\dh$-split, and similarly for quotient modules. \end{defn}

Note that any isomorphism, or more generally any genuinely split map (over $U$), is automatically $\dh$-split. However, $U\spmod$ is not an additive category, as the sum of
$\dh$-split morphisms need not be $\dh$-split. After all, $\id + \id = 2 \id$. Nonetheless, if $f$ is any $\dh$-split morphism then its kernel and cokernel naturally live in
$U\spmod$, so $U\spmod$ shares some features with an abelian category.

When restricting to $\dh$-split morphisms, verma and coverma modules have features which resemble those in the familiar category $\OC$.

\begin{prop} The modules $\Delta(k)$ and $\nabla(k)$ are indecomposable for all $k \in \Z$. When $k > 0$, the modules $\Delta(k)$ and $\nabla(k)$ have no $\dh$-split submodules or quotients. There are $\dh$-split short exact sequences for all $k \le 0$ given by
\begin{subequations} \label{eq:ses}
\begin{equation} 0 \to \nabla(-k+2) \to \Delta(k) \to W(k) \to 0, \end{equation} 
\begin{equation} 0 \to W^\vee(k) \to \nabla(k) \to \Delta(-k+2) \to 0, \end{equation}
\end{subequations}
where $W(k)$ and $W^\vee(k)$, the so-called \emph{Weyl} and \emph{dual Weyl modules}, are defined by these short exact sequences. Aside from those given in \eqref{eq:ses}, there are no other $\dh$-split submodules or quotients of $\Delta(k)$ and $\nabla(k)$ for $k \le 0$. \end{prop}

\begin{proof} Since all weight spaces are free of rank $1$ over $\Z$, an $\dh$-split submodule is determined by which weight spaces it includes. Now the usual arguments (analyzing
which arrows have zero coefficient) imply that only certain collections of weight spaces can give a submodule. One need only verify that when $k \le 0$ the submodule of $\Delta(k)$ generated in degree $-k+2$ is $\nabla(-k+2)$, which one can do directly from the pictures \eqref{vermacovermapictures} or the corresponding formulas. \end{proof}

What is more interesting is the following observation.

\begin{prop} If $\Bbbk$ contains $\Q$ as a subring, then $\nabla(k) \cong \Delta(k)$ if and only if $k > 0$. Otherwise, $\nabla(k) \cong \Delta(k)$ if and only if $k=1$.
\end{prop}

Note a major difference between this setting and the usual category $\OC$ in characteristic zero. Normally, $\Delta(k) \cong \nabla(k)$ for all $k > 0$, and one reverses the
placement of $\Delta(-k+2)$ and $\nabla(-k+2)$ in \eqref{eq:ses}; they are isomorphic, so it doesn't matter, but it is often done to help illustrate the verma and coverma
resolutions of simple modules. Over $\Z$, one is not permitted to swap $\Delta(-k+2)$ and $\nabla(-k+2)$, and the Weyl and dual Weyl modules do not have Verma or coVerma resolutions. This fact will be of great importance in the sequel.

\begin{proof} Clearly $\nabla(1) \cong \Delta(1)$ by sending $w_{1,m} \mapsto v_{1,m}$, since the formulas \eqref{vermadef} and \eqref{covermadef} agree. When $k \le 0$ clearly
$\nabla(k) \not\cong \Delta(k)$ since one has a finite rank submodule and the other does not. So suppose $k > 0$. What happens over $\Q$ is well-known, so assume that there is some
$a \ge 1$ such that $\{1, \ldots, a\}$ is invertible in $\Bbbk$ and $a+1$ is not invertible. We claim that the image of $\ddown^{(a)}$ in degree $k$ is different when comparing
$\nabla(k)$ to $\Delta(k)$, which implies they are non-isomorphic. Clearly $w_k$ is in the image of $\ddown^{(a)}$ inside $\nabla(k)$, whereas only the span of $(a+1)v_k$ is in the
image of $\ddown^{(a)}$ inside $\Delta(k)$. \end{proof}
	
\begin{rem} If $M$ is any module with a lowest weight vector $v$ in weight $k$, then there is a morphism $\Delta(k) \to M$ sending $v_k \mapsto v$. This morphism will not
necessarily be $\dh$-split. For example, the natural morphism $\Delta(k) \to \nabla(k)$ is not $\dh$-split except when $k=1$. \end{rem}
	
\subsection{Rank-one modules over polynomial rings} \label{subsec:rank1}

Let $R_1 = \Bbbk[x]$, equipped with its standard $\sltwo$ structure from \eqref{actiononRn}, where $\dup = x^2 \frac{\pa}{\pa x}$ and $\ddown = \frac{\pa}{\pa x}$. By our choice of convention, $\dh=2x\frac{\pa }{\pa x}$.

\begin{prop} As an $\sltwo$-module, $R_1 \cong \nabla(0)$. In particular, divided powers are well-defined, so $R_1$ is a $U$-module. \end{prop}

\begin{proof} If we send $w_{0,m} \mapsto x^m$, we can confirm \eqref{covermadef} easily. \end{proof}

\begin{rem} One is tempted to say that the ideal $(x)$ is isomorphic to $\Delta(2)$, but one must be careful with this statement. Ideals are usually thought of as submodules, so
that $(x)$ is identified with a subset of $R_1$. This subset is an $R_1$-submodule but not a $U$-submodule, because it is not preserved by $\ddown$. This was discussed at more length in the introduction, see the toy example from \S\ref{subsec:NHintro}. \end{rem}


\begin{defn} Let $a \in \Bbbk$ and let $R_1\ang{a}$ denote the free rank one graded $R_1$ module with generator $1_a$. We give it an $(R_1,\sltwo)$-module structure by setting
\begin{equation} \label{actionon1a} \dup(1_a) = a x \cdot 1_a, \qquad \ddown(1_a) = 0, \qquad \dh(1_a) = a 1_a, \end{equation}
and extending these operators to all of $R_1\ang{a}$ by the Leibniz rule. \end{defn}

That $R_1\ang{a}$ is well-defined is a consequence of the following proposition.
	
\begin{prop} \label{prop:rank1}
\begin{enumerate}
\item[(i)] Any $(R_1,\sltwo)$-module structure on the rank-one free module $R_1$ is isomorphic to $R_1(a)$ for a unique $a\in \Bbbk$.
\item[(ii)] As an $\sltwo$-module, $R_1\ang{a} \cong \nabla(a)$. In particular, divided powers are well-defined when $a \in \Z$.
\end{enumerate}
\end{prop}

\begin{proof} 
For a rank-one module, denote by $v$ a generator as an $R_1$-module. Then, the module is also generated by $v$ as an $(R_1,\sltwo)$-module. For degree reasons, there must be $a,a^\prime\in \Bbbk$ such that
\[
\dup(v)=axv, \quad  \quad \dh(v)=a^\prime v ,\quad \quad \ddown(v)=0.
\]
The commutation relation $[\dup, -\ddown]=\dh$ applied to the generator $v$ shows that
\[
a^\prime v = \dh v= \ddown \dup v - \dup \ddown v = a v,
\]
implying that $a^\prime=a$. 

For the second statement, if we send $w_{a,m} \mapsto x^m \cdot 1_a$, we can confirm \eqref{covermadef} easily. 
\end{proof}

Now we let $R_n = \Bbbk[x_1, \ldots, x_n]$ be the polynomial $\sltwo$-algebra, as in Definition \ref{defn:polySLTDGA}.

\begin{lem} \label{lem:polyringiscool} As an algebra and an $\sltwo$-module we have $R_n \cong R_1 \ot \cdots \ot R_1$. In particular, divided powers are well-defined. \end{lem}

\begin{proof} The isomorphism with the tensor product is obvious. See Corollary \ref{cor:tensorproduct} for the rest. \end{proof}

We recall a definition from the introduction.

\begin{defn} Let $p = \sum a_i x_i$ be a linear polynomial in $R_n$, and let $\smsum(p) = \sum a_i \in \Bbbk$. Note that $\smsum(p) = \ddown_R(p)$, though $\smsum$ is more descriptive notation. There is a free rank one graded $R_n$-module $R_n\ang{p}$ with generator $1_p$. We define an $(R_n,\sltwo)$-module structure on $R_n\ang{p}$ by setting
\begin{equation} \label{actionon1p} \dup(1_p) = p \cdot 1_p, \qquad \ddown(1_p) = 0, \qquad \dh(1_p) = \smsum(p) 1_p, \end{equation}
and extending by the Leibniz rule.
\end{defn}

\begin{prop} \label{prop:Rnp} As an $\sltwo$-module, 
\begin{equation}  \label{Rnangpisom} R_n\ang{p} \cong R_1\ang{a_1} \ot \cdots \ot R_1\ang{a_n} \cong \nabla(a_1) \ot \cdots \ot \nabla(a_n). \end{equation}
In particular, divided powers are well-defined when $a_i \in \Z$ for all $i$. \end{prop}

\begin{proof} It is easily verified that the map $1_p \mapsto 1_{a_1} \ot \cdots \ot 1_{a_n}$ induces an isomorphism. \end{proof}

\begin{prop} Every $(R_n,\sltwo)$-module which is free of rank one as a graded $R_n$-module is isomorphic to $R_n\ang{p}$ for some $p$. \end{prop}

\begin{proof} Let $M$ be such a module, and name the generator $1_M$. Then $\ddown(1_M) = 0$, and $\dup(1_M) = p \cdot 1_M$ for some linear polynomial $p$, both for
degree reasons. The remaining structure follows from the Leibniz rule. The fact that the degree of $1_M$ must be $\smsum(p)$ follows from the fact that 
\begin{equation} \ddown \dup(1_M) = \ddown(p) \cdot 1_M = \smsum(p) 1_M. \end{equation} \end{proof}

\begin{rem} For more examples of $\sltwo$-algebras, see \S\ref{subsec-sympoly} and \S\ref{subsec-symfunc}. \end{rem}

\subsection{The core} \label{subsec:core}

We will eventually be interested in settings where the base ring $\Bbbk$ is $\Z$ or $\Q$ or $\Z[y]$ or of finite characteristic. Thus we are careful in this chapter to make general statements.

\begin{defn} An $\sltwo(\Bbbk)$ module is called \emph{bounded} if \begin{itemize}
	\item it is a weight module, with weights in $\Z$,
	\item the set of weights with non-zero weight spaces is bounded below,
	\item and each non-zero weight space is free of finite rank over $\Bbbk$. \end{itemize}
\end{defn}

\begin{prop} Suppose that $\Z \to \Bbbk$ is injective and $\Bbbk$ is Noetherian. Any bounded $\sltwo$ representation $M$ will have a maximal submodule which is finitely generated over $\Bbbk$, which we call
the \emph{core} of $M$ and denote $\core(M)$. It satisfies
\begin{equation} \label{usetodefine}
\core(M) = \{m \in M \mid \dup^N(m) = 0 \text{ for some } N \in \N  \}. \end{equation}
Note that the core may be zero.
 \end{prop}

\begin{proof} For the sake of this proof, define $\core(M)$ using \eqref{usetodefine}. Note that $\core(M)$ is a weight module. Clearly $\dup$ acts nilpotently on any submodule
which is finitely generated as a $\Bbbk$-module, so $\core(M)$ contains all such submodules. We need only prove that $\core(M)$ is finitely generated over $\Bbbk$.

We now argue that $\core(M)$ is locally finitely generated, i.e. any element is contained in an $\sltwo$ submodule which is finitely generated over $\Bbbk$. Suppose that $m \in M$
homogeneous is acted upon nilpotently by $\dup$. By the PBW theorem, the span of $\{\ddown^a \dh^b \dup^c \cdot m\}_{a,b,c \ge 0}$ is an $\sltwo$ subrepresentation containing $m$.
Only finitely many pairs $(a,c)$ will give a nonzero result, and the span of $\{\ddown^a \dh^b \dup^c \cdot m\}_{b \ge 0}$ agrees with the span of the single vector $\ddown^a \dup^c
\cdot m$. Thus $m$ is contained in a subrepresentation which is finite rank over $\Bbbk$.

If $\core(M)$ is not finitely generated over $\Bbbk$, then by the Noetherian hypothesis it must have nonzero elements in infinitely many weight spaces (lest it be contained in a
finite rank $\Bbbk$-module). Because $\core(M)$ is locally finitely generated, it must have highest weight vectors in arbitrarily high weights. We now argue that $\core(M)$ also has
(lowest weight) vectors in arbitrary low weights, a contradiction because $M$ is bounded. This last argument is standard, we are just being careful because we have placed very few
assumptions on $\Bbbk$.

Suppose that $\core(M)$ has a highest weight vector $v$ in weight $k > 0$. Let $\Delta'(k)$ be the bounded-above Verma module with highest weight $k$. Then $v$ induces a nonzero map
of $\sltwo(\Bbbk)$-modules $\phi \co \Delta'(k) \to \core(M)$. Since the target is a submodule of a free $\Bbbk$-module it is torsion free over $\Bbbk$. Since each weight space in
$\Delta'(k)$ is free of rank 1 over $\Bbbk$, the kernel of $\phi$ is divisible, so it is spanned by a subset of the weight spaces in $\Delta'(k)$. Now standard arguments imply that
the kernel of $\phi$ is zero in weight $-k$, so $\core(M)$ is nonzero in weight $-k$. \end{proof}

\begin{rem} When $\Bbbk$ is a PID, one can use Smith normal form to deduce that the kernel of a map between free (finite rank) modules is a split summand of the source. One can use
this to prove that $\core(M)$ is a $\dh$-split submodule of $M$. We imagine this is always true (even when $\Bbbk$ is a Dedekind domain or something unusual) but we do not have the
knowledge or examples to say one way or another. \end{rem}

\begin{lem} There is a left-exact functor $\core$ from bounded $\sltwo$-modules to $\sltwo$-modules which are finitely generated over $\Bbbk$, sending $M$ to $\core(M)$, and restricting any morphism to the core.
\end{lem}

\begin{proof} The image under an $\sltwo$-intertwiner of a finitely generated module is a finitely generated submodule. Thus the core is sent to the core under any $\sltwo$-intertwiner. Injective maps restrict to injective maps, so $\core$ is left-exact. \end{proof}

\begin{rem} \label{rem:coreofquotient}  Note that $\core$ is not right exact. For example, for $a \le 0$ the canonical map $\Delta(a) \to W(a)$ is surjective, but $\core(\Delta(a)) = 0$ while $\core(W(a)) = W(a)$. In this sense, the core of a quotient can be larger than the core of the original module. Of course, the core of a quotient can be smaller too, as the quotient map could kill part of the core. \end{rem}

\begin{prop} \label{prop:coreRnp} Let $p = \sum a_i x_i \in R_n$ with $a_i \in \Z$. If $a_i > 0$ for some $i$ then $\core(R_n\ang{p}) = 0$. If $a_i \le 0$ for all $i$ then
\begin{equation} \core(R_n\ang{p}) \cong W^\vee_{a_1} \ot \cdots \ot W^\vee_{a_n}. \end{equation}
\end{prop}

\begin{proof} We use the identification \eqref{Rnangpisom} of $R_n\ang{p}$ with a tensor product of covermas.
	
Let $M$ be any $\sltwo$ module and consider $\nabla(a) \ot M$. Any vector $v$ in this tensor product can be written uniquely as
\begin{equation} \sum_{j \ge 0} w_{a,j} \ot m_j \end{equation}
for some $m_j \in M$, with finitely many $m_j$ being nonzero. Let $J$ be the maximal value of $j$ for which $m_j \ne 0$. Then
\begin{equation} \dup^N(\sum_{j \ge 0} w_{a,j} \ot m_j) = \dup^N(w_{a,J}) \ot m_J + \sum_{j' < J+N} w_{a,j'} \ot n_{j'} \end{equation}
for some $n_{j'} \in M$. So long as $\dup^N(w_{a,J}) \ne 0$, $\dup^N(v) \ne 0$. But by the formula from \eqref{actiononcoverma}, $\dup^N(w_{a,J})$ is never zero when $a > 0$. In fact, it is never zero when $J > -a$.

As a consequence we deduce that $\nabla(a) \ot M$ has zero core when $a > 0$. Similar arguments prove that $\Delta(a) \ot M$ has no core for any $a \in \Z$. This also follows from
the familiar idea that tensoring a verma module with anything yields a module with a verma filtration, and modules with a verma filtration have no finite-dimensional submodules.

This proves that $\core(R_n\ang{p}) = 0$ if any $a_i > 0$. If $a_i \le 0$ for all $i$, then the tensor product of the submodules $W^\vee_{a_i}$ clearly lives inside $\core(R_n\ang{p})$. Since $\nabla(a_i) / W^\vee_{a_i} \cong \Delta(-a_i + 2)$, the quotient of $R_n\ang{p}$ by the tensor product of coweyl modules has a filtration whose subquotients are isomorphic to $M \ot \Delta(a)$ for some $a$. If $\dup$ has no nilpotents on the associated graded of a filtered module, it has no nilpotents on the whole module. This proves that the core of $R_n\ang{p}$ is not bigger than expected. \end{proof}

\begin{cor} \label{cor:corebasis} Let $p = \sum a_i x_i$. Then the elements
\begin{equation} \BM := \{ x_1^{b_1} x_2^{b_2} \cdots x_n^{b_n} \mid 0 \le b_i \le -a_i \} \end{equation}
form a basis for the core of $R_n\ang{p}$. \end{cor}

\begin{proof} This follows from the case $n=1$, where it is straightforward. \end{proof}

Finally, we mention one more useful result. The same argument will show that the core of a tensor product contains the tensor product of the cores.

\begin{prop} \label{prop:coresubalgebra} In an $\sltwo$-category where all morphism spaces are bounded $\sltwo$ representations, the core is a subcategory. In a monoidal $\sltwo$-category with this property, the core is a monoidal subcategory.
\end{prop}

\begin{proof} We need to prove that the core contains all identity maps (which is obvious) and is closed under composition (both horizontal and vertical). We use the description of
the core as those elements on which $\dup$ acts nilpotently. If $\dup^{N_1}(f) = 0$ and $\dup^{N_2}(g) = 0$ then $\dup^{N_1 + N_2}(f \circ g) = 0$ and $\dup^{N_1 + N_2}(f \ot g) =
0$, by the Leibniz rule. \end{proof}

\section{On decompositions in the presence of $\sltwo$} \label{sec:blathering}

In this quasi-introductory chapter we discuss the nature of direct sum decompositions in the presence of an $\sltwo$ action, and the possible implications for the Grothendieck group.
Afterwards, we will reconnect these ideas to the concept of the core. For simplicity, we work with categories over a field $\Bbbk$ rather than the integers.

The idea has already arisen in categorification at a root of unity (which studies categories with a derivation $\dup$) to study a direct sum decomposition not with the traditional use of idempotents, but by studying filtrations on representable functors instead. This is a guiding principle below.

\subsection{Actions on partial idempotent completions}

It is traditional when studying (pre)additive categories over a field to take (the additive closure and then) the Karoubi envelope. One obtains a Krull-Schmidt category, whose split
Grothendieck group is spanned by isomorphism classes of indecomposable objects. Then one uses direct sum decompositions to find relations in the Grothendieck group. One way to motivate the Karoubi envelope is that it naturally appears when studying modules over the category: representable modules $\Hom(X,-)$ are projective, and thus so are their direct summands $\Hom(X,-) \circ e$ for any idempotent $e \in \End(X)$.

However, for a category with an $\sltwo$ action, there is no induced action of $\sltwo$ on the Karoubi envelope!

Let $X$ be an object and $e \in \End(X)$ an idempotent. The subspace \[
\Hom(X,-) \circ e \subset \Hom(X,-) \] is not closed under $\sltwo$ in general. For each $\dx \in \{\dup, \ddown, \dh\}$ one can try to define an operator $\underline{\dx}$ on $\Hom(X,-) \circ
e$ by truncating the usual action, to wit \begin{equation} \label{eq:inducedaction} \underline{\dx}(f) := \dx(f) e \qquad \text{for } f \in \Hom(X,-) \circ e. \end{equation} Each operator $\underline{\dx}$ satisfies the Leibniz rule. However,
these operators will not necessarily satisfy the $\sltwo$ relations.

\begin{defn} An \emph{$\sltwo$ submodule idempotent} is an idempotent satisfying $e \dx(e) = 0$ for all $\dx \in \{\dup, \dh, \ddown\}$. An \emph{$\sltwo$ quotient idempotent} instead satisfies $\dx(e) e = 0$. An object which has no nontrivial $\sltwo$ submodule idempotents is called \emph{$\sltwo$-indecomposable}. \end{defn}

By applying the Leibniz rule to $\dx(e^2) = \dx(e)$, one deduces that $e$ is an $\sltwo$ submodule idempotent if and only if $\dx(e) = \dx(e) e$. Hence $\Hom(X,-) \circ e$ is closed
under $\dx$, and the formula \eqref{eq:inducedaction} is another way to write the restriction of the $\sltwo$ action to this submodule.

It is also easy to deduce that $e$ is an $\sltwo$ submodule idempotent if and only if $(1-e)$ is an $\sltwo$ quotient idempotent. In this case, $\Hom(X,-) \circ (1-e)$ is not closed under the $\sltwo$ action. However, using the identification
\[ \Hom(X,-) \circ (1-e) \cong \Hom(X,-) / \Hom(X,-) \circ e, \]
the $\sltwo$ action induced on this quotient agrees with the formulas from \eqref{eq:inducedaction}.

Let $e$ be an $\sltwo$ submodule idempotent. Looking at representable functors out of the partial idempotent completion which adjoins the images of $e$ and $1-e$, we have a short exact sequence of $\sltwo$ representations
\begin{equation} \label{eq:idempotentses} 0 \to \Hom(X,-) \circ e \to \Hom(X,-) \to \Hom(X,-) \circ (1-e) \to 0, \end{equation}
which splits as modules for the underlying additive category. In this way, the direct sum decomposition $X \cong \Im(e) \oplus \Im(1-e)$ in the additive category becomes a filtration for $\sltwo$. We still refer to the idempotents $e$ and $(1-e)$ as complements.

\begin{defn} An \emph{$\sltwo$ decomposition} of an object $X$ is an ordered sequence of idempotents $e_1, e_2, \ldots, e_n$ such that
\begin{enumerate} \item The set $\{e_1, \ldots, e_n\}$ is a complete collection of orthogonal idempotents: $e_i e_j = e_j e_i = 0$ for $i \ne j$, and $\sum_{i=1}^n e_i = \id_X$. 
\item For any $k \le n$, the idempotent $e_1 + e_2 + \ldots + e_k$ is an $\sltwo$ submodule idempotent. \end{enumerate}
In this case, the idempotents $e_i$ are called \emph{$\sltwo$ subquotient idempotents}. 
\end{defn}

Similarly, the formula \eqref{eq:inducedaction} defines an $\sltwo$ action on $\Hom(X,-) \circ e$ for any $\sltwo$ subquotient idempotent. For the analogous definition in the context of $p$-dg algebras, see \cite[\S 4]{EQpDGbig}.

\subsection{Future goals, and cofibrance}

The goal is to find a setting where certain idempotent decompositions give rise to relations on the Grothendieck group, much as they did in the split Grothendieck group of the
underlying additive category. Another goal is for this Grothendieck group to be naturally a module over the Grothendieck ring of finite-dimensional $\sltwo$ representations. Not every
idempotent decomposition should give rise to a relation. For a general idempotent, the subspaces $\Hom(X,-) \circ e$ and $\Hom(X,-) \circ (1-e)$ are entangled by the $\sltwo$ action,
and no filtration \eqref{eq:idempotentses} is present; one does not expect a relation in the Grothendieck group. However, even for submodule idempotents, one does not always expect a
relation, because of a technical point we now address.

The program of categorification at a root of unity, which accounts for the action of $\dup$ but not $\ddown$, has an established framework for Grothendieck groups like the one we seek.
There is a triangulated category associated to a category with an action of $\dup$, constructed using Hopfological algebra \cite{KhoHopf,QiHopf}. In this triangulated category, the
short exact sequence \eqref{eq:idempotentses} does not always give rise to a distinguished triangle, but it does when the module $\Hom(X,-) \circ e$ is \emph{cofibrant}. In particular,
the appropriate category to study is not the partial idempotent completion which adjoins all $\sltwo$ subquotient idempotents, but only those which are cofibrant.

There are
typically two ways to prove that the image of an idempotent is cofibrant. The first is to show that it is isomorphic to a representable module, which happens when the idempotent factors compatibly with $\dup$:
\begin{equation} \label{eq:factoringidempotent} e = i \circ p, \qquad p \circ i = \id_Y \text{ for some object } Y, \qquad \dup(p) \circ i = 0. \end{equation}
The second method, the \emph{two-out-of-three principle}, uses the fact that if two out of three terms in \eqref{eq:idempotentses} are cofibrant, then so is the third. Thus the complement of a representable module in a representable module is cofibrant. Adjoining its image, one can assume it is representable for future decompositions.

So to summarize, one expects relations on the Grothendieck group coming from cofibrant short exact sequences like \eqref{eq:idempotentses}, and cofibrance can be confirmed by factoring the idempotent as a projection and inclusion, and checking $\dup(p) \circ i = 0$. For a longer discussion of these thorny issues, see \cite[\S 4]{EQpDGbig}.

Everything below in this section should be taken with a grain of salt. The authors' expectation is that there is a triangulated category one can build from any category with
$\sltwo$ action, where short exact sequences coming from cofibrant $\sltwo$ submodule idempotents give rise to distinguished triangles. If one can factor an idempotent as $e = i
\circ p$, where $i$ and $p$ have degree zero, then it will be cofibrant if $\dx(p) \circ i = 0$ for all $\dx \in \sltwo$. However, there are additional cofibrant objects coming from
representable modules in the enriched category.

\begin{defn} Let $\CC$ be an $\sltwo$-category. The \emph{$\sltwo$-enriched category} is the $\sltwo$-category $\CC_{\enrich}$ defined as follows. Its objects are formal tensor products $B \sqot V$ (see Remark \ref{rem:modoverggalgebraasHmod} for an explanation of the notation), for an object $B$ of $\CC$ and a finite-dimensional representation $V$ of $\sltwo$. We set
\begin{equation} \Hom_{\CC_{\enrich}}(B \sqot V, B' \sqot V') := \Hom_{\CC}(B,B') \ot \Hom_{\sltwo}(V,V') \end{equation}
as an $\sltwo$-module. Composition is given by $(f \ot \phi) \circ (g \ot \psi) = (f \circ g) \ot (\phi \circ \psi)$. There is a monoidal action of the category of finite-dimensional $\sltwo$ representations on $\CC_{\enrich}$, where on objects one has $(B \sqot V) \ot V' := B \sqot (V \ot V')$. \end{defn}

One should think about $B \sqot V$ as though it were the direct sum of $\dim(V)$ copies of $B$ in $\CC$ (with the appropriate grading shifts). Then $\Hom_{\CC_{\enrich}}(B \sqot V, B'
\sqot V')$ can be viewed appropriately as matrices of morphisms between $B$ and $B'$. However, the action of $\sltwo$ on $V$ and $V'$ is used to twist the action on these matrices of
morphisms. If $\dim(V) > 1$ then $B \sqot V$ is decomposable as an object in the underlying category $\CC$, but can be $\sltwo$-indecomposable when $V$ is indecomposable.


We develop the enriched category more formally in the sequel. This method of enrichment is similar to techniques one encounters when studying the
equivariantization-deequivariantization principle, see e.g. \cite[\S 4]{DGNO}. In the next section we discuss a practical method to find an $\sltwo$ submodule idempotent whose image is
$B \sqot V$, generalizing the equations in \eqref{eq:factoringidempotent}.

Representable modules from the enriched category should also be considered as cofibrant. Taking the $\sltwo$-enriched category is analogous to taking the additive closure of a
preadditive category, and adjoining the images of cofibrant $\sltwo$ subquotient idempotents is analogous to taking the Karoubi envelope. Both are operations one typically does
before attempting to take a Grothendieck group. Below, we refer to the partial idempotent completion which adds cofibrant $\sltwo$ subquotient idempotents as the \emph{cofibrant Karoubi envelope}.

At the moment there are only conjectural ideas for how to build a triangulated category from an $\sltwo$-enriched category with an appropriate Grothendieck group. These ideas
($\sltwo$-omological algebra) will be explored in the sequel. However, one can explicitly state conjectures for the expected behavior of the $2$-categories studied in this paper, which
would be required for this triangulated category to be well-behaved, such as our $\sltwo$-adapted Soergel categorification conjecture in \S\ref{subsec:Soergelconjecture} below.

\subsection{Decompositions and the core}

For sanity, throughout this section we work over a field of characteristic zero. Recall the divided powers operators $\dx^{(k)} := \frac{\dx^k}{k!}$ for $k \ge 1$. It may help the reader to read this section concurrently with the following one, which gives examples of the general phenomena discussed here.

Suppose we have objects $X$ and $B$ in an $\sltwo$ category, and we wish to construct an $\sltwo$ submodule idempotent whose image is $B \sqot V$ where, for example, $V$ is the three-dimensional irreducible module for $\sltwo$. In the underlying additive category, we would need to prove that
\[ B(-2) \oplus B(0) \oplus B(+2) \sumset X, \]
which we would do by constructing inclusion and projection maps. Specifically, we would find maps $p_j \co X \to B$ of degree $j$ for $j \in \{-2,0,2\}$, and inclusion maps $i_j$ of degree $-j$, satisfying
\begin{equation} \label{eq:orthonormal} p_j \circ i_k = \delta_{jk} \id_B. \end{equation}
We set $e_j = i_j \circ p_j$ and $e = \sum e_j$.

Suppose can find maps $p_{-2}$ and $i_{2}$ of degree $-2$, satisfying the equations
\begin{equation} \label{eq:stringreq} \ddown(p_{-2}) = 0 = \ddown(i_{-2}), \qquad \dup^{3}(p_{-2}) = 0, \qquad \dup^{(2)}(p_{-2}) \circ i_{-2} = \id_B. \end{equation}
Let us set 
\begin{equation} p_0 = \dup(p_{-2}), \qquad p_2 = \dup^{(2)}(p_{-2}).\end{equation}
Then $\dup(p_2) = 0$, and $\{p_{-2}, p_0, p_2\}$ is the standard basis for $V$ (styled as a Weyl module) sitting isomorphically inside the core of $\Hom(X,B)$. Similarly one can define
\begin{equation} i_0 = -\frac{1}{2} \dup(i_2), \qquad i_{-2} = \dup^{(2)}(i_2). \end{equation}
However, we make no assumption that $\dup^{3}(i_2) = 0$! Were this true, then $\{i_2, -i_0, i_{-2}\}$ would be the dual standard basis for $V$ (styled as a dual Weyl module) sitting inside $\Hom(B,X)$. In practice, $\Hom(B,X)$ may have no finite-dimensional $\sltwo$ subrepresentations.

Applying $\dup$ to the equality $p_2 i_2 = \id_B$ we obtain
\[ -2 p_2 i_0 = \dup(p_2) i_2 + p_2 \dup(i_2) = 0, \]
so $p_2 i_0 = 0$. Similarly, applying $\dup$ and $\ddown$ repeatedly we can deduce that \eqref{eq:orthonormal} holds.

Moreover, let $\phi = \dup(i_{-2})$, which measures the failure of the span of $\{i_2, i_0, i_{-2}\}$ to be a finite-dimensional subrepresentation. Again using the Leibniz rule, we deduce that
\begin{equation} p_j \circ \phi = 0 \end{equation}
for all $j \in \{-2,0,2\}$.

We leave the reader to verify that
\[ e_2 \dup(e_2) = i_2 p_2 \dup(i_2) p_2 = -2 i_2 p_2 i_0 p_2 = 0.\]
\[ e_2 \ddown(e_2) = i_2 p_2 i_2 \ddown(p_2) = i_2 p_0 \ne 0. \]
In particular, $e_2$ is a $\dup$-submodule idempotent and a $\ddown$-quotient idempotent, but not an $\sltwo$ subquotient idempotent. Similarly, $e_{-2}$ is a $\dup$-quotient idempotent and a $\ddown$-submodule idempotent, and $e_0$ is a subquotient for both, but neither of these are $\sltwo$ subquotient idempotents. However, recalling that $e = \sum e_j$, we have
\[ e \dup(e) = 0, \qquad e \ddown(e) = 0. \]
Thus $e$ is an $\sltwo$ submodule idempotent. We claim that $\Hom(X,-) \circ e$ is a submodule of $\Hom(X,-)$ which is isomorphic as an $\sltwo$-module to $\Hom(B,-) \sqot V$.

\begin{rem} The inclusion and projection maps $i_j$ and $p_j$ form an Fc-filtration on $B \sqot V$ for $\dup$ in the sense of \cite[\S 4.7]{EQpDGbig}. They also form an Fc-filtration for
$\ddown$, with the opposite order. \end{rem}

The reader may imagine how this generalizes to other irreducible representations of $\sltwo$. Finding projection and inclusion maps of negative-most degree which satisfy the analogue of
\eqref{eq:stringreq} is our primary tool for constructing cofibrant $\sltwo$-submodule idempotents. As one can see, this is intimately related with finding the core of $\Hom(X,B)$, and
proving that the maps within are split.

\begin{rem} \label{rem:technical} It is not obvious that every $\sltwo$-submodule idempotent of $X$ whose image is isomorphic to $B \sqot V$ should be constructible in this fashion, and
this is a slightly thorny and technical issue, not suitable for this brief overview. We explore this question in the sequel. \end{rem}

\begin{rem} That projection and inclusion maps of negative-most degree are in the kernel of $\ddown$ is a relatively ubiquitous phenomenon. Maps of negative-most degree are constructed
from diagrams without polynomials, which by the formulas in this paper are in the kernel of $\ddown$. Also, maps of negative-most degree often split when this degree is non-positive
(for the diagrammatic Hecke category in characteristic zero, this is related to the Soergel conjecture, and further to hard Lefschetz on multiplicity spaces, see \cite{EWRel}). \end{rem}

\begin{rem} One can define a \emph{Lefschetz form} on degree $-k$ elements of $\Hom(X,B)$, which sends $(f,g) \mapsto f \circ \dup^k(\bar{g})$, where $\bar{g}$ is the vertical flip of
$g$. This form is valued in $\End(B)$. If $B$ is absolutely indecomposable then $\End(B)$ is local with maximal ideal $\mathfrak{m}$ and quotient field $\End(B)/\mathfrak{m} \cong \Bbbk$. In this case
one can define a Lefschetz form valued in $\Bbbk$.

In the prototypical example above, when $i_2 = \overline{p_{-2}}$, then the assumption that $i_{2}$ pairs nontrivially against $\dup^{2}(p_{-2})$ is the nondegeneracy of a Lefschetz form. In the examples we have computed, this Lefschetz form satisfies the Hodge-Riemann bilinear relations. Thus the existence of cofibrant $\sltwo$ subquotient idempotents is related to the Hodge theory of the operator $\dup$, when restricted to the core of $\Hom(X,B)$. \end{rem}

The entire discussion above relates to proving that an $\sltwo$ subquotient idempotent has the appropriate factorization behavior, which implies it is cofibrant. What about the other $\sltwo$ subquotient idempotents in an $\sltwo$ decomposition?

If we wish to continue decomposing $X$, we need to examine the complement of $B \sqot V$. If it is not $\sltwo$-indecomposable, we seek a cofibrant $\sltwo$ submodule summand isomorphic
to $B' \sqot V'$ in this complement, which would be an $\sltwo$ subquotient summand of $X$. To apply our techniques above one does not necessarily need the core of $\Hom(X,B')$ to be
nontrivial; one needs the core of $\Hom(\Im(1-e),B')$ to be nontrivial instead. In other words, we want a projection map $p'_j : X \to B'$ of non-positive degree for which $\ddown(p'_j)
= 0$, and while $\dup^{-j+1}(p'_j)$ might not be zero, it must be zero modulo $\Hom(X,-) \circ e$. These secondary projection maps $p'_j$ live in the core relative to $e$, while $e$ was
generated by projection maps in a genuine core. Thus these projection maps live in an \emph{iterated core}: the core of some morphism space, considered modulo an $\sltwo$-invariant
ideal generated by elements in the core of some other morphism space (and so forth).

In the sequel we will carefully define the notion of an \emph{$\sltwo$-antastic filtration}, which is effectively an $\sltwo$ decomposition where each subquotient idempotent $e_i$ is
constructed using projection maps which are in the core modulo the ideal generated by $e_1 + \ldots + e_{i-1}$. For an $\sltwo$-antastic filtration, the subquotient idempotents $e_i$
have cofibrant images. We also permit $\sltwo$-antastic filtrations where all but one of the idempotents factors as above. The remaining idempotent will be cofibrant by the two-out-of-three principle.

\subsection{Behavior in examples} \label{subsec:iteratedcoreexamples}

To elaborate on the ideas of the previous section, and on the relationship between the core and the semisimplification, we discuss several miracles which occur in the diagrammatic Hecke
category $\HC(S_n)$, and try to avoid distracting the reader with the details. All these examples were computed explicitly in \cite{EQHecke}, and we give more precise references for
each example. Only the raising operator $\dup$ was studied in \cite{EQHecke}, but adding $\ddown$ into the mix is relatively easy (often it is zero for degree reasons). We
work over a field $\Bbbk$.

We begin with an example where everything is as easy as possible, studied in \cite[\S 4.3 and 4.4]{EQHecke}.

\begin{example} \label{ex:BsBs} There is an object $B_s$ (for any simple reflection $s$) satisfying \begin{equation} \label{BsBs} B_s B_s \cong B_s(-1) \oplus B_s(+1) \end{equation} in $\HC(S_n)$.
Since the endomorphism ring of $B_s$ is a local ring with quotient field $\Bbbk$, the semisimplification of $\End(B_s B_s)$ is isomorphic to $\Mat_2(\Bbbk)$. The degree $-1$ projection map $p_1 \co B_s B_s \to B_s(-1)$ is uniquely determined up to scalar, and the miracle is that $\dup^2(p_1) = 0$.

In close analogy to Example \ref{ex:NH2}, the core of the $\sltwo$ action on $\End(B_s B_s)$ is a subalgebra isomorphic to $\Mat_2(\Bbbk)$, which is isomorphic to $V \ot V^*$ as an
$\sltwo$ representation. The matrix entries have the form $i_j p_k$ for $j,k \in \{1,2\}$, where $i_j$ and $p_k$ are particular inclusion and projection maps for the decomposition
\eqref{BsBs}. Moreover, the core of $\Hom(B_s B_s, B_s)$ is precisely the span of $p_1$ and $p_2$, with $\dup(p_1) = p_2$ and $\dup(p_2) = 0$. Similarly, the core of $\Hom(B_s, B_s
B_s)$ is the span of $i_1$ and $i_2$ with $\dup(i_2) = i_1$ and $\dup(i_1) = 0$.

Note however that neither primitive idempotent $e_1 = i_1 p_1$ or $e_2 = i_2 p_2 = 1 - e_1$ is an $\sltwo$ submodule idempotent. We have $e_1 \ddown(e_1) = 0$ but $e_1 \dup(e_1) \ne 0$, and vice versa for $e_2$. Thus $B_s B_s$ is $\sltwo$-indecomposable.

In the enriched category, there is an isomorphism $B_s B_s \cong B_s \sqot \Bbbk^2$, where $\Bbbk^2$ is the standard representation of $\sltwo$. \end{example}

The next example was done in great detail in \cite[\S 6.4 and 6.5, clarified further in \S 6.8]{EQHecke}.

\begin{example} In the Hecke category of $S_3$, there is a specific object\footnote{Here $X = B_s B_t B_s$, $Y = B_s$, and $Z = B_{sts}$.} $X$ which splits as a direct sum of two non-isomorphic indecomposable objects $Y$ and $Z$. Let $J$ denote the Jacobson radical of the category. It is the case
that $\End(Y)/J \cong \Bbbk$ and $\End(Z)/J \cong \Bbbk$, so that $\End(X)/J$ is two-dimensional, spanned by the two idempotents which project to these summands. However, the core
of $\End(X)$ is one-dimensional, spanned by the identity map, and neither primitive idempotent is killed by $\dup$.

The decomposition $X \cong Y \oplus Z$ implies that $\Hom(X,-) \cong \Hom(Y,-) \oplus \Hom(Z,-)$, or more precisely
\begin{equation} \label{repfunctordecomp} \Hom(X,-) \cong \Hom(Y,-) \circ p_Y \oplus \Hom(Z,-) \circ p_Z, \end{equation}
where $p_Y$ and $p_Z$ are the projection maps.

Consider $\Hom(X,Y)$ and $\Hom(X,Z)$. Both are supported in non-negative degrees, and are spanned in degree zero by their respective projection maps $p_Y$ and $p_Z$. Thus for degree
reasons, $\ddown(p_Y) = 0$ and $\ddown(p_Z) = 0$, and if either Hom space has a nonzero core, it must be a trivial module spanned by the projection map.

The first miracle is that
\begin{equation} \dup(p_Y) = 0. \end{equation}
As a consequence, $\Hom(Y,-) \circ p_Y \subset \Hom(X,-)$ is preserved by the $\sltwo$ action, making the decomposition \eqref{repfunctordecomp} a filtration with respect to $\sltwo$. In particular, this submodule $\Hom(Y,-) \circ p_Y$ is generated by the core of $\Hom(X,Y)$.

Meanwhile, $\dup(p_Z) \ne 0$, and $\core(\Hom(X,Z)) = 0$. The second miracle is that
\begin{equation} \label{miracle2} \dup(p_Z) \in \Hom(Y,Z) \cdot p_Y. \end{equation}
We have already observed by $\Hom(Y,Z) \cdot p_Y$ is preserved by $\sltwo$ inside $\Hom(X,Z)$, and \eqref{miracle2} implies that the image of $p_Z$ spans the core of the quotient module $\Hom(X,Z)/(\Hom(Y,Z) \cdot p_Y)$.

Hence the splitting of \eqref{repfunctordecomp} as modules over $\HC$ becomes a filtration with respect to the $\sltwo$ action, with $\Hom(Y,-) \circ p_Y$ being the submodule. Moreover, each layer of the filtration is generated by its core modulo the ideal generated by the previous part of the filtration.

In contrast, if $i_Y$ is the corresponding inclusion map in $\Hom(Y,X)$, we have $\dup(i_Y) \ne 0$. This explains why the corresponding idempotent $e_Y = i_Y p_Y$ is not killed by
$\dup$. However we have the miracles
\begin{equation} \dup(i_Z) = 0, \qquad \dup(i_Y) \in i_Z \circ \Hom(Y,Z). \end{equation}
Working with contravariant representable functors instead, we see that the decomposition
\begin{equation} \Hom(-,X) = i_Y \circ \Hom(-,Y) \oplus i_Z \circ \Hom(Z,-) \end{equation}
of functors is an filtration for $\sltwo$, where this time $i_Z \circ \Hom(Z,-)$ is the submodule.

Another consequence of these miracles is that $e_Y \dup(e_Y) = 0$, so $e_Y$ is an $\sltwo$ subquotient idempotent with complement $e_Z$. This equips $X$ with an $\sltwo$-antastic
filtration. \end{example}

In the above example, one could not find sufficiently many projection maps in the core to split $X$ into indecomposable direct summands. However, one could construct \emph{iterated
cores} by taking the core, quotienting by the ideal it generates, taking the core of what remains, quotienting by the ideal it generates, and so forth. One could find sufficiently
many projection maps for $X$ within the iterated cores. Similar to Conjecture \ref{conj:coreissplit}, we conjecture that iterated cores do not intersect the Jacobson radical (in the
$2$-categories we consider in this paper). A thorough development of iterated cores and a formal statement of this conjecture will be found in the sequel to this paper.

%

Iterated cores are still insufficient to split any object into indecomposable direct summands, as seen in the next example, which is studied in \cite[\S 9]{EQHecke}.

\begin{example} \label{ex:S8} In the Hecke category of $S_8$, there is an object\footnote{Here $X$ is the Bott-Samelson bimodule associated to the sequence $35246135724635$, with top summand $Z$, and
$Y$ is the indecomposable Soergel bimodule associated to the element $232565$.} $X$ which splits as a direct sum of two non-isomorphic indecomposable objects $Y$ and $Z$. The object $Z$
is technically not in the category but in the Karoubi envelope. However, there is no splitting of $\Hom(X,-)$ or $\Hom(-,X)$ into direct summands as a left module over $\HC$ which is
filtered with respect to $\sltwo$; neither $e_Y$ nor $e_Z$ is an $\sltwo$ subquotient idempotent. These two additive summands are entangled by the $\sltwo$ action. At the same time, the
cores of both $\Hom(X,Y)$ and $\Hom(Y,X)$ are zero. We can not even define an $\sltwo$ structure on $\Hom(Z,-)$ or $\Hom(-,Z)$. \end{example}

In the example above, the $\sltwo$ structure can not be used to find a complement for the Jacobson radical, and $X$ is $\sltwo$-indecomposable. This may be a feature rather than a bug! In characteristic $2$ the splitting $X \cong Y \oplus Z$ does not hold, and $X$ is indecomposable; perhaps this is detected by the $\sltwo$ action.

\subsection{The $\sltwo$-adapted Soergel categorification conjecture} \label{subsec:Soergelconjecture}

We assume the reader is familiar with the ordinary Soergel categorification theorem, see \cite[\S 5.5 and 11]{EMTW} or \cite[Theorems 3.15 and 6.25]{EWGr4sb}.

\begin{conj} \label{conj:sl2soergel} Let us work in the $\sltwo$-enriched category of $\HC(S_n)$. There are objects $D_y$ in the cofibrant Karoubi envelope, which satisfy and are inductively defined by the following properties. \begin{enumerate}
\item Let $\un{w}$ be a reduced expression for an element $w \in S_n$. The Bott-Samelson object $\BS(\un{w})$ has an $\sltwo$-antastic filtration, and all but one subquotient in the filtration is isomorphic to $D_y \sqot V$ for some finite-dimensional $\sltwo$ representation $V$ and some $y < w$ in $S_n$. The remaining subquotient is therefore cofibrant. This last subquotient is $\sltwo$-indecomposable and not isomorphic to $D_y \sqot V$ for any $y < w$ and any $V$, and we denote it by $D_{\un{w}}$.
\item The morphisms known as rex moves produce isomorphisms between $D_{\un{w}}$ and $D_{\un{w}'}$ for any two reduced expressions for the same element $w$. We let $D_w$ denote an object in this isomorphism class.
\item Any Bott-Samelson object $\BS(\un{x})$, whether $\un{x}$ is reduced or not, has a $\sltwo$-antastic filtration whose subquotients are isomorphic to $D_y \sqot V$ for various finite-dimensional $V$ and various $y \in S_n$.
\end{enumerate}
\end{conj}

The ordinary Soergel categorification theorem can be stated in almost identical fashion, except using ordinary direct sum decompositions (with grading shifts) rather than
$\sltwo$-antastic filtrations. This is used to define the indecomposable Soergel bimodule $B_w$ for $w \in S_n$. Note that $D_w$ need not be isomorphic to $B_w$ in the underlying
additive category. In example \ref{ex:S8}, letting $w = 35246135724635$ which is short for $s_3 s_5 s_2 \cdots$, the Bott-Samelson bimodule $X = \BS(\un{w})$ is $\sltwo$-indecomposable,
so $D_w = X$. However, $D_w$ is not indecomposable, and its top indecomposable summand $B_w$ was called $Z$ above.

If this conjecture holds then, in some conjectural triangulated category, we expect that the objects $D_w$ will descend to a basis in the triangulated Grothendieck group, as a free
module over the Grothendieck ring of finite-dimensional $\sltwo$ representations. Because an $\sltwo$-antastic filtration descends to an Fc-filtration for $\dup$, we will obtain an
analogous result for the $p$-dg Hecke category, as a free module over the cyclotomic ring $\Z[\zeta]$ where $\zeta^2$ is a primitive $p$-th root of unity.

The Grothendieck ring of this triangulated category should provide a $\Z[q+q^{-1}]$-form of the Hecke algebra. The structure coefficients for the basis $[D_w]$ are unimodal polynomials, because they are the graded dimensions of finite-dimensional $\sltwo$ representations. This does not imply the unimodality of Kazhdan-Lusztig structure coefficients, since $D_w \ne B_w$.

We think of Conjecture \ref{conj:sl2soergel} as stating that the Hecke category ``has enough'' $\sltwo$-antastic filtrations, or that it ``has enough'' split maps in the iterated
cores. When trying to generalize this conjecture to the other $2$-categories studied in this paper, it is not easy to rigorously state what it means to ``have enough,'' without pinning
down implicitly some suspected collection of objects, like $\{D_w\}$, which generate the category. We have failed to find a concise way to discuss the generalizations of Conjecture \ref{conj:sl2soergel}, so we leave further discussion to the next paper where these ideas are developed more formally.

\section{$\sltwo$ action on the KLR algebra}
\label{sec:KLR}

\subsection{Definition of the action}

\begin{defn} Let $\UC^+(\gg)$ denote the Khovanov-Lauda-Rouquier category associated to an oriented simply laced root datum. A presentation by generators and relations can be found
in \cite[\S 4.1]{KQ}. \end{defn}

For the reader familiar with KLR algebras: a KLR algebra is determined by the polynomials $Q_{ij}$ in the double crossing relation, for each pair of vertices $i$ and $j$ in the
(oriented) Dynkin diagram. If $i$ and $j$ are not connected by an edge, $Q_{ij} = 0$ as usual. If an edge is oriented from $i$ to $j$, then $Q_{ij} = x_i - x_j$. In our action of
$\sltwo$ below, we will have 
\begin{equation} \label{actonxiinKLR} \dup(x_i) = x_i^2 ,\qquad \ddown(x_i) = 1, \end{equation}
for all\footnote{Do not confuse the index $i$, which represents a node in the Dynkin diagram, with the indices
in $R_n = \Bbbk[x_1, \ldots, x_n]$. For each node $i$ one has an $i$-colored nilHecke algebra, containing a polynomial ring $R_n$ which is the $n$-fold tensor product of
$\Bbbk[x_i]$.} $i$. It is crucial here that $\ddown(Q_{ij}) = 0$, so one can not use the alternative choice $Q_{ij} = x_i + x_j$.

\begin{rem} The original definition of Khovanov-Lauda \cite{KhoLau09} used $Q_{ij} = x_i + x_j$, while the signed variant was introduced in \cite[page 15]{KhoLau11}. It was correctly
predicted on \cite[page 17]{KhoLau11} that this signed variant matches Lusztig's geometric categorification of the quantum group. As shown explicitly by Lauda in
\cite[Proposition 3.4]{LaudaParameters}, different choices of $Q_{ij}$ yield isomorphic categories when the Dynkin diagram is simply connected. This isomorphism will rescale the polynomials $x_i$ by
various signs. Evidently one could transfer our $\sltwo$ action along Lauda's isomorphisms, but they would no longer satisfy \eqref{actonxiinKLR}. \end{rem}

\begin{defn} \label{defn:sl2KLR} Equip $\UC^+(\gg)$ with an $\sltwo$ action as follows. The derivation $\dup$ is the derivation $\partial_1$ defined in \cite[Definition 4.13]{KQ}.
The derivation $\ddown$ sends a dot to $1$, and kills all crossings. In formulas, where red and blue are adjacent colors and red and green are distant, we have
\begin{subequations}
\begin{equation} \dup \left(~\ig{1}{dot}~ \right) = {
\labellist
\small\hair 2pt
 \pinlabel {$2$} [ ] at 7 19
\endlabellist
\centering
\ig{1}{dot}
}, \qquad \dup \left(~\ig{1}{Xii} ~\right) = - \ig{1}{XiiNW} - \ig{1}{XiiSE}, \end{equation}
\begin{equation} \dup \left(~\ig{1}{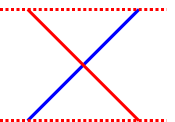}~ \right) = \ig{1}{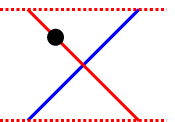}, \qquad \dup \left(~\ig{1}{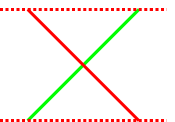} ~\right) = 0, \end{equation}
\begin{equation} \ddown \left(~\ig{1}{dot} ~\right) = \ig{1}{line}, \qquad \ddown \left(~\ig{1}{Xii}~ \right) = 0, \end{equation}
\begin{equation} \ddown \left(~\ig{1}{Xij}~ \right) = 0, \qquad \ddown \left(~\ig{1}{Xik}~ \right) = 0. \end{equation}
\end{subequations}
One extends the derivation to all diagrams using the monoidal Leibniz rule.
\end{defn}

\begin{thm} \label{thm:sl2KLR} The action of $\sltwo$ on $\UC^+(\gg)$ given in Definition \ref{defn:sl2KLR} is well-defined. The divided power operators $\dup^{(k)}$ and
$\ddown^{(k)}$ are well-defined in the integral form for all $k \ge 0$, making $\UC^+(\gg)$ into a divided powers $\sltwo$-algebra. \end{thm}
	
\begin{proof} First we need to check that the derivations $\dup$ and $\ddown$ are well-defined. By Lemma \ref{lem:defderivongen}, we need only check that they satisfy the relations of $\UC^+(\gg)$. That $\dup$ preserves the relations was checked\footnote{In \cite[end of \S 4.1]{KQ}, a multiparameter family of
degree $+2$ derivations is defined, and it is checked throughout the entirety of that section that the relations are preserved by each derivation in this family. The differentials in
\cite[Definition 4.13]{KQ} are special members of this family, c.f. \cite[Proposition 4.11]{KQ}.} in \cite[\S 4.1]{KQ}. Let us check that $\ddown$ preserves the relations.
	
Any diagram without dots is sent to zero by $\ddown$. Any relation which is a linear combination of diagrams without dots is therefore preserved by $\ddown$. The remaining relations
are these.

\begin{equation} \label{eq:dotslidesame} \ig{1}{line} \ig{1}{line} \quad = \quad \ig{1}{XiiNW} \; - \; \ig{1}{XiiSE}. \end{equation}
Applying $\ddown$ will kill the identity map, and send the other diagrams to the crossing,
\begin{equation} 
\ddown \left(~ \ig{1}{XiiNW}~ \right) \quad = \quad \ig{1}{Xii} \quad = \quad \ddown \left(~ \ig{1}{XiiSE}~ \right), 
\end{equation}
so this relation is preserved. The vertical flip is also a relation, and preserved for the same reason.

Similarly, $\ddown$ will preserve the relation
\begin{equation} \label{eq:dotslideadj} \ig{1}{XijNW} \quad = \quad \ig{1}{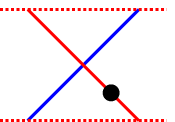}, \end{equation}
because
\begin{equation}
 \ddown \left(~ \ig{1}{XijNW}~ \right) \quad = \quad \ig{1}{Xij} \quad = \quad \ddown \left(~ \ig{1}{XijSE} ~\right).
  \end{equation}
Variants of this relation (e.g. when colors are distant, or the vertical flip) are preserved by $\ddown$ for the same reason.

The final relation with dots involves two adjacent colors.
\begin{equation} \label{eq:R2ij} 
\ig{1}{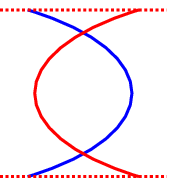} \quad = \quad \pm \left(~
 \ig{1}{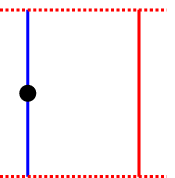} \quad - \quad \ig{1}{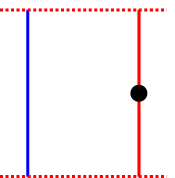} ~\right).
\end{equation}
The sign is determined by the orientation of the edge connecting these two colors, and will not affect the rest of the computation. Applying $\ddown$ to the LHS yields zero. The RHS is $Q_{ij} = \pm (x_i - x_j)$, and $\ddown(x_i - x_j) = 0$.

Hence $\ddown$ is a well-defined derivation. Now we wish to check that $(\dup, \dh, -\ddown)$ is an $\sltwo$ triple, and also that divided powers are well-defined integrally. By Lemmas \ref{lem:checksl2ongen} and \ref{lem:checkdivpowerongen}, we need only check this on the generators. 

Certain Hom spaces are clearly identified as $(R, \sltwo)$-modules for various polynomial rings $R$. For example, when $i$ and $j$ are adjacent colors, $\Hom(\EC_i \EC_j, \EC_j \EC_i) \cong \Z[x_i, x_j]$ as left $\Z[x_i, x_j]$-modules. By construction $\ddown$ kills the generator of the free rank $1$ module, and $\dup$ multiplies it by $x_j$, so the result is precisely the module $\Z[x_i, x_j]\ang{x_j}$. (We should also confirm that the generator lives in degree $+1 = \smsum(x_j)$.) We have already proven in Proposition \ref{prop:Rnp} that this is a well-defined $\sltwo$-representation with integrally-defined divided powers. Here
are several such isomorphisms (here $i$ and $j$ are adjacent colors, and $i$ and $k$ are distant): 
\begin{equation} 
\End(\EC_i) \cong \Z[x_i]\ang{0}, \qquad \Hom(\EC_i \EC_k,
\EC_k \EC_i) \cong \Z[x_i,x_k]\ang{0}, \qquad \Hom(\EC_i \EC_j, \EC_j \EC_i) \cong \Z[x_i, x_j] \ang{x_j}.
 \end{equation}
This takes care of all generators but the crossing in $\End(\EC_i \EC_i)$.

It is easy to verify that
\begin{equation} 
\ddown \dup \left(~\ig{1}{Xii} ~\right) \quad = \quad - 2 \ig{1}{Xii} \end{equation}
as desired. It is easy to verify that
\begin{equation}
 \dup^2\left(~ \ig{1}{Xii}~ \right) \quad = \quad 2 \ig{1}{XiiNWSE}, \end{equation}
\begin{equation} 
\dup^3\left( ~\ig{1}{Xii}~ \right) = 0. 
\end{equation}
One can see this also from the matrix of \eqref{NH2matrix}. In particular, when acting on the $(i,i)$ crossing, $\dup^{(2)}$ is well-defined integrally and $\dup^{(k)} = 0$ for all $k \ge 3$.
\end{proof}

\subsection{Action on the basis of crossings}

It was proven by in \cite[Theorem 2.5]{KhoLau09} that each morphism space in $\UC^+(\gg)$ is a free left (or right) module over the polynomial subring generated by dots, and that a
basis can be constructed using diagrams built entirely out of crossings.

More precisely, let $\ui = (i_1, \ldots, i_n)$ and $\uj = (j_1, \ldots, j_n)$ be sequences of colors, both having length $n$. Let $(S_n)_{\ui}^{\uj}$ denote the subset of $S_n$
consisting of those permutations $w$ for which $i_k = j_{w(k)}$ for all $1 \le k \le n$. This subset is empty unless the number of strands of each color agrees in $\ui$ and $\uj$.
If this happens, and there are $a_m$ strands with the color $m$, then restricting $w \in (S_n)_{\ui}^{\uj}$ to the strands of color $m$ we get a permutation in $S_{a_m}$. This
induces a bijection between $(S_n)_{\ui}^{\uj}$ and $S_{a_1} \times \cdots \times S_{a_d}$, and this bijection preserves the Bruhat order. In each case, the Bruhat order is
generated by the operation which removes a single crossing between two same-colored strands, when the result is still a reduced expression.

Pick a reduced expression of $w \in (S_n)_{\ui}^{\uj}$, and let $\psi_w$ denote the crossing diagram for that reduced expression, with the strands colored to represent a morphism in $\Hom(\EC_{\ui},
\EC_{\uj})$. Then, as a left module over $R_{\uj}$ (dots acting on the target sequence $\uj$), or as a right module over $R_{\ui}$ (dots acting on the source sequence $\ui$), we
have \begin{equation} \{\psi_w\}_{w \in (S_n)_{\ui}^{\uj}} \quad \text{is a basis for} \quad \Hom(\EC_{\ui},\EC_{\uj}). \end{equation}

\begin{thm} \label{thm:KLRassgr} For any $w \in (S_n)_{\ui}^{\uj}$ we have
\begin{equation} \ddown(\psi_w) = 0, \qquad \dup(\psi_w) \in \Span \{\psi_v\}_{v \in (S_n)^{\uj}_{\ui}, v \le w}. \end{equation}
Thus $\{\psi_w\}$ induces a downfree filtration of $\Hom(\EC_{\ui},\EC_{\uj})$, parametrized by $(S_n)^{\uj}_{\ui}$ with its Bruhat order. In the associated graded, the span of $\psi_w$ is a copy of $R_n \ang{p(w)}$, where a formula for $p(w)$ is given in \eqref{pwleft} or \eqref{pwright}, depending on whether we study the left $R_{\uj}$-action or the right $R_{\ui}$-action. \end{thm}

\begin{proof} Let us first prove the statement in the context of left modules over $R_{\uj}$. Clearly $\ddown$ kills $\psi_w$ as desired. The result is easiest to see using a different formula for the action of $\dup$ on
same-colored crossings, namely 
\begin{equation} \label{dXalt} 
\dup \left(~ \ig{1}{Xii}~ \right) \quad = \quad \ig{1}{line} \ig{1}{line} \; - 2 \; \ig{1}{XiiNW}. 
\end{equation}
We also recall
\[ 
\dup \left(~\ig{1}{Xij} ~\right) = \ig{1}{XijNW}, \qquad \dup \left(~\ig{1}{Xik} ~\right) = 0. 
\]
When we apply $\dup$ to $\psi_w$, we take the sum of $\dup$ applied to each crossing: up to linear combinations, this will either add a dot to the northwest of a crossing, or remove the crossing.

The first term on the RHS of \eqref{dXalt} says to remove a same-colored crossing. This produces an expression for an element $v < w$ which is smaller in the Bruhat order. Note that
this may not be a reduced expression for $v$, and even if it is, it may not be the chosen reduced expression $\psi_v$, so relations (such as \eqref{eq:R2ij}) must be applied to
rewrite this as a linear combination of basis diagrams. This rewriting process may produce diagrams which are even lower in the Bruhat order (c.f. the adjacent-colored Reidemeister
III relation \cite[(2.8)]{KhoLau09}), but it is well-known that the relations of the KLR algebra can be used to simplify an arbitrary diagram to a basis diagram without ever going
upwards in the Bruhat order. We state this as Lemma \ref{lem:simplifygreedy} below. So whenever we remove a crossing, we get something in the
$R_{\uj}$-span of $\psi_{v'}$ for $v' < w$ in the Bruhat order.

Now suppose a dot is placed to the northwest of a crossing. This dot should be forced to the top of the diagram (using \eqref{eq:dotslidesame} or \eqref{eq:dotslideadj}). On
its way up, \eqref{eq:dotslidesame} can produce an error term where a same-colored crossing is removed; once again, this will be in the span of diagrams which are smaller in the
Bruhat order. To see what polynomial remains on the top, we need only count the number of times a dot appeared on each strand.

For the $k$-th strand on top (colored $j_k$), we count the ways a dot can be added to that strand. We get coefficient $-2$ each time that strand is the northwest strand of a
same-colored crossing, and $+1$ each time that strand is the northwest strand of an adjacent-colored crossing. Note that the $k$-th strand on top will cross the $\ell$-th strand on
top, and be the northwest strand, if and only if $k < \ell$ and $w^{-1}(\ell) < w^{-1}(k)$. Consequently, let
\begin{subequations} \label{pwleft} 
\begin{equation} p(w) = \sum_{k=1}^n x_k \cdot \left( -2 \cdot \# \{\ell > k \mid w^{-1}(\ell) < w^{-1}(k), j_\ell = j_k\} \; + \; \# \{\ell > k \mid w^{-1}(\ell) < w^{-1}(k), j_\ell \text{ adjacent to } j_k\} \right). \end{equation}
Hence $p(w)$ represents the action of $\dup$ on the associated graded in the downfree filtration. We may write $p(w)$ more succinctly using the dot product on the set of colors (see \cite[top of p3]{KhoLau09}) by the formula
\begin{equation} p(w) = \sum_{k=1}^n \left( \sum_{\ell > k, w^{-1}(\ell) < w^{-1}(k)} -j_k \cdot j_\ell \right) x_k . \end{equation}
\end{subequations}

Suppose we instead study this morphism space as a right $R_{\ui}$-module. Now it helps to use the following formulas instead.
\begin{equation} \label{dXaltalt} \dup \left(~ \ig{1}{Xii}~ \right) \quad = \quad - \; \ig{1}{line} \ig{1}{line} \; - 2 \; \ig{1}{XiiSE}, \qquad \dup \left(~ \ig{1}{Xij} ~\right) \quad = \quad \ig{1}{XijSE}. \end{equation}
The $k$-th strand on bottom receives a dot whenever it is the southeast strand in a crossing, and if it crosses the $\ell$-th strand on bottom, then $\ell < k$ and $w(\ell) > w(k)$. So let
\begin{subequations} \label{pwright} 
\begin{equation}  p(w) = \sum_{k=1}^n x_k \cdot \left( -2 \cdot \# \{\ell < k \mid w(\ell) > w(k), i_\ell = i_k\} \; + \; \# \{\ell < k \mid w(\ell) > w(k), i_\ell \text{ adjacent to } i_k\} \right), \end{equation}
\begin{equation} p(w) = \sum_{k=1}^n \left( \sum_{\ell < k, w(\ell) > w(k)} -i_{\ell} \cdot i_k \right) x_k. \end{equation} \end{subequations}
Then this $p(w)$ describes the associated graded as a right $R_{\ui}$-module. \end{proof}

\begin{rem} Applying the anti-involution which flips each diagram upside-down, we get the analogous result for the $\sltwo$-triple $(\bardup, \dh, -\ddown)$. This will swap the left and right action, so it will also swap \eqref{pwright} and \eqref{pwleft}. \end{rem}
	
In the proof above, we used the following statement.

\begin{lem} \label{lem:simplifygreedy} Let $D$ be a crossing diagram corresponding to an subexpression of a reduced expression for $w \in S_n$, where we remove one crossing. Then
$D$ is in the left $R_{\uj}$-span (or the right $R_{\ui}$-span) of $\{\psi_{v}\}_{v \le w}$. \end{lem}

This lemma is considered obvious by most people in the field (including us), for which reason it is difficult to cite. The result is stated without proof in \cite[(2.33)]{KhoLau09},
for example. It would not be hard to prove this lemma directly, but it would be several annoying pages without much payoff. We will continue the tradition and not prove the lemma.
One could also deduce the result from the main theorem in \cite{EDiamond}. 

	
While Theorem \ref{thm:KLRassgr} gives control over the associated graded, it does not describe the upper-triangular terms in $\dup$ explicitly. This was because rewriting non-reduced expressions or the wrong reduced expression in terms of the basis can be complicated. However, for the nilHecke algebra it is easy: nonreduced expressions are zero, and all reduced expressions for the same element are equal. For the nilHecke algebra it is not hard to give a precise formula for $\dup(\psi_w)$.

\begin{thm} Inside $\NH_n$ we have
\begin{equation} \label{NHdupformula} \dup(\psi_w) = p(w) \psi_w + \sum_{v <_1 w} (1 + 2 m_{v,w})\psi_v. \end{equation}
Here, $v <_1 w$ means that $v < w$ in the Bruhat order and $\ell(v) = \ell(w) - 1$. The integers $m_{v,w}$ will be described in the proof.
\end{thm}

Note that we study the left action of $R_n$ in this proof. We leave the adaptations for the right action to the reader.
	
\begin{proof} To any element $w \in S_n$ one can associate its inversion set, the set of pairs $(a, b)$ with $1 \le a < b \le n$ such that $w(a) > w(b)$. Each crossing in a reduced
expression for $w$ produces a single inversion $(a,b)$, and we refer to it as the \emph{$(a,b)$-crossing}. In other words, if $(a,b)$ is in the inversion set of $w$ then the $a$th
strand and the $b$th strand (on bottom) will eventually cross exactly once, and this crossing can be identified with the pair $(a,b)$. Note that $b$ will be the strand on the northwest of the $(a,b)$ crossing; we always write our inversions $(a,b)$ in order so that $a < b$.

Pick a triple $a < b < c$ and consider the ordered set $\{(a,b) < (a,c) < (b,c)\}$, called the \emph{packet} of the triple $(a,b,c)$. A quick examination (try to draw it) should
convince the reader that the inversion set of $w$ intersected with this packet is either a prefix or a suffix. Moreover, for any reduced expression for $w$, the crossings in a
packet either appear in lexicographic order (if a prefix) or antilexicographic order (if a suffix). Only when the entire packet is contained in the inversion set of $w$ can they
appear in either lexicographic or antilexicographic order (because the full set is both a prefix and a suffix), and the braid relation $sts = tst$ swaps lexicographic for
antilexicographic. For example, if $(b,c)$ is not an inversion, but $(a,c)$ is, then $(a,b)$ must also be an inversion, and $(a,b)$ must occur below $(a,c)$ in a reduced expression
for $w$. These ideas are the start of Manin-Schechtmann's theory of higher Bruhat orders, see \cite{ManSch}.

If $v <_1 w$ then the inversion set of $v$ is equal to that of $w$ with exactly one inversion removed. Let $(a,b)$ be this inversion. Then (with regards to the left action of $R_n$)
we set \begin{equation} m_{v,w} := \# \{c \mid a < c < b \text{ and } w(b) < w(a) < w(c)\}. \end{equation} In this formula, the condition that $w(c) > w(a)$ is equivalent to the condition that $w(c) > w(b)$. After all, if $w(b) < w(c) < w(a)$ then the packet of $(a,c,b)$ will be entirely contained in the inversion set of $w$, and removing $(a,b)$ will not yield either a prefix or a suffix, so it will not yield a reduced expression.

We now argue that \eqref{NHdupformula} holds.

Applying $\dup$ to $\psi_w$ as in the proof of Theorem \ref{thm:KLRassgr}, we get a sum of diagrams where either a dot is added or a
crossing is removed. We have already computed that when all dots reach the top, the overall polynomial is $p(w)$. We need to compute the coefficient with which the $(a,b)$ crossing is removed, for some $(a,b)$ in the inversion set of $w$. Any crossing which occurs lower down in $\psi_w$, and for which the $a$ or $b$ strand is the northwest-southeast strand, will produce a dot which will eventually be forced through the $(a,b)$ crossing. Thus the coefficient involved in the removal of the $(a,b)$ crossing is \begin{itemize} \item $+1$ from \eqref{dXalt}, \item $+2$ if some $(c,b)$ crossing occurs below the $(a,b)$ crossing, \item $-2$ if some $(c,a)$ crossing occurs below the $(a,b)$ crossing. \end{itemize}

So let $c$ be any other strand. If some $(c,a)$ crossing occurs below the $(a,b)$ crossing then clearly $c < a < b$. Since both $(c,a)$ and $(a,b)$ crossings appear, so must
$(c,b)$, and since $(c,a)$ happens before $(a,b)$, they come in lexicographic order. So $(c,b)$ also occurs below $(a,b)$. Then the overall contribution to the coefficient is $-2 +
2 = 0$.  Thus no contribution of $-2$ can occur without being canceled by a contribution of $+2$. Conversely, if $c < a < b$ and a $(c,b)$ crossing occurs before an $(a,b)$ crossing, then the crossings appear in lexicographic order, meaning that a $(c,a)$ crossing must have come first. 

If some $(c,b)$ crossing occurs below the $(a,b)$ crossing, but a $(c,a)$ crossing does not also occur below the $(a,b)$ crossing, then we must have $a < c < b$. This will
contribute $+2$ to the coefficient, and $m_{v,w}$ is exactly counting such contributions. \end{proof}

\subsection{The nilHecke algebra and the matrix algebra} \label{subsec-NH}

We wish to justify some of the claims made in the introduction. For sake of sanity we work in characteristic zero. Let $n \ge 1$.

\begin{thm} \label{thm:NHmatrix} Let $\Bbbk = \Q$. As an $\sltwo$ representation, $\Mat_{n!}(\Bbbk)$ is a finite-dimensional subrepresentation of $\NH_n$, and is isomorphic to $\End_{\Q}(L_0 \ot L_1 \ot \cdots \ot L_{n-1})$, where $L_k$ is the irreducible representation of $\sltwo$ of dimension $k+1$. It is the core of $\NH_n$. \end{thm}

Most aspects of this theorem were proven in \cite[Proposition 3.24 and preceding]{KQ}. We give three proofs, mostly for pedagogical reasons. The example of $n=3$ is done explicitly
after the proofs, and it may help the reader to look at the proofs and the example simultaneously.

\begin{proof} Let $V$ denote the cohomology ring of the flag variety, thought of as the quotient of $R_n$ by the ideal generated by positive degree elements of $R_n^{S_n}$. Then $V$
is an $n!$-dimensional graded vector space, equipped with a natural action of $\dup$. By \cite[Example 2.2.3]{KhoLau09} there is a ring isomorphism $\NH_n \cong \Mat_{n!}(R_n^{S_n})$, and
hence a vector space isomorphism $\NH_n \cong \Mat_{n!}(\Bbbk) \ot_{\Bbbk} R_n^{S_n}$, where we identify $\Mat_{n!}(\Bbbk)$ with $\End_{\Bbbk}(V)$. It was proven in \cite[Proposition 3.24]{KQ}
that the action of $\dup$ respects this tensor product decomposition. From this we can immediately deduce that $\dup$ acts nilpotently on $\Mat_{n!}(\Bbbk) \subset \NH_n$, so this
subalgebra is contained in the core of $\NH_n$.

There are a number of ways to confirm that the core of $\NH_n$ is not bigger than $\Mat_{n!}(\Bbbk)$. Using the tensor decomposition $\NH_n \cong \End(V) \ot R_n^{S_n}$, one can
verify that nothing else is acted on nilpotently by $\dup$. Alternatively, the core of a (bounded below) $\sltwo$-algebra is a subalgebra by Proposition \ref{prop:coresubalgebra},
and any subalgebra of $\Mat_{n!}(R_n^{S_n})$ which properly contains $\Mat_{n!}(\Bbbk)$ is infinite-dimensional.

Note that $V$ is not naturally an $\sltwo$-representation: $\sltwo$ acts on $R_n$ but does not preserve the ideal generated by positive degree elements of $R_n^{S_n}$, so it does
not act on the quotient. However, we know that the core is some finite-dimensional $\sltwo$-representation, so its isomorphism class is determined by its graded dimension.
\end{proof}

Here is a second, more computational proof.

\begin{proof} Let $\delta = (n-1)x_1 + (n-2)x_2 + \cdots + 2 x_{n-2} + 1 x_{n-1} \in R_n$. Consider the $(R_n, \sltwo)$-module $R_n\ang{-\delta}$. As an $\sltwo$ representation, Proposition \ref{prop:coreRnp} states that
\begin{equation} \label{coredelta} \core(R_n\ang{-\delta}) \cong L_{n-1} \ot L_{n-2} \ot \cdots \ot L_0, \end{equation}
which is a representation of dimension $n!$. The same is true for $R_n\ang{-\delta'}$ where $\delta' = 1 x_2 + 2 x_3 + \ldots + (n-1) x_n$, though it may be nicer to order the tensor products in a fashion respecting the indices on the polynomials:
\begin{equation} \label{coredeltaprime} \core(R_n\ang{-\delta'}) \cong L_0 \ot \cdots \ot L_{n-2} \ot L_{n-1}. \end{equation}

In fact, a basis for $\core(R_n\ang{-\delta})$ is also a basis for $R_n$ as a free module over $R_n^{S_n}$. From Schubert theory, it is well-known that the polynomials
\begin{equation} \BM = \{x_1^{a_1} \cdots x_n^{a_n} \mid 0 \le a_i \le n - i \text{ for all } i\} \end{equation}
form a basis for $R_n$ over $R_n^{S_n}$. Meanwhile, by Proposition \ref{prop:coreRnp} and Corollary \ref{cor:corebasis}, they also form a basis for the core of $R_n\ang{-\delta}$.  In similar fashion, for $\core(R_n\ang{-\delta'})$ we can use the basis
\begin{equation} \BM' = \{x_1^{b_1} \cdots x_n^{b_n} \mid 0 \le b_i \le i-1 \text{ for all } i\}. \end{equation}

There is a map of $R_n$-bimodules
\begin{equation} R_n \ot_{\Bbbk} R_n \to \NH_n, \qquad f \ot g \mapsto f \psi_{w_0} g. \end{equation}
To make this a map of graded $R_n$-bimodules, we need to shift the source so that $1 \ot 1$ lives in degree $-n(n-1)$. Keeping track of the action of $\dup$, we get an $\sltwo$-intertwiner
\begin{equation} \label{R2nintertwiner} \phi \co R_n \ot_{\Bbbk} R_n\ang{-\delta_l -\delta'_{r}} \to \NH_n. \end{equation}
Here we think of $R_n \ot_{\Bbbk} R_n$ as a polynomial ring in $2n$ variables, the left variables and the right variables, where $\delta_l$ uses the left variables and $\delta'_r$ uses the right variables. This map $\phi$ is known to be surjective, see e.g. the matrix basis of $\NH_n$ described in \cite[Proposition 2.16]{KLMS} or \cite[Proposition 3.3]{KQ}. It is not injective, because it factors through the quotient $R_n \ot_{R_n^{S_n}} R_n$, though $\Ker \phi$ must be an $\sltwo$-submodule. In fact, by counting graded dimensions, one can verify that
\begin{equation} R_n \ot_{R_n^{S_n}} R_n\ang{-\delta_l - \delta'_r} \to \NH_n \end{equation}
is an isomorphism.

Again by Proposition \ref{prop:coreRnp} and Corollary \ref{cor:corebasis}, we know that
\begin{equation} \core(R_n \ot_{\Bbbk} R_n\ang{-\delta_l - \delta'_r}) \cong \End(L_{n-1} \ot L_{n-2} \ot \cdots \ot L_0) \end{equation}
as $\sltwo$-representations, being spanned by $\BM \ot \BM'$. However, we are interested in the core of the quotient $R_n \ot_{R_n^{S_n}} R_n$, and the core of a quotient module can be both bigger and smaller than the original, see Remark \ref{rem:coreofquotient}. Thankfully, we also know that $\BM \ot \BM'$ is a basis of $R_n \ot_{R_n^{S_n}} R_n$ as an $R_n^{S_n}$ module, and goes to a basis of $\NH_n$ over $R_n^{S_n}$, see again \cite[Proposition 2.16]{KLMS}. So the map $\phi$ from \eqref{R2nintertwiner} is injective on the core.

As in the other proof, once one knows that the core of $\NH_n$ contains $\BM \ot \BM'$ or $\Mat_{n!}(\Bbbk)$, there are a number of ways to confirm that it is not bigger. \end{proof}

\begin{rem} Note that $\BM \ot \BM'$ is not quite the matrix basis of $\Mat_{n!}(\Bbbk) \subset \Mat_{n!}(R_n^{S_n})$, though it has the same span. To get the matrix basis one must use dual bases for $R_n$ over $R_n^{S_n}$, such as the Schubert and dual Schubert bases. \end{rem}

Here is a sketch of a third proof that the core contains $\End_{\Q}(L_0 \ot L_1 \ot \cdots \ot L_{n-1})$. We ignore the remaining parts of the theorem.

\begin{proof}(Sketch) By restricting from $R_n$ to $R_n^{S_n}$, $R_n\ang{-\delta}$ is an $(R_n^{S_n}, \mathfrak{sl}_2)$-module of rank $n!$ with the $\sltwo$-stable basis $\BM$. As was shown in \cite[Section 3.1]{KQ} (see the discussion there around equations (63) and (64)), the isomorphism
\begin{equation}
\NH_n \cong \End_{R_n^{S_n}}(R_n\ang{-\delta}),
\end{equation}
equips the nilHecke algebra with an $\sltwo$-action, which agrees with Definition \ref{defn:sl2KLR}.

We have already argued in \S\ref{subsec:core} that the core of a tensor product contains the tensor product of the cores. Since endomorphism rings are particular kinds of tensor
products (though one must be careful when taking tensor products of infinite-dimensional representations in this way), one deduces that the core of an endomorphism ring contains the
endomorphism ring of the cores. Thus $\core(\NH_n)$ contains $\End_{\Q}(\core(R_n\ang{-\delta}))$, which is $\End_{\Q}(L_0 \ot L_1 \ot \cdots \ot L_{n-1})$. \end{proof}

\begin{example}
 Let $n=3$. Here is a basis of $\core(R_3\ang{-2x_1 - x_2})$, making it clear the isomorphism with $L_2 \ot L_1$.
 \begin{equation}
 \begin{gathered}
 \xymatrix{
 1 \ar@/^/[rr]^{\dup=-2} \ar@/^/[dd]^{\dup=-1} && x_1 \ar@/^/[rr]^{\dup = -1} \ar@/^/[dd]^{\dup = -1} \ar@/^/[ll]^{\ddown=1} &&  x_1^2 \ar@/^/[dd]^{\dup=-1} \ar@/^/[ll]^{\ddown=2} \\ 
 && && \\
 x_2  \ar@/^/[rr]^{\dup = -2} \ar@/^/[uu]^{\ddown=1} && x_1x_2 \ar@/^/[rr]^{\dup =-1} \ar@/^/[ll]^{\ddown=1} \ar@/^/[uu]^{\ddown=1} && x_1^2x_2 \ar@/^/[uu]^{\ddown=1} \ar@/^/[ll]^{\ddown=2}
 }
 \end{gathered}
\end{equation}  
Here is another basis making clear the isomorphism with $L_3 \oplus L_1$.
 \begin{subequations}
\begin{equation*} 
\xymatrix{
1 \ar@/^/[rr]^{\dup=1} && -2x_1 - x_2  \ar@/^/[ll]^{\ddown=-3} \ar@/^/[rr]^{\dup=2}  &&  x_1^2 + 2 x_1 x_2  \ar@/^/[ll]^{\ddown=-2}  \ar@/^/[rr]^{\dup=3} && - x_1^2 x_2 \ar@/^/[ll]^{\ddown=-1}
}
\end{equation*}
\begin{equation}
\bigoplus
\end{equation}
\begin{equation*}
\xymatrix{
x_1 - x_2 \ar@/^/[rr]^{\dup=1} && x_1^2 - x_1 x_2  \ar@/^/[ll]^{\ddown=1} 
}
\end{equation*}
\end{subequations}
For those interested in integral structure, this second basis is not a $\Z$-basis, as the determinant of the change of basis matrix is $3$.

To obtain analogous bases of $R_3\ang{-x_2 - 2 x_3}$ just swap $x_1$ and $x_3$.

Now $\NH_3$ is isomorphic to a $6 \times 6$ matrix algebra over $R_3^{S_3}$, and the matrix entries correspond to $\psi_{w_0}$ with certain polynomials on top and on bottom, see
e.g. \cite[Proposition 3.3]{KQ}. The polynomials on top and those on bottom must be dual bases for $R_3$ over $R_3^{S_3}$ with respect to the Demazure operator $\partial_{w_0}$.
However, one can get a basis over $R_3^{S_3}$ (not necessarily a matrix basis) but choosing any two bases for $R_3$ over $R_3^{S_3}$, and placing them on top and bottom of
$\psi_{w_0}$. Choosing your favorite bases for $\core(R_3\ang{-2x_1 - x_2})$ and $\core(R_3\ang{-x_2 - 2 x_3})$ respectively for the top and bottom, the $\sltwo$ structure is
transparent. \end{example}

\section{$\sltwo$ action on the Hecke category} \label{sec:hecke}

\subsection{Definition of the action}

\begin{defn} Let $\HC = \HC(S_n)$ denote the diagrammatic Hecke category associated to the action of $S_n$ on $R_n = \Z[x_1, \ldots, x_n]$. A presentation by generators
and relations can be found in \cite{EKho}. \end{defn}

\begin{defn} \label{defn:sl2hecke} Equip $\HC$ with an $\sltwo$ action as follows. The derivation $\dup$ is the derivation defined in \cite[Theorem 2.5]{EQHecke}, where $g_i = x_i$.
The derivation $\ddown$ kills all diagrams without polynomials, and sends $x_i$ to $1$. In the formulas below, blue represents the simple reflection $s_i = (i,i+1)$. Red represents $s_{i+1}$, and green represents some color distant from blue.
\begin{subequations}
\begin{equation} \label{eq:ddots} \dup \left( \ig{1}{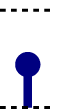} \right) = \; \poly{x_i} \ig{1}{enddotblue}, \qquad \dup \left( \ig{1}{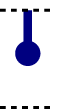} \right) = \; \poly{x_{i+1}} \ig{1}{startdotblue}, \end{equation}
\begin{equation} \label{eq:dtri} \dup \left( \ig{1}{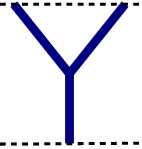} \right) = - \; {
\labellist
\tiny\hair 2pt
 \pinlabel {$x_i$} [ ] at 19 35
\endlabellist
\centering
\ig{1}{splitpoly}
}, \qquad \dup \left( \ig{1}{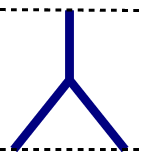} \right) = - \; {
\labellist
\tiny\hair 2pt
 \pinlabel {$x_{i+1}$} [ ] at 20 5
\endlabellist
\centering
\ig{1}{mergepoly}
}, \end{equation}
\begin{equation} \label{eq:dcups} \dup \left( \ig{1}{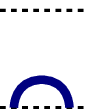} \right) = \; \ig{1}{enddotblue} \ig{1}{enddotblue}, \qquad \dup \left( \ig{1}{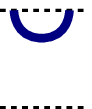} \right) = \; -\; \ig{1}{startdotblue} \ig{1}{startdotblue}, \end{equation}
\begin{equation} \label{eq:drex0} \dup \left( \ig{1}{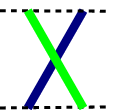} \right) = 0, \qquad \dup \left( \ig{1}{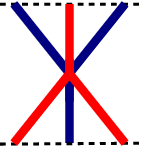} \right) = 0, \end{equation}
\begin{equation} \label{eq:d6} \dup \left( \ig{1}{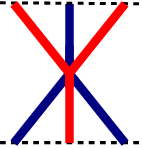} \right) = \; \ig{1}{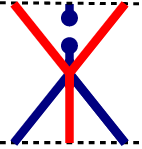} \; - \; \ig{1}{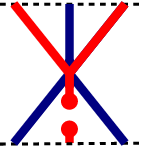}, \end{equation}
\begin{equation} \dup \left( \poly{x_i} \right) = \; \poly{x_i^2}, \qquad \ddown \left( \poly{x_i} \right) = \; \poly{}, \end{equation}
\begin{equation} \ddown \left( \ig{1}{enddotblue} \right) = 0, \qquad \ddown \left( \ig{1}{startdotblue} \right) = 0, \qquad \ddown \left( \ig{1}{splitpoly} \right) = 0, \qquad \ddown \left( \ig{1}{mergepoly} \right) = 0, \end{equation}
\begin{equation} \ddown \left( \ig{1}{capblue} \right) = 0, \qquad \ddown \left( \ig{1}{cupblue} \right) = 0, \end{equation}
\begin{equation} \ddown \left( \ig{1}{Xbg} \right) = 0, \qquad \ddown \left( \ig{1}{other6valent} \right) = 0, \qquad \ddown \left( \ig{1}{6valent} \right) = 0. \end{equation}
\end{subequations}

One extends the derivation to all diagrams using the monoidal Leibniz rule.
\end{defn}

\begin{thm} \label{thm:sl2hecke} The action of $\sltwo$ on $\HC$ given in Definition \ref{defn:sl2hecke} is well-defined. The divided power operators $\dup^{(k)}$ and
$\ddown^{(k)}$ are defined integrally for all $k \ge 0$, making $\HC$ into a divided powers $\sltwo$-algebra. \end{thm}

\begin{proof} As in the proof of Theorem \ref{thm:sl2KLR}, we will use the lemmas of \S\ref{subsec-leibniz} to reduce the amount of work we need to do. This time we do it tacitly.
	
First we need to check that the derivations $\dup$ and $\ddown$ preserve the relations of $\HC$. That $\dup$ preserves the relations was checked in \cite[Theorems 2.1, 2.3]{EQHecke}.
Let us check that $\ddown$ preserves the relations.

Any diagram without polynomials is sent to zero by $\ddown$. Any relation which is a linear combination of diagrams without polynomials is therefore preserved by $\ddown$. There is
only one relation that remains, the polynomial forcing relation. Below, blue represents $s_i$, and $j$ is arbitrary.
\begin{equation} \label{polyforce} \ig{1}{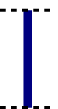} \poly{x_j} \; - \; \ig{1}{space} \poly{s_i(x_j)} \ig{1}{lineblue} \; = \; \pa_i(x_j) \ig{1}{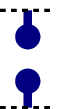}. \end{equation}
Note that $\pa_i(x_j)$ is a scalar. This relation is preserved by $\ddown$ since $\ddown(x_i) = 1$ for all $i$, so $\ddown$ kills the LHS, and the RHS is a diagram without polynomials so is also killed by $\ddown$. Hence $\ddown$ is a well-defined derivation.

Now we wish to check that $(\dup, \dh, -\ddown)$ is an $\sltwo$ triple, and also that divided powers are defined integrally. We can check these properties on the generators. Since $\dup$ and $\ddown$ raise or lower the degree appropriately, we need only check that $[\ddown, \dup] = \dh$. For each of the non-polynomial generators $\phi$, $\ddown(\phi) = 0$, and it is very easy to confirm that 
\begin{equation} \ddown(\dup(\phi)) = (\deg \phi) \cdot \phi = \dh(\phi). \end{equation}
Meanwhile, for $\End(\1) = R_n$, the action of $\sltwo$ is the standard one. This confirms that the $\sltwo$ action is well-defined. That $\dup^{(k)}$ is defined integrally was checked in \cite[\S 8]{EQHecke}. That $\ddown^{(k)}$ is defined integrally on the polynomial ring $\End(\1)$ was checked in Lemma \ref{lem:polyringiscool}.
That $\ddown^{(k)}$ is defined integrally on the other generators is easy, since it is zero for $k \ge 1$. \end{proof}
	
\subsection{Reminders: rex moves and lower terms} \label{subsec:rexmoves}

\begin{defn} For $w \in S_n$, let $I_{< w}$ denote the two-sided ideal spanned by all morphisms which factor through reduced expressions for elements $v \in S_n$ with $v < w$. \end{defn}

Given two reduced expressions $\uw, \uw'$ for the same element $w \in S_n$, Matsumoto's theorem states that they can be connected by a sequence of braid relations. To such a sequence
of braid relations there is a corresponding diagram built from $4$-valent and $6$-valent vertices, having source $\uw$ and target $\uw'$, which we call a \emph{rex move}. There are
many potential sequences of braid relations which go from $\uw$ to $\uw'$, and the corresponding rex moves in $\HC$ are not equal.

\begin{lem} Let $\uw, \uw'$ be reduced expressions for $w \in S_n$. For any two rex moves $\uw \to \uw'$, their difference lies in $I_{< w}$. \end{lem}

\begin{proof} This is proven in \cite[Lemma 7.4]{EWGr4sb}. \end{proof}

Rex moves are always in the kernel of $\ddown$, as any diagram is, but are not always in the kernel of $\dup$, see \eqref{eq:d6}. Two different rex moves with the same source and target can have different values of $\dup$. For example, let $s = s_i$ and $t = s_{i+1}$, and consider the reduced expression $(s,t,s)$. The identity map of this reduced expression is a rex move which is killed by $\dup$. Meanwhile, the path $(s,t,s) \to (t,s,t) \to (s,t,s)$ gives the doubled 6-valent vertex, which is not in the kernel of $\dup$. The reader versed with this diagrammatic calculus should have no trouble verifying that
\begin{equation} \label{ddouble6} \dup \left( \ig{1}{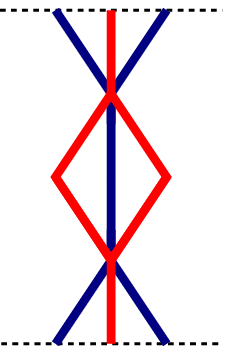} \right) = \; - \; \ig{1}{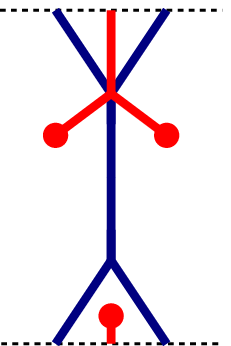} \; = \; - \left(\; \ig{1}{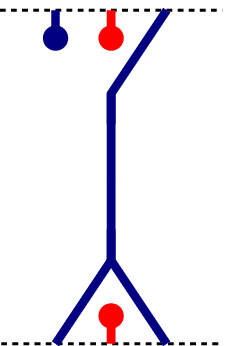} \; + \; \ig{1}{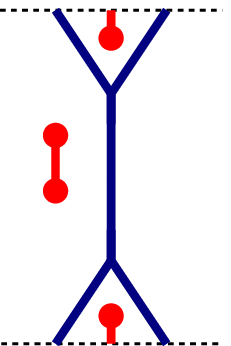}\; \right). \end{equation}
Note at least that $\dup$ sends the doubled 6-valent vertex to the ideal $I_{<sts}$ of lower terms.

\begin{rem} The right hand side of \eqref{ddouble6} is a left $R_n$-linear combination of double leaves, and double leaves form a basis for morphisms as a left (or right) $R_n$-module. Note the non-trivial root $\alpha_t$ which appears. This gives an example of
the kind of behavior discussed in Remark \ref{rem:rootsapain}. \end{rem}

\subsection{Reminders: subexpressions and light leaves}

Let $\ux = (s_{x_1}, \ldots, s_{x_d})$ be an expression of length $d$, and $\eb \subset \ux$ be a subexpression. To $\eb$ we associate a Bruhat stroll $(1 = w_0, w_1, \ldots, w_d)$
as in \cite[\S 2.4]{EWGr4sb}, where to get from $w_{i-1}$ to $w_i$ we multiply by either $s_{x_i}$ or by $1$ depending on $\eb$. Note that the Bruhat stroll determines the
subexpression, and vice versa. We refer to $w_d$ as the \emph{terminus} of $\eb$. We let $E(\ux,w)$ denote the set of subexpressions $\eb \subset \ux$ with terminus $w$. Following
ideas of Libedinsky \cite{LibLL}, to each $\eb \in E(\ux,w)$ we associate in \cite[Construction 6.1]{EWGr4sb} a morphism $\LL_{\eb} \co \ux \to \uw$ called a \emph{light leaf},
whose target is some reduced expression for $w$ (depending on $\eb$). Flipping this light leaf upside-down, we get a morphism $\GG_{\eb} \co \uw \to \ux$.

Let us remark on some important features of the light leaf construction. If $\uv$ is  a reduced expression for some $v \in S_n$, then it has a unique subexpression with
terminus $v$, the \emph{top subexpression}. Any rex move starting at $\uv$, including the identity map, is a valid light leaf for the top subexpression, and all light leaves for the top subexpression are rex moves. Let us note that any other  subexpression of $\uv$ has terminus $v' < v$, so its light leaf lies in the ideal $I_{< v}$.

If $\ux = \uv \uz$ is a concatenation of two smaller sequences, we can restrict a subsequence $\eb \subset \ux$ to a subsequence $\fb \subset \uv$. Suppose that $\uv$ is a reduced
expression for some $v \in S_n$, and $\eb$ restricts to the top subsequence $\fb \subset \uv$. Then we refer to $\LL_{\eb}$ as a \emph{light tail}, and since it is determined by the
restriction of $\eb$ to $\uz$, we may use the notation $\TL_{\eb \setminus \fb}$.

Finally, suppose that $\ux = \uy \uz$ is a concatenation of two smaller sequences, and $\eb \in E(\ux,w)$ restricts to $\fb \in E(\uy,v)$. Then $\LL_{\fb}$ is a map from $\uy$ to $\uv$ for some reduced expression for $v$. Meanwhile, there is a light tail $\uv \uz \to \uw$ determined by the subexpression $\eb \setminus \fb$ of $\uz$. The inductive construction of light leaves states that
\begin{equation} \label{LLinduct} \LL_{\eb} = \TL_{\eb \setminus \fb} \circ (\LL_{\fb} \ot \id_{\uz}). \end{equation}
Schematically, we draw
\begin{equation} \label{LLtail} {
\labellist
\small\hair 2pt
 \pinlabel {$\LL_{\eb}$} [ ] at 50 19
\endlabellist
\centering
\ig{1}{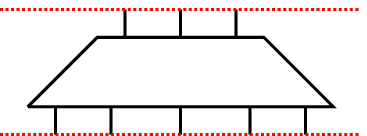}
}  \quad = \quad {
\labellist
\small\hair 2pt
 \pinlabel {$\uy$} [ ] at 6 9
 \pinlabel {$\uz$} [ ] at 71 9
 \pinlabel {$\uw$} [ ] at 69 67
 \pinlabel {$\uv$} [ ] at 7 42
 \pinlabel {$\LL_{\fb}$} [ ] at 23 26
 \pinlabel {$\TL_{\eb \setminus \fb}$} [ ] at 43 52
\endlabellist
\centering
\ig{1}{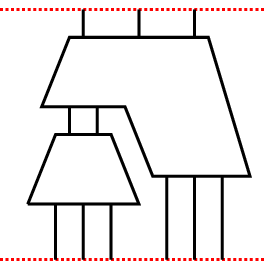}
} \end{equation}
When $\uz$ has length $1$, we think of $\TL_{\eb \setminus \fb}$ as being a single \emph{tier} of the light leaf algorithm. When $\uz$ has length $d$, the light tail is built inductively from the last $d$ tiers. This is discussed in \cite[Remark 6.4]{EWGr4sb}.

Because light tails will be important in some proofs below, let us introduce some terminology.

\begin{defn} Let $w \in S_n$. Relative to $w$, we call a sequence of simple reflections $\uz = (s_{z_1}, \cdots, s_{z_m})$ a \emph{tail expression}, and a sequence $\eb \subset \uz$
of $0$s and $1$s a \emph{tail subexpression}. A tail (sub)expression is the same thing as a (sub)expression, but we interpret its Bruhat stroll differently, and we decorate the
subsequence with $U$s and $D$s accordingly. There is a unique $v \in S_n$ such that $v s_{z_1}^{\eb_1} s_{z_2}^{\eb_2} \cdots s_{z_k}^{\eb_m} = w$. The \emph{tail Bruhat stroll}
associated to $w$ and $\eb \subset \uz$ is the Bruhat stroll which starts at $v$ rather than $1$, and wends its way to $w$. We call $v$ the \emph{start} of the tail Bruhat stroll.

Fix $v, w \in S_n$ and a reduced expression $\uv$ for $v$. For a tail expression $\uz$, there is a bijection between ordinary subexpressions $\eb \subset \uv \uz$ whose restriction
to $\uv$ is the top subsequence $\fb \subset \uv$, and tail subexpressions $(\eb \setminus \fb) \subset \uz$ with start $v$ and terminus $w$. This bijection is natural over the
choice of reduced expression for $\uv$, in the obvious sense. Henceforth (and unlike the previous paragraph), we always use $(\eb \setminus \fb) \subset \uz$ as notation for a tail
subexpression, even though $\eb$ and $\fb$ themselves need not have been chosen. Whenever we choose another sequence $\uy$ and consider the concatenation $\uy\uz$, and
whenever $\fb \in E(\uy,v)$, we will then set $\eb$ to be the subexpression of $\uy \uz$ whose restriction to $\uy$ is $\fb$ and whose restriction to $\uz$ is $\eb \setminus \fb$.
\end{defn}

\subsection{Reminders: double leaves}

For two sequences $\ux$ and $\ux'$, subexpressions $\eb \subset \ux$ and $\fb \subset \ux'$ are called \emph{coterminal} if they have the same terminus $w$, and we refer to $(\eb, \fb, w)$ as a \emph{coterminal triple} subordinate to $(\ux,\ux')$. Sometimes we omit $w$ from the triple, writing only $(\eb, \fb)$. To each coterminal triple we associate a \emph{double leaf} morphism in $\Hom(\ux, \ux')$, 
\begin{equation} \label{eq:DL} \DL_{\eb, \fb} := \GG_{\fb} \circ N \circ \LL_{\eb}, \end{equation}
where $N$ is some rex move from the target of $\LL_{\eb}$ to the source of $\GG_{\fb}$, both being reduced expressions for $w$. By taking one double leaf for each coterminal triple, one obtains a basis for $\Hom(\ux, \ux')$ as a left (or right) $R_n$-module. Double leaves are always diagrams without polynomials, so they are killed by $\ddown$.

Light leaves and double leaves are not determined only by the subexpressions $\eb$ and $\fb$. There are many choices of rex moves in the construction of each light leaf, as well as
the choice of rex move $N$ in \eqref{eq:DL}. In particular, the composition $N \circ \LL_{\eb}$ of a rex move with a light leaf is, itself, another valid choice of a light leaf
associated to $\eb$. In this way, we could remove the rex move $N$ from \eqref{eq:DL}, absorbing it into the light leaf $\LL_{\eb}$ (or into the upside-down light leaf $\GG_{\fb}$).

When one speaks about the double leaves \emph{basis}, one must choose one amongst the many possible double leaves for each coterminal triple to be a basis element. There are many
different double leaves bases. When we speak of light leaves or double leaves, we typically refer to the set of all possible maps produced by the algorithm, with the flexibility of
using arbitrary rex moves. When we speak of a \emph{distinguished} light leaf or double leaf, we must have fixed one for each subexpression, and we refer to that one. In this way we
can separate in our language between the rigid choices one must make to get a basis, and the flexible choices which are sufficient for a spanning set.

There is a filtration by the spans of certain double leaves with regards to a certain partial order, for which the image of a double leaf in the associated graded does not depend on
the choice of rex moves! This result will be proven in the next section. Thankfully, $\dup$ will preserve this (downfree) filtration, and one can easily compute its action on the
associated graded. One can not expect much more from the combinatorics of subexpressions: since $\dup$ acts nontrivially on rex moves, one should not expect a formula for the $\dup$
action on a double leaf which is independent of the choice of rex move.

\subsection{The lexicoBruhat order}

In this section we develop some technology for working with light leaves and double leaves. This technology is not original to this paper: it is part of work in progress \cite{EMSCI}
by the first author, and much was known (but not written in the literature) previously to the experts\footnote{The technology we present is closely related to ideas developed by
Elias and Williamson in their early attempts to prove that double leaves span. The eventual proof in \cite{EWGr4sb} used a different inductive proof of spanning. In the proof that
double leaves are linearly independent, \cite{EWGr4sb} used the path dominance order on triples, which is a stronger partial order than the lexicoBruhat order.}.

To study double leaves it helps to first study light leaves. From the fact that distinguished double leaves form a basis for $\Hom(\ux,\uw)$, we deduce\footnote{This was actually
proven first, and used to deduce that distinguished double leaves form a basis, see \cite[Proposition 7.6]{EWGr4sb}.} that distinguished in $E(\ux,w)$ form a basis for
$\Hom(\ux,\uw)/I_{< w}$. Until further notice, we will be studying this Hom space, modulo lower terms. In the following proposition we discuss not the span of particular light
leaves, but the span of all light leaves which are constructible by the non-deterministic algorithm of \cite[Construction 6.1]{EWGr4sb}. In other words, in this proposition we are
agnostic to the choice of rex moves.

\begin{prop} \label{prop:kform} Let $\uw$ be a reduced expression for some $w \in S_n$. Let $\ux = \uy \uz$ be the concatenation of a sequence $\uy$ of length $k$ and a sequence
$\uz$ of length $d-k$. Each subexpression of $\ux$ restricts to a subexpression of $\uy$. For $v \in S_n$, let $X_{\le v} \subset E(\ux,w)$ denote the subset which restricts to
$E(\uy,v')$ for some $v' \le v$. Let $H_{\le v}$ denote the subspace of $\Hom(\ux,\uw)$ spanned (over the left action of $R_n$) by all possible light leaves corresponding to subexpressions in $X_{\le v}$.
	
Let $f \co \uy \to \uv$ be any morphism whose target is a reduced expression $\uv$ for $v \in S_n$. Then
\begin{equation} \label{kform} g \circ (f \ot \id_{\uz}) \in H_{\le v} \text{ modulo } I_{< w} \end{equation}
for all $g \in \Hom(\uv \uz, \uw)$.
Schematically, we have
\begin{equation} {
\labellist
\small\hair 2pt
 \pinlabel {$f$} [ ] at 23 19
 \pinlabel {$g$} [ ] at 38 47
 \pinlabel {$\uv$} [ ] at 5 35
 \pinlabel {$\uy$} [ ] at 4 6
 \pinlabel {$\uz$} [ ] at 63 6
 \pinlabel {$\uw$} [ ] at 6 61
\endlabellist
\centering
\ig{1}{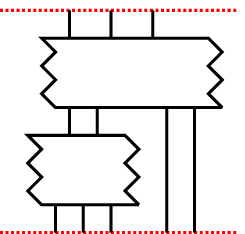}
} \quad = \sum R_n \cdot \;\;{
\labellist
\small\hair 2pt
 \pinlabel {$\uy$} [ ] at 6 9
 \pinlabel {$\uz$} [ ] at 71 9
 \pinlabel {$\uw$} [ ] at 69 67
 \pinlabel {$\uv'$} [ ] at 7 42
 \pinlabel {$\LL$} [ ] at 23 26
 \pinlabel {$\TL$} [ ] at 43 52
\endlabellist
\centering
\ig{1}{LLtail}
} + I_{< w}, \end{equation}
where $v' \le v$.
\end{prop}

\begin{proof} We suppose the result is true for all elements less than $v$ in the Bruhat order, and deduce the result for $v$. Then, by induction, the result will be true for any $v
\in S_n$ (without the need to check any base case). Note that everything else (i.e. $\uy$, $\uz$, $\uw$, $k$) is unchanging in this induction.

We know that $f \in \Hom(\uy,\uv)$ and $g \in \Hom(\uv \uz, \uw)$ are both in the left $R_n$-span of double leaves. So, up to taking left linear combinations over $R_n$, our morphism has the following form.
\begin{equation} {
\labellist
\small\hair 2pt
 \pinlabel {$\LL_1$} [ ] at 23 17
 \pinlabel {$\GG_1$} [ ] at 23 33
 \pinlabel {$\LL_2$} [ ] at 36 57
 \pinlabel {$\GG_2$} [ ] at 36 74
 \pinlabel {$\uv$} [ ] at 2 45
 \pinlabel {$\uv'$} [ ] at 6 27
 \pinlabel {$\uw'$} [ ] at 3 65
 \pinlabel {$\uw$} [ ] at 3 84
\endlabellist
\centering
\ig{1}{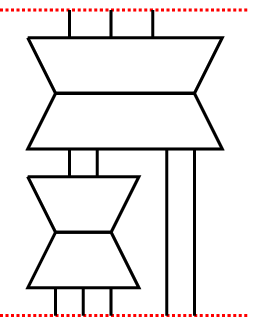}
} \end{equation}
The source $\uw'$ of $\GG_2$ is a reduced expression for some element $w' \le w$. Thus either $\GG_2 \in I_{< w}$, or $\GG_2$ is a rex move, and can be absorbed into $\LL_2$. We assume henceforth that $\GG_2$ is the identity map, and $\uw' = \uw$.

The source $\uv'$ of $\GG_1$ is a reduced expression for some element $v' \le v$. If $v' < v$, then we can replace $g$ with $\LL_2 \circ (\GG_1 \ot \id_{\uz})$ and $f$ with $\LL_1$, and we have a diagram of the form \eqref{kform} but for $v'$ instead of $v$. By induction, this lives in $H_{\le v'}$ modulo $I_{< w}$. So we reduce to the case when $v' = v$, in which case $\GG_1$ is a rex move, and can be absorbed into $\LL_1$. Thus we have reduced to the following case.
\begin{equation} \label{LLstack} {
\labellist
\small\hair 2pt
 \pinlabel {$\LL_1$} [ ] at 23 16
 \pinlabel {$\LL_2$} [ ] at 33 50
 \pinlabel {$\uv$} [ ] at 4 28
 \pinlabel {$\uw$} [ ] at 4 68
\endlabellist
\centering
\ig{1}{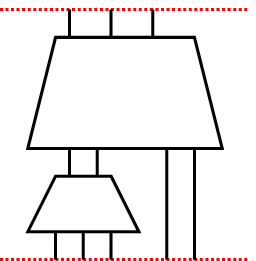}
} \quad = \quad {
\labellist
\small\hair 2pt
 \pinlabel {$\LL_1$} [ ] at 23 17
 \pinlabel {$\LL_{\fb}$} [ ] at 22 40
 \pinlabel {$\TL_{\eb \setminus \fb}$} [ ] at 42 56
\endlabellist
\centering
\ig{1}{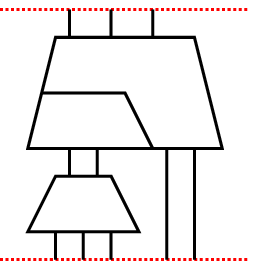}
} \end{equation}

Now $\LL_2$ is associated to some subexpression $\eb \subset \uv \uz$, which restricts to some subexpression $\fb \subset \uv$. If $\fb$ is the top subexpression (the unique
subexpression with terminus $v$) then \eqref{LLstack} is itself a light leaf by \eqref{LLtail}. If $\fb$ is any other expression,
then $\LL_{\fb}$ sends $\uv$ to a reduced expression for an element $v' < v$. Now we again can refactor our diagram, letting $f = \LL_{\fb} \circ \LL_1$, and letting $g$ denote $\TL_{\eb \setminus \fb}$, and can use the inductive hypothesis for $v'$. \end{proof}

Our main application of Proposition \ref{prop:kform} will be to a \emph{mistaken light leaf}, or a \emph{light leaf with error}. Suppose that $\ux = \uy \uz$ is a concatenation, and
$\eb \subset \ux$. In order to construct $\LL_{\eb}$ as in \eqref{LLinduct}, we have already constructed the light leaf $\LL_{\fb}$ associated to the restriction of $\eb$ to $\uy$,
and the light tail $\TL_{\eb \setminus \fb}$ associated to the restriction of $\eb$ to $\uz$. However, instead of gluing these together (along the reduced expression $\uv$) as in
\eqref{LLinduct}, we make an error and insert some morphism $E \in I_{< v}$ as in the following picture.
\begin{equation} \label{mistaken} {
\labellist
\small\hair 2pt
 \pinlabel {$\LL_{\fb}$} [ ] at 23 15
 \pinlabel {$\TL_{\eb \setminus \fb}$} [ ] at 51 47
 \pinlabel {$E$} [ ] at 23 37
\endlabellist
\centering
\ig{1}{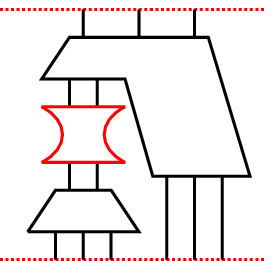}
} \end{equation}
The result is a mistaken light leaf; if $\uy$ has length $k$, we say the mistake happened in the $k$-th tier of the light leaf, and we denote the mistaken light leaf by $\LL_{\eb}^{\oops, k}$.

\begin{lem} \label{lem:propkformtolemma} With the same setup as for \eqref{mistaken}, the mistaken light leaf lies in the subspace $H_{< v}$ of $\Hom(\ux, \uw)/I_{< w}$. \end{lem}

\begin{proof} This is immediate from Proposition \ref{prop:kform}. \end{proof}

This lemma implies that $\LL_{\eb}^{\oops, k}$ is in the span of light leaves $\LL_{\eb'}$ which, at the $k$-th step, factor through elements $< v$ rather than through $v$. However,
we have not yet proven that the Bruhat stroll of $\eb'$ does not go much higher than $\eb$ ever went. Below we will get more control on mistaken light leaves, asserting that $\eb'$
is less than $\eb$ in some partial order on subexpressions. To get this additional control, we can not merely apply Proposition \ref{prop:kform}, but must produce a more subtle
version. The reader should think of Proposition \ref{prop:kform} as a warm-up exercise; the same ideas factor into the proof of Theorem \ref{thm:mistakelexico}.

\begin{defn} Let $\ux = (s_{x_1}, \ldots, s_{x_d})$ be an expression of length $d$, and let $\eb$ and $\eb'$ be two subexpressions of $\ux$. Let $(1 = w_0, w_1, \ldots, w_d)$ and
$(1 = w'_0, \ldots, w'_d)$ be their associated Bruhat strolls, see \cite[\S 2.4]{EWGr4sb}. If $\eb \ne \eb'$, then let $1 \le k \le d$ be the index such that $w_k \ne w_k'$ and $w_j = w_j'$ for all $j > k$; we call $k$ the \emph{index of last difference}. We say that $\eb \prec \eb'$ if $w_k < w_k'$ in the Bruhat order. We write $\eb \preceq \eb'$ if either $\eb = \eb'$ or $\eb \prec \eb'$. We call this the \emph{lexicoBruhat order on subexpressions}. \end{defn}

\begin{lem} The lexicoBruhat order is a total order on the set of subexpressions of $\ux$. \end{lem}

\begin{proof} Clearly this relation is transitive, and $\eb \preceq \eb' \preceq \eb$ implies $\eb = \eb'$. Suppose that $\eb \ne \eb'$, and let $k$ be the index of last difference.
Then $w_k = w_k' s_{x_k}$, so $w_k$ and $w_k'$ are comparable in the Bruhat order. Hence $\eb$ and $\eb'$ are comparable in the lexicoBruhat order. \end{proof}

\begin{rem} Fix $w$, and consider a tail expression $\uz$. Clearly the lexicoBruhat order on subexpressions can be extended to a total order on tail subexpressions of $\uz$ relative to $w$, in the obvious way. \end{rem}

\begin{thm} \label{thm:mistakelexico} Fix $w \in S_n$ with a rex $\uw$, and let $\uz$ be a tail expression of length $m$ relative to $w$. Let $\eb \setminus \fb \subset \uz$ be a tail subexpression with start $v$; choose a rex $\uv$ for $v$, and a light tail $\TL_{\eb \setminus \fb} \co \uv \uz \to \uw$. Let $f \co \uy \to \uv$ be any morphism to $\uv$ which lives in $I_{< v}$. Then
\begin{equation} \TL_{\eb \setminus \fb} \circ (f \ot \id_{\uz}) \in H_{\prec \eb \setminus \fb} \text{ modulo } I_{< w}, \end{equation}
where $H_{\prec \eb \setminus \fb}$ is the left $R_n$-span of light leaves $\LL_{\eb'} \co \uy \uz \to \uw$ whose tail subexpressions $\eb' \setminus \fb' \subset \uz$ satisfy $(\eb' \setminus \fb') \prec (\eb \setminus \fb)$ in the lexicoBruhat order. \end{thm}

Pictorially, we have
\begin{equation} \label{inductivestatementFTW} {
\labellist
\small\hair 2pt
 \pinlabel {$f$} [ ] at 23 19
 \pinlabel {$\TL_{\eb \setminus \fb}$} [ ] at 50 32
 \pinlabel {$\uy$} [ ] at 7 7
 \pinlabel {$\uz$} [ ] at 75 6
 \pinlabel {$\uv$} [ ] at 7 35
 \pinlabel {$\uw$} [ ] at 13 52
 \pinlabel {$< v$} [ ] at 7 20
\endlabellist
\centering
\ig{1}{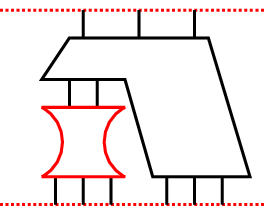}
} \quad \in H_{\prec \eb \setminus \fb} \text{ modulo } I_{< w}. \end{equation}
Before giving the proof, here is the immediate corollary.

\begin{cor} \label{cor:mistakelexico} Fix $\ux = \uy \uz$, where $\uy$ has length $k$. Let $\eb \in E(\ux,w)$, and let $\fb$ be the restriction of $\eb$ to $\uy$. Let $\LL_{\eb}^{\oops, k}$ be any mistaken light leaf where the mistake happens in the $k$-th tier of the light leaf. Then
\begin{equation} \LL_{\eb}^{\oops, k} \in H_{\prec \eb \setminus \fb} \text{ modulo } I_{< w}.\end{equation} \end{cor}

\begin{proof}[Proof of Theorem \ref{thm:mistakelexico}] We prove this by induction on the length $m$ of $\uz$, and within each $m$, by induction on $v$. If $m=0$ then $f \in I_{<
w}$ and the result is trivial. If $m = 1$ then the result actually follows easily Proposition \ref{prop:kform}, but we will use essentially the same argument below for the inductive
step.

Suppose that the result holds for $\uz$, and let us extend $\uz$ by prepending a simple reflection $s$. Suppose we have $(\eb \subset \fb) \subset s \uz$, with start $v$. Let $\eb \setminus \fb'$ be its restriction to $\uz$, which has start $x$. The associated light tail has a first tier $\TL_1$ associated to $s$, and the remaining light tail $\TL_2$ associated to $\uz$. Write $f = f_1 \circ f_2$ where the source of $f_1$ is a reduced expression for some element $v' < v$. We are analyzing the composition
\[ {
\labellist
\small\hair 2pt
 \pinlabel {$f_2$} [ ] at 23 16
 \pinlabel {$f_1$} [ ] at 22 28
 \pinlabel {$\TL_1$} [ ] at 33 48
 \pinlabel {$\TL_2$} [ ] at 67 50
 \pinlabel {$\uy$} [ ] at 3 7
 \pinlabel {$v'$} [ ] at 7 21
 \pinlabel {$v$} [ ] at 2 38
 \pinlabel {$\uw$} [ ] at 5 68
 \pinlabel {$s$} [ ] at 56 8
 \pinlabel {$\uz$} [ ] at 86 10
\endlabellist
\centering
\ig{1}{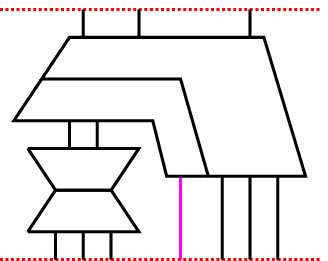}
}. \]
Now apply Proposition \ref{prop:kform} to the subdiagram which is $\TL_1 \circ (f_1 \ot \id_s)$. We get
\begin{equation} \label{inductiveproofFTW} {
\labellist
\small\hair 2pt
 \pinlabel {$f_2$} [ ] at 23 16
 \pinlabel {$f_1$} [ ] at 22 28
 \pinlabel {$\TL_1$} [ ] at 33 48
 \pinlabel {$\TL_2$} [ ] at 67 50
 \pinlabel {$\uy$} [ ] at 3 7
 \pinlabel {$v'$} [ ] at 7 21
 \pinlabel {$v$} [ ] at 2 38
 \pinlabel {$\uw$} [ ] at 5 68
 \pinlabel {$s$} [ ] at 56 8
 \pinlabel {$\uz$} [ ] at 86 10
\endlabellist
\centering
\ig{1}{inductiveFTW}
} \; = \; 
{
\labellist
\small\hair 2pt
 \pinlabel {$f_2$} [ ] at 23 16
 \pinlabel {$\LL$} [ ] at 23 36
 \pinlabel {$\TL_1'$} [ ] at 30 56
 \pinlabel {$\TL_2$} [ ] at 68 56
 \pinlabel {$\uz$} [ ] at 87 9
 \pinlabel {$s$} [ ] at 58 9
 \pinlabel {$\uy$} [ ] at 5 8
 \pinlabel {$v'$} [ ] at 5 26
 \pinlabel {$v''$} [ ] at 5 44
 \pinlabel {$x$} [ ] at 5 64
 \pinlabel {$\uz$} [ ] at 4 78
\endlabellist
\centering
\ig{1}{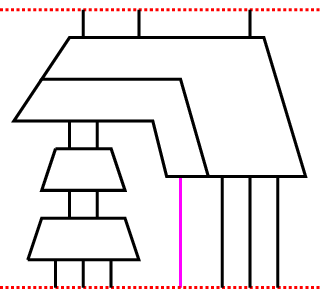}
} \; + \;  {
\labellist
\small\hair 2pt
 \pinlabel {$f_2$} [ ] at 24 15
 \pinlabel {$\TL_2$} [ ] at 67 56
 \pinlabel {$x$} [ ] at 4 59
 \pinlabel {$< x$} [ ] at 7 44
\endlabellist
\centering
\ig{1}{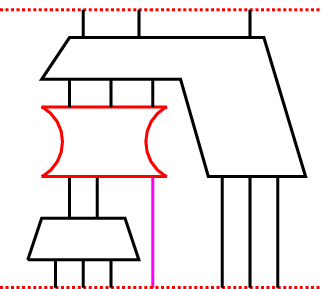}
}. \end{equation}
The RHS is supposed to schematically represent a linear combination (over the left action of $R_n$) of morphisms of two kinds. In the second diagram on the RHS of \eqref{inductiveproofFTW}, the curvy red morphism is supposed to represent the ``lower terms'' of Proposition \ref{prop:kform}. Note that $f_2$ can be absorbed into the lower terms, and the result exactly has the form \eqref{inductivestatementFTW} but for $\uz$ rather than $s \uz$. By induction, this morphism is in the span of $H_{\prec \eb \setminus \fb'}$, which is a subspace of $H_{\prec \eb \setminus \fb}$, modulo $I_{< w}$.

In the first diagram in \eqref{inductiveproofFTW}, $v''$ is some element $< v$, but since it is the target of a light leaf whose source is a rex for $v'$, we must also have $v'' \le v' < v$. Meanwhile, the composition of the two light tails $\TL_2 \circ (\TL_1' \ot \id_{\uz})$ is a light tail whose associated tail subexpression agrees with $\eb \setminus \fb$ on $\uz$, but disagrees on $s$, going to a lower term. Let us resolve $\LL \circ f_2$, noting that light leaves span the maps to $\uv''$ modulo lower terms.
\begin{equation} \label{almostdone} {
\labellist
\small\hair 2pt
 \pinlabel {$f_2$} [ ] at 23 16
 \pinlabel {$\LL$} [ ] at 23 36
 \pinlabel {$\TL_{\prec \eb \setminus \fb}$} [ ] at 57 56
 \pinlabel {$v''$} [ ] at 7 44
\endlabellist
\centering
\ig{1}{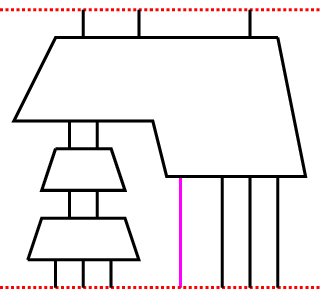}
} \; = \; {
\labellist
\small\hair 2pt
 \pinlabel {$\LL'$} [ ] at 23 24
 \pinlabel {$\TL_{\prec \eb \setminus \fb}$} [ ] at 65 55
 \pinlabel {$v''$} [ ] at 7 44
\endlabellist
\centering
\ig{1}{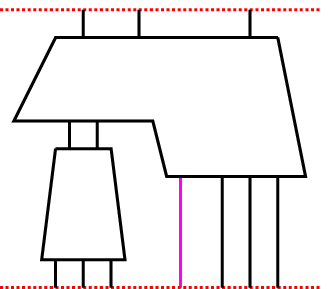}
} \; + \; {
\labellist
\small\hair 2pt
 \pinlabel {$\TL_{\prec \eb \setminus \fb}$} [ ] at 62 54
 \pinlabel {$v''$} [ ] at 2 46
 \pinlabel {$< v''$} [ ] at 2 24
\endlabellist
\centering
\ig{1}{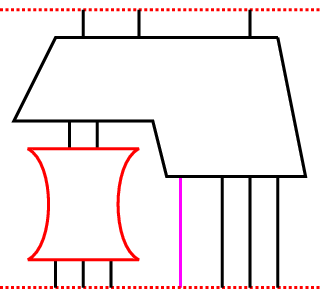}
}. \end{equation}
The first diagram on the RHS of \eqref{almostdone} is a light leaf and lives in $H_{\prec \eb \setminus \fb}$. The second diagram is in $H_{\prec \eb \setminus \fb}$ by induction on $v$, since $v'' < v$. This proves \eqref{inductivestatementFTW} for $s \uz$, completing the inductive step.\end{proof}

\begin{cor} \label{cor:LLequalmodulo} Any two light leaves associated to the same subexpression $\eb \in E(\ux,w)$ are equal in the space $\Hom(\ux, \uw)/I_{< w}$, modulo the left $R_n$-span
of light leaves for smaller subexpressions in the lexicoBruhat order. \end{cor}

\begin{proof} Let $\LL$ and $\LL'$ be the two light leaves associated to $\eb$. The only difference between $\LL$ and $\LL'$ is the choice of rex moves at various tiers in the light
leaves algorithm. Let $\ux$ have length $m$. We can write $\LL - \LL'$ as a telescoping sum
\begin{equation} \LL - \LL' = \LL_0 - \LL_1 + \LL_1 - \LL_2 + \ldots + \LL_{m-1} - \LL_m, \end{equation}
where $\LL_0 = \LL$, $\LL' = \LL_m$, each $\LL_k$ is a light leaf for $\eb$, and $\LL_{k-1}$ differs from $\LL_{k}$ in the choice of rex move made at the $k$-th tier of the algorithm. For example, $\LL_{m-1}$ uses the first $m-1$ tiers of $\LL'$, but uses the last tier from $\LL$.

The difference between two rex moves (between two reduced expressions for some $v \in S_n$) consists of lower terms (i.e. lives in the ideal $I_{< v}$). Thus each difference $\LL_k
- \LL_{k-1}$ is a mistaken light leaf, where the mistake happened in the $k$-th tier. By Corollary \ref{thm:mistakelexico}, $\LL_k - \LL_{k-1}$ is in the span of light leaves
$\LL_{\eb'}$ with $\eb' \prec \eb$. \end{proof}

\begin{rem} This corollary gives an alternate route to proving that distinguished light leaves span the space $\Hom(\ux,\uw)/I_{< w}$. Once one proves that light leaves span, this corollary
proves that a single distinguished light leaf $\LL_{\eb}$ will have the same span as all light leaves associated to $\eb$, modulo $I_{< w}$ and modulo light leaves for $\eb' \prec
\eb$. By induction on the lexicoBruhat order, one deduces that distinguished light leaves have the same span as all light leaves modulo $I_{< w}$. \end{rem}

Let us quickly remark on the difference between the left $R_n$-span and the right $R_n$-span. Modulo lower terms, there is none! Thus all the results above also apply to the right $R_n$-span.

\begin{lem} \label{lem:spansagree} The left $R_n$-span of any set of morphisms in $\Hom(\ux, \uw)$ agrees with the right $R_n$-span modulo $I_{<w}$. \end{lem}
	
\begin{proof} Any polynomial $p$ can be forced through the reduced expression for $w$ at the top of the diagram, using \eqref{polyforce}. The result is $w^{-1}(p)$ on the right hand side, plus error terms where strands are broken. These error terms are all in $I_{< w}$. \end{proof}

Now we bootstrap these results about light leaves to results about double leaves.

\begin{defn} Let $\ux$ and $\uy$ be expressions. Suppose that $(\eb, \fb, w)$ and $(\eb', \fb', w')$ are two coterminal triples subordinate to $(\ux,\uy)$. We say that $(\eb, \fb,
w) \preceq (\eb',\fb',w')$ if $\eb \preceq \eb'$ and $\fb \preceq \fb'$. We call this the \emph{lexicoBruhat order on triples}. \end{defn}

Suppose $w < w'$. If $\ux$ has length $m$ then the index of last difference between $\eb$ and $\eb'$ is $m$, and $w < w'$, so $\eb \prec \eb'$. Similarly, $\fb \prec \fb'$. Thus
$(\eb, \fb, w) \prec (\eb', \fb', w')$ whenever $w < w'$. If particular, the span of all double leaves associated to triples less than $(\eb, \fb, w)$ in the lexicoBruhat order will
have $I_{< w}$ as a subspace.

The definition of a mistaken double leaf $\DL_{\eb, \fb}^{\oops}$ is similar to that of a mistaken light leaf. At one tier in either $\LL_{\eb}$ or $\GG_{\fb}$, one inserts an error term which goes lower in the Bruhat order than it should.

\begin{thm} \label{thm:mistakenDL} Any mistaken double leaf $\DL_{\eb, \fb}^{\oops}$ is in the left $R_n$ span of double leaves associated to triples $(\eb', \fb', w')$ which are smaller than $(\eb, \fb, w)$ in the lexicoBruhat order. \end{thm}

\begin{proof} The error in $\DL_{\eb, \fb}^{\oops}$ is made either in $\LL_{\eb}$ or in $\GG_{\fb}$, and the arguments are the same either way, so let us assume the error is made in $\LL_{\eb}$. Thus
\begin{equation} \DL_{\eb, \fb}^{\oops} = \GG_{\fb} \circ \LL_{\eb}^{\oops}. \end{equation}
By Corollary \ref{thm:mistakelexico}
\begin{equation} \label{eq:lowersumoops}\LL_{\eb}^{\oops} \in \left( \sum_{\eb' \prec \eb} R_n \cdot \LL_{\eb'} \right) + I_{< w}.\end{equation}
Since $I_{< w}$ is a two-sided ideal, $\GG_{\fb} \circ I_{< w} \subset I_{<w}$. Also, $I_{< w}$ lies in the span of smaller double leaves (as discussed a few paragraphs ago). Meanwhile,
\begin{equation} \GG_{\fb} \circ \LL_{\eb'} = \DL_{\fb, \eb'}, \end{equation}
and $(\fb, \eb', w) \prec (\fb, \eb, w)$ when $\eb' \prec \eb$. Thus every term in the sum \eqref{eq:lowersumoops}, when composed with $\GG_{\fb}$, is in the span of lower double leaves in the lexicoBruhat order on triples. \end{proof}

\begin{cor} \label{cor:DLequalmodulo} Any two double leaves associated to the same triple $(\eb, \fb, w)$ are equal in $\Hom(\ux, \uy)$ modulo the span of double leaves for smaller triples in the lexicoBruhat order. \end{cor}

\begin{proof} The proof is the same as Corollary \ref{cor:LLequalmodulo}. \end{proof}

\begin{rem} Once one proves that double leaves span, one deduces from Corollary \ref{cor:DLequalmodulo} that distinguished double leaves have the same span. This gives an
alternative proof of some of the results from \cite[Chapter 7]{EWGr4sb}. \end{rem}

To obtain an analog of Theorem \ref{thm:mistakenDL} for the right $R_n$-action, one can either use Lemma \ref{lem:spansagree} and modify the proof of the theorem accordingly, or one can use the following lemma.

\begin{lem} \label{lem:morespansagree} For any coterminal triple $(\eb, \fb, w)$, the left and right $R_n$ spans of $\{\DL_{\eb', \fb'}\}_{(\eb', \fb') \prec (\eb, \fb)}$ agree.
\end{lem}

\begin{proof} Consider any coterminal triple $(\eb', \fb', w') \prec (\eb, \fb, w)$. Clearly $w' \le w$. Using \eqref{polyforce} to push a polynomial $p$ across the reduced expression for $w'$ in the middle of the double leaf, we get
\begin{equation} \label{eq:leftvsright} \DL_{\eb',\fb'} \cdot p - w'(p)\cdot  \DL_{\eb',\fb'} \in I_{< w'}. \end{equation}
However, $I_{< w'} \subset I_{< w}$ is already in the (right or left) span of lower double leaves. \end{proof}

\subsection{Double leaves and the raising operator}

The key result of this section will be that $\dup(\DL_{\eb, \fb})$ is equal to $p_{\DL}(\eb, \fb) \cdot \DL_{\eb, \fb}$ for some polynomial $p(\eb, \fb)$, plus a linear combination
of mistaken double leaves. By Theorem \ref{thm:mistakenDL}, this means that $\dup$ preserves the span of double leaves (or distinguished double leaves) associated to triples
$\preceq (\eb, \fb, w)$.

\begin{defn} \label{defn:assgrDL} Let $\uy = (s_{y_1}, \ldots, s_{y_d})$ have length $d$, and let $\eb \subset \uy$ be a subexpression, with associated Bruhat stroll $(1 = w_0, w_1, \ldots, w_d)$. Recall from \cite[\S 2.4]{EWGr4sb} that each index of $\eb$ can be decorated as either $U0$, $U1$, $D0$, or $D1$. Define linear polynomials
\begin{subequations} \label{eq:assgrDL}
\begin{equation} p_{\LL}(\eb) = \sum_{\eb_k = U0} w_k(x_{y_i}) - \sum_{\eb_k = D0} w_k(x_{y_k + 1}), \end{equation}
\begin{equation} p_{\GG}(\eb) = \sum_{\eb_k = U0} w_k(x_{y_k+1}) - \sum_{\eb_k = D0} w_k(x_{y_k}). \end{equation}
For a coterminal triple $(\eb, \fb, w)$, let
\begin{equation} p_{\DL}(\eb, \fb) = p_{\LL}(\eb) + p_{\GG}(\fb). \end{equation}
\end{subequations}
\end{defn}

\begin{prop} \label{prop:Heckeassgr} Let $(\eb, \fb, w)$ be a coterminal triple subordinate to $(\ux, \uy)$, and let $\DL_{\eb, \fb}$ be any double leaf for this triple. Then
\begin{subequations}
\begin{equation} \ddown(\DL_{\eb, \fb}) = 0, \qquad \dup(\DL_{\eb, \fb}) = p_{\DL}(\eb, \fb) \cdot \DL_{\eb, \fb} + E, \end{equation}
where $E$ is a linear combination of mistaken double leaves for $(\eb, \fb, w)$. Working instead with right $R_n$-modules we have
\begin{equation} \ddown(\DL_{\eb, \fb}) = 0, \qquad \dup(\DL_{\eb, \fb}) = \DL_{\eb, \fb} \cdot w^{-1}(p_{\DL}(\eb,\fb)) + E. \end{equation}
\end{subequations}
\end{prop}

\begin{proof} Clearly $\ddown$ kills any double leaf, because it is a diagram without polynomials. Let us apply $\dup$ to a double leaf. This is a linear combination of terms where we have applied $\dup$ to each generator in the double leaf, and we analyze each term individually. We prove the result for the left $R_n$ action, as the result for the right action follows by the proof of Lemma \ref{lem:morespansagree}, see \eqref{eq:leftvsright}.

When we apply $\dup$ to a $2m$-valent vertex in a rex move, it produces either zero \eqref{eq:drex0}, or a rex move with broken strands \eqref{eq:d6}. A rex move with broken strands, in the context of the larger diagram, yields a mistaken double leaf.

Whenever $U0$ appears in $\eb$, there is a dot colored $s_i$ in $\LL_{\eb}$ for some $1 \le i \le n-1$. When we apply $\dup$ to this dot, it gets multiplied by $x_i$, see
\eqref{eq:ddots}. If $U0$ appears at the $k$-th step, then this polynomial $x_i$ can then be forced to the left through the reduced expression for $w_k$, using \eqref{polyforce}.
The result will be $w_k(x_i)$ on the left, plus terms with broken strands. These broken strands yield mistaken double leaves.

Similarly, whenever $U0$ appears in $\fb$, there is a dot colored $s_i$ in $\GG_{\fb}$. This time $\dup$ multiplies the dot by $x_{i+1}$, see \eqref{eq:ddots}. The rest of the
argument is the same.

Whenever $D0$ appears in $\eb$, there is a trivalent vertex colored $s_i$ in $\LL_{\eb}$. The raising operator $\dup$ places $-x_{i+1}$ below this trivalent vertex, see \eqref{eq:dtri}. This can then be pulled left through a reduced expression for $w_k$, yielding $-w_k(x_{i+1})$ on the left, together with broken strands. Whenever $D0$ appears in $\fb$, the argument is the same except that $\dup$ places $-x_i$ above the trivalent vertex in $\GG_{\fb}$ instead.

Together, these contributions from $U0$ and $D0$ add up to $p_{\DL}(\eb, \fb)$.

Whenever $U1$ appears in $\eb$ or $\fb$, the corresponding part of the light leaf is the identity map, which is killed by $\dup$.

Whenever $D1$ appears in $\eb$ or $\fb$ (at the $k$-th step), the corresponding part of the light leaf is a cap or cup colored $s_i$. Caps and cups are sent to broken caps and cups
by $\dup$, see \eqref{eq:dcups}. Equivalently, one can break the $s_i$-colored line at the end of the reduced expression for $w_{k-1}$, and then apply a cap or cup. Thus the result is a
mistaken double leaf. \end{proof}

\begin{thm} \label{thm:Heckeassgr} Pick a distinguished double leaf for each coterminal triple subordinate to $(\ux, \uy)$. For any such triple $(\eb, \fb, w)$ we have
\begin{equation} \ddown(\DL_{\eb, \fb}) = 0, \qquad \dup(\DL_{\eb, \fb}) \in \Span \{\DL_{\eb', \fb'}\}_{(\eb', \fb') \preceq (\eb, \fb)}. \end{equation} Note that, by Lemma
\ref{lem:morespansagree}, this span does not depend on whether we work with left or right $R_n$-modules. Thus distinguished double leaves induce a downfree filtration of $\Hom(\ux,
\uy)$ as a left or right $R_n$-module, parametrized by coterminal triples with the lexicoBruhat order. In the associated graded, the left span of $\DL_{\eb, \fb}$ is a copy of $R_n
\ang{p_{\DL}(\eb, \fb)}$, see Definition \ref{defn:assgrDL}. Working instead with right modules, the span of $\DL_{\eb, \fb}$ is a copy of $R_n \ang{w^{-1}(p_{\DL}(\eb, \fb))}$. \end{thm}

\begin{proof} This follows immediately from Proposition \ref{prop:Heckeassgr} and Theorem \ref{thm:mistakenDL}. \end{proof}

\begin{rem} Replacing $\dup$ with $\bardup$ will effectively swap the roles of $\eb$ and $\fb$, obtaining the polynomial $p_{\LL}(\fb) + p_{\GG}(\eb)$ for the left action of $R_n$.
\end{rem}

\section{$\sltwo$ action on Lauda's categorification of $U_q(\sltwo)$}

\subsection{Actions on symmetric polynomials} \label{subsec-sympoly}

The standard action of $\sltwo$ on $R_n$ (see \eqref{actiononRn}) is $S_n$ equivariant, so it descends to an action of $\sltwo$ on $R_n^{S_n}$. Recall that $\pa_i$ is the Demazure operator, 
\begin{equation} \pa_i(f) = \frac{f - s_i f}{x_i - x_{i+1}}, \end{equation}
which sends $R_n$ to $R_n^{s_i}$.

\begin{lem} For any $f \in R_n$ and any $1 \le i \le n-1$ we have
\begin{equation} \label{paddowncommute} \pa_i(\ddown(f)) = \ddown(\pa_i(f)). \end{equation}
\end{lem}

\begin{rem} Note that the corresponding statement for $\dup$ is false! This lemma holds for $\ddown$ effectively because $\ddown$ kills the crossing in the nilHecke algebra, while $\dup$ does not kill the crossing. \end{rem}

\begin{proof} Let us compute in $\NH_n$, with its $\sltwo$ structure from Definition \ref{defn:sl2KLR}. If $X_i$ denotes the $i$-th crossing in $\NH_n$, then we have
\begin{equation} \label{eq:NHforce} X_i f - s_i(f) X_i = \pa_i(f). \end{equation}
Applying $\ddown$ to both sides, and using the fact that $\ddown$ kills $X_i$, we get
\begin{equation} X_i \ddown(f) - \ddown(s_i(f)) X_i = \ddown(\pa_i(f)). \end{equation}
However, $\ddown$ commutes with $s_i$, so we can apply \eqref{eq:NHforce} with $f$ replaced by $\ddown(f)$ to get
\begin{equation} X_i \ddown(f) - s_i(\ddown(f)) X_i = \pa_i(\ddown(f)). \end{equation}
Combining these two equations, we deduce \eqref{paddowncommute}.
\end{proof}

\begin{cor} The operator $\ddown$ commutes with the Demazure operator $\pa_{w_0}$ associated to the longest element of $S_n$. \end{cor}
	
\begin{proof} The Demazure operator $\pa_{w_0}$ can be defined as a composition of operators $\pa_i$ along a reduced expression for $w_0$. Now apply \eqref{paddowncommute}.
\end{proof}

\begin{lem} Let $e_k$ (resp. $h_k$, $p_k$) denote the elementary (resp. complete, power sum) symmetric polynomials in $n$ variables, living in $R_n^{S_n}$. Set $e_{n+1} = 0$. Then
\begin{subequations} \label{eq:sympolyformulas} 
\begin{equation} \label{eq:dupsym} \dup(e_k) = e_k e_1 - (k+1) e_{k+1}, \qquad \dup(h_k) = (k+1) h_{k+1} - h_k h_1, \qquad \dup(p_k) = k p_{k+1}, \end{equation}
\begin{equation} \label{eq:ddownsym} \ddown(e_k) = (n+1-k) e_{k-1}, \qquad \ddown(h_k) = (n-1+k) h_{k-1}, \qquad \ddown(p_k) = \begin{cases} k p_{k-1} & \text{if } k \ne 1 \\ n p_0 & \text{if } k=1 \end{cases}. \end{equation}
Moreover, let $s_{\lambda}$ denote the Schur polynomial, where $\lambda$ is a partition with at most $n$ rows, and let $s_{\lambda} = 0$ when $\lambda$ has more than $n$ rows. Then
\begin{equation} \label{eq:schur} \dup(s_{\lambda}) = \sum_{\mu = \lambda + \square} c(\square) s_{\mu}, \qquad \ddown(s_{\lambda}) = \sum_{\lambda = \mu + \square} (c(\square) + n) s_{\mu}. \end{equation}
Here, $\lambda + \square$ is any partition obtained by adding a single box to $\lambda$, and $c(\square)$ is the content of that box.
\end{subequations}
\end{lem}

\begin{proof} The formulas for $\dup$ can be found in \cite[(2.7)-(2.9)]{EQpDGbig}, except for $\dup(p_k)$ which is easy to check. The formula for $\ddown(p_k)$ is also easy to check.

Let us compute $\ddown(h_k)$. Suppose that $\sum_{i=1}^n a_i = k-1$. There are contributions to the coefficient of $x_1^{a_1} \cdots x_n^{a_n}$ inside $\ddown(h_k)$ coming from $x_i \cdot x_1^{a_1} \cdots x_n^{a_n}$ for each $1 \le i \le n$. We have
\begin{equation} \ddown(x_i \cdot x_1^{a_1} \cdots x_n^{a_n}) = (a_i + 1) x_1^{a_1} \cdots x_n^{a_n} + \text{ terms with different monomials}. \end{equation}
So the overall contribution to the coefficient of $x_1^{a_1} \cdots x_n^{a_n}$ is
\begin{equation} \sum_{i=1}^{n} (a_i + 1) = k-1 + n. \end{equation}	
We leave the computation of $\ddown(e_k)$ to the reader.

To compute $\ddown(s_{\lambda})$, recall that one can define the Schur polynomial as
\begin{equation} s_{\lambda} := \pa_{w_0}(x_1^{\lambda_1} x_2^{\lambda_2} \cdots x_n^{\lambda_n} \cdot x_1^{n-1} x_2^{n-2} \cdots x_n^{0}). \end{equation}
By the above Corollary, $\ddown$ commutes with $\pa_{w_0}$. Applying $\ddown$ to the monomial inside the Demazure operator, we obtain
\begin{equation} \label{ddownonmonomial} \sum_{i=1}^n (\lambda_i + n - i) x_1^{\lambda_1} \cdots x_i^{\lambda_i - 1} \cdots x_n^{\lambda_n} \cdot x_1^{n-1} x_2^{n-2} \cdots x_n^{0}. \end{equation}
If removing a box from the $i$-th row of $\lambda$ yields a partition $\mu$, then $\lambda_i - i$ is the content of that box. Applying $\pa_{w_0}$ to the sum in \eqref{ddownonmonomial}, such terms yield
\[ \sum_{\lambda = \mu + \square} (c(\square) + n) s_{\mu} \]
which is the desired answer. We need only prove that all the other terms in \eqref{ddownonmonomial} are killed by $\pa_{w_0}$. 

If removing a box from the $i$-th row of $\lambda$ is not a partition, it is because $\lambda_i = \lambda_{i+1}$, in which case $\lambda_i - 1 + n - i = \lambda_{i+1} + n - (i+1)$.
Thus the coefficients of $x_i$ and $x_{i+1}$ in the $i$-th monomial from the sum in \eqref{ddownonmonomial} will be equal. Thus this monomial is invariant under $s_i$, and killed by
$\pa_i$ and hence by $\pa_{w_0}$. \end{proof}

\begin{lem} Divided powers $\dup^{(k)}$ and $\ddown^{(k)}$ are defined integrally on $R_n^{S_n}$ for all $k \ge 0$. \end{lem}

\begin{proof} These divided powers are defined on $R_n$, and $R_n^{S_n}$ is an $\sltwo$-invariant subring. \end{proof}

\subsection{Actions on symmetric functions} \label{subsec-symfunc}

It was observed in \cite[end of \S 2.2]{EQpDGbig} that the formulas for $\dup$ acting on symmetric polynomials were independent of the number of variables $n$, and could be extended to
a derivation on the ring $\Lambda$ of symmetric functions (in infinitely many variables). However, it is clear from \eqref{eq:ddownsym} that the formulas for $\ddown$ acting on
symmetric polynomials are not independent of $n$. If we let the number of variables go to infinity, where should $\ddown$ send $e_1 = x_1 + x_2 + \ldots$, given that $\ddown(x_i) =
1$ for all $i$? The appropriate answer is to treat $n$ as a formal variable; this was suggested to the authors by M. Khovanov.

\begin{defn} \label{defn:sl2lambda} Let $\Lambda[y]$ denote the ring of symmetric functions, extended by a formal variable $y$ of degree zero. That is, $\Lambda[y] = \Z[y,e_1,e_2,\ldots]$ is a polynomial ring in infinitely many variables, where $e_k$ has degree $2k$.  We define an $\sltwo$ action on the generators of $\Lambda[y]$ as follows.
\begin{subequations}
\begin{equation} \dup(e_k) = e_k e_1 - (k+1) e_{k+1}, \qquad \dup(y) = 0, \end{equation}
\begin{equation} \label{eq:ddowneklambda} \ddown(e_k) = (y+1-k) e_{k-1}, \qquad \ddown(y) = 0. \end{equation}
\end{subequations}
We extend this action to all of $\Lambda[y]$ using the Leibniz rule.
\end{defn}

\begin{lem} Definition \ref{defn:sl2lambda} gives a well-defined $\sltwo$ action by derivations. The $\sltwo$ action also satisfies the following formulas on familiar elements of $\Lambda[y]$.
\begin{subequations} \label{eq:formulasforLambda}
\begin{equation} \label{eq:duplambda} \dup(e_k) = e_k e_1 - (k+1) e_{k+1}, \qquad \dup(h_k) = (k+1) h_{k+1} - h_k h_1, \qquad \dup(p_k) = k p_{k+1}, \end{equation}
\begin{equation} \label{eq:ddownlambda} \ddown(e_k) = (y+1-k) e_{k-1}, \qquad \ddown(h_k) = (y-1+k) h_{k-1}, \qquad \ddown(p_k) = \begin{cases} k p_{k-1} & \text{if } k \ne 1 \\ y p_0 & \text{if } k=1 \end{cases}. \end{equation}
\begin{equation} \label{eq:schurlambda} \dup(s_{\lambda}) = \sum_{\mu = \lambda + \square} c(\square) s_{\mu}, \qquad \ddown(s_{\lambda}) = \sum_{\lambda = \mu + \square} (c(\square) + y) s_{\mu}. \end{equation}
\end{subequations}
For all $n \ge 1$ there is a surjective $\sltwo$-equivariant quotient map
\begin{equation} \Psi_n \co \Lambda[y] \to R_n^{S_n}, \qquad e_k \mapsto \begin{cases} e_k & \text{if } k \le n, \\ 0 & \text{if } k > n \end{cases}, \qquad y \mapsto n. \end{equation}
\end{lem}

\begin{proof} The existence of $\Psi_n$ for all $n$ is immediate from the formulas for the actions of $\dup$ and $\ddown$. We can use this to efficiently check \eqref{eq:formulasforLambda} using \eqref{eq:sympolyformulas}. For an example of this argument, consider $\ddown(h_k) \in \Lambda[y]$. This will be some $\Z[y]$-linear combination of polynomials in the variables $\{h_i\}$. Because $\Psi_n$ is equivariant, the coefficient of $h_{k-1}$ is some polynomial which evaluates to $n-1+k$ after specializing $y \mapsto n$ for all $n \ge k$. Thus this coefficient is precisely $y-1+k$. Similarly, the coefficient of $h_1 h_{k-2}$ is some polynomial which evaluates to zero after specializing $y \mapsto n$ for all $n \ge k$. Thus this coefficient is zero. \end{proof}
	
\begin{rem} For any $n \in \Z$ there is a quotient of $\Lambda[y]$ which sends $y \mapsto -n$. It is an $\sltwo$-algebra over $\Z$. When $n > 1$, the ideal in this quotient generated by $e_k$ for $k > n$ is $\sltwo$-invariant, and the further quotient by this ideal gives the map $\Psi_n$. When $n \le 0$, there are no nontrivial $\sltwo$-invariant ideals in the quotient. \end{rem}

Sadly, divided powers of $\ddown$ are not defined integrally on $\Lambda[y]$. For example, $\ddown^{(2)}(e_5) = \frac{(y-3)(y-4)}{2} e_3$, and $\binom{y-3}{2} :=
\frac{(y-3)(y-4)}{2}$ is not in $\Z[y]$. The following lemma will imply that divided powers of $\ddown$ are defined integrally whenever $y$ is specialized to an integer, even a
negative integer.

\begin{lem} On $\Lambda[y]$, the divided powers $\dup^{(m)}$ are defined integrally, but $\ddown^{(m)}$ are not defined integrally. However, if we base change over $\Z[y]$ to a larger base ring containing $\binom{y+k}{m}$ for any $k \in \Z$ and $m \ge 0$, then $\ddown^{(m)}$ is defined after base change. \end{lem}

\begin{proof} We need only check that divided powers are well-defined on the generators. The formula for $\dup$ does not involve the variable $y$, so the fact that $\dup^{(m)}(e_k)$ is defined integrally can be checked after applying $\Psi_n$ for sufficiently large $n$. Using \eqref{eq:ddowneklambda} it is easy to verify that
\begin{equation} \ddown^{(m)}(e_k) = \binom{y+1-k}{m} e_{k-m}. \end{equation}   \end{proof}
	
\begin{rem} Note that the divided powers $\ddown^{(m)}$ are always defined integrally on power sums $p_k$, without the need for binomial coefficients. However, power sums generate a different integral form for $\Lambda[y]$ than do the elementary or complete symmetric functions. \end{rem}

\subsection{Rank-one modules and their cores}

In \S\ref{subsec:rank1} we examined $(R_n, \sltwo)$-modules which were free of rank $1$ over $R_n$. Now we do the same for $(R_n^{S_n}, \sltwo)$-modules and $\Lambda[y]$-modules. As before, the rank-one free modules are parametrized by degree two elements in the corresponding ring.

For each $a \in \Z$ define a rank-one module over $R_n^{S_n}$ called $R_n^{S_n}\ang{a e_1}$, with generator $1_a$ living in degree $an$. Equip it with an $\sltwo$ action by the formulas
\begin{equation} \label{eq:actiononrank1forsym}
\dup (1_a) = ae_1 1_a, \quad \quad \ddown (1_a)=0.
\end{equation}

For each $a \in \Z$ define a rank-one module over $\Lambda[y]$ called $\Lambda[y]\ang{a e_1}$, with generator $1_a$ living in degree $ay$. Equip it with an $\sltwo$ action by the same formulas \eqref{eq:actiononrank1forsym}. Note that this $\sltwo$-representation has weights not in $\Z$ but in $\Z[y]$. After specializing $y$ to an integer, this yields an $\sltwo$-representation with weights in $\Z$.


\begin{prop}
If $a> 0$, then $\core(R_n^{S_n} \ang{a e_1})= 0$. If $a\leq 0$ then 
\begin{subequations}
\begin{equation} \label{eq:coreRnSna}
\core(R_n^{S_n}\ang{a e_1}) =  \Z \left\langle
e_{1}^{i_1}\cdots e_n^{i_n}\big| i_1+\cdots +i_n \leq -a
\right\rangle . 
\end{equation}

On the other hand, $\core(\Lambda[y] \ang{a e_1}) = 0$ whenever $a \neq 0$, and 
\begin{equation} \label{eq:coreLambdaa} 
\core(\Lambda[y]) = \Z [y].
\end{equation}
\end{subequations}
\end{prop}

\begin{proof} That $\core(R_n^{S_n}\ang{a e_1}) = 0$ when $a > 0$ follows from Proposition \ref{prop:coreRnp} and the left exactness of taking cores:
\begin{equation}
\core(R_n^{S_n}\ang{a e_1}) \subset \core(R_n \ang{e_1}) = 0.
\end{equation}

The core computation \eqref{eq:coreRnSna} when $a<0$ essentially follows from the proof of \cite[Corollary 2.12]{EQpDGbig}. There, we worked only with $\dup$ in the context of
$p$-dg algebras, but exactly the same computation shows that the right-hand side of \eqref{eq:coreRnSna} is closed under $\dup$. It is clearly closed under $\ddown$, making it a
finite rank submodule. We also showed that the remainder of $R_n^{S_n}\ang{a e_1}$ is acyclic as a $p$-complex, or in other words, after specialization to finite characteristic it
splits into free modules over $\F_p[\dup]/(\dup^p)$. We get the result here by taking the limit $p \to \infty$. More concretely, let $v$ be some element in $R_n^{S_n}\ang{a e_1}$
not in the suspected core; it will be in the core if and only if there is some $N>0$ such that $\dup^N(v) = 0$. However, when $p$ is sufficiently large (e.g. $p > \deg(v) + N$, $p$
does not divide any coefficients of $v$), $v$ descends to a nonzero element in a free module over $\F_p[\dup]/(\dup^p)$, and hence $\dup^N(v) \ne 0$.

Similarly, \eqref{eq:coreLambdaa} in the context of $p$-dg algebras was studied in \cite[Proposition 3.8]{EQpDGsmall}. There, it is shown that $\Lambda \ang{a e_1}$ is acyclic as a $p$-complex whenever $a\neq 0$. When $a = 0$, the augmentation ideal was proven to be acyclic. A similar argument to the previous paragraph, letting $p \to \infty$, will imply the desired results. \end{proof}

\subsection{Action on $\UC$}

\begin{defn} Let $\UC = \UC(\sltwo)$ denote the categorification of quantum $\sltwo$, defined by Lauda in \cite{LauSL2}. We will follow the review given in  \cite[\S 4.1]{EQpDGsmall}. \end{defn}

\begin{defn} \label{defn:sl2U} Place an $\sltwo$ action by derivations on $\UC$ as follows. The raising operator $\dup$ was defined\footnote{In \cite[Definition 5.9]{EQpDGsmall} the raising operator was
called $\pa_1$.} in \cite[Definition 5.9]{EQpDGsmall}. The lowering operator $\ddown$ kills all generators except for dots, and as usual, $\ddown$ sends a dot to the identity map.
\end{defn}

\begin{thm} \label{thm:sl2U} The action of $\sltwo$ on $\UC$ given in Definition \ref{defn:sl2U} is well-defined. The divided power operators $\dup^{(k)}$ and
$\ddown^{(k)}$ are well-defined in the integral form for all $k \ge 0$, making $\UC$ into a divided powers $\sltwo$-algebra. Moreover, inside a region labeled by $\lambda \in \Z$, the subalgebra of $\End(\1_{\lambda})$ generated by bubbles is preserved by the $\sltwo$ action, and is isomorphic as an $\sltwo$-algebra to the specialization $\Lambda[y]/(y-\lambda)$. Under this isomorphism, the clockwise bubble of degree $2k$ is matched with the symmetric function $h_k$, and a counterclockwise bubble of degree $2k$ is matched with $(-1)^k e_k$.  \end{thm}

\begin{proof} First we confirm the action of $\sltwo$ on bubbles. That $\dup$ acts by the formulas in \eqref{eq:duplambda} was already proven in \cite[Corollary 4.8]{EQpDGsmall}. Inside a region labeled $\lambda$, the clockwise bubble with no dots has degree $2(1-\lambda)$. Thus a degree $2k$ bubble has $\lambda + k - 1$ dots. Applying $\ddown$ to a bubble with $\lambda + k - 1$ dots, we get $(\lambda + k - 1)$ times a bubble with $\lambda + k - 2$ dots. This matches the formula
\[ \ddown(h_k) = (y + k - 1) h_{k-1} \]
after specializing $y = \lambda$. Note that a real bubble of negative degree is sent to another real bubble of negative degree, which is still zero.

The computation for counterclockwise bubbles is similar. Consequently, $\ddown$ preserves the bubble relations: positivity of bubbles, normalization of bubbles, and the infinite
Grassmannian relations. We already know that $\ddown$ preserves the nilHecke relations, since the action on the nilHecke algebra is the same as that in Definition \ref{defn:sl2KLR}. We need to check the remaining relations.

The biadjointness relations \cite[(4.1)]{EQpDGsmall} are easy to check.

Consider the ``reduction to bubbles'' or the ``curl relation'' \cite[(4.6a)]{EQpDGsmall}. Up to a sign, the right hand side is
\begin{equation} \sum_{a + b = -\lambda} h_a x^b \end{equation}
where the clockwise bubble $h_a$ appears in region $\lambda$. We compute that
\begin{equation} \ddown(\sum_{a + b = \lambda} h_a x^b) = \sum_{a+b = -\lambda} b h_a x^{b-1} + (\lambda+a-1) h_{a-1} b = \sum_{a+b = -\lambda - 1} h_a x^b (b+1 + \lambda +a) = 0. \end{equation}
Since $\ddown$ kills the left hand side (it has no dots), $\ddown$ preserves the reduction to bubbles relation. The other curl relation \cite[(4.6b)]{EQpDGsmall} is proved similar.

Consider the identity decomposition relation \cite[(4.7b)]{EQpDGsmall}. Using very similar arguments, we compute that
\begin{equation} \ddown(\sum_{a + b + c = -\lambda - 1} x_1^a h_b x_2^c) = \sum_{a + b + c = -\lambda - 2} x_1^a h_b x_2^c (a+1 + \lambda + b + c + 1) = 0. \end{equation}
Hence $\ddown$ preserves the identity decomposition relation. That handles all the relations of $\UC$.

We need to check that $[\ddown,\dup] = \dh$ on each of the generating morphisms of $\UC$, which is immediate from the formulas of \cite[Definition 5.9]{EQpDGsmall}.

We need to check that divided powers are defined integrally on the generating morphisms of $\UC$. For the generators inside a nilHecke algebra, this was already done in Theorem
\ref{thm:sl2KLR}. The cups and caps are killed by $\ddown$. All that remains to check is that $\dup^{(k)}$ is well-defined on the cups and caps. A related question was pursued in
\cite[Lemma A.3 and preceding]{EQpDGsmall}, which proves that $\dup^p = 0$ in characteristic $p$. The proof was to argue that all the coefficients appearing in $\dup^k$ were multinomial
coefficients times $k!$, and this is the same proof\footnote{Admittedly, the proof in \cite[Lemma A.3]{EQpDGsmall} is rather hand-wavey, but the result is still a relatively easy exercise.} needed to show that $\dup^{(k)}$ is defined integrally. \end{proof}

\subsection{Remarks on a downfree filtration}

Let $\EC$ denote the upward strand, an object in $\UC$, and $\FC$ the downward strand. Throughout this section we fix $\lambda \in \Z$ and let $\Lambda$ denote the specialization of $\Lambda[y]$ at $y = \lambda$.

In \cite[\S 8]{LauSL2}, Lauda proves that $\Hom_{\UC}(\1_{\lambda} \EC^n, \1_{\lambda} \EC^n)$ is isomorphic to $\Lambda \ot \NH_n$, where bubbles appear on the left of crossing
diagrams. We can view crossing diagrams as a basis of this Hom space over $\Lambda \ot R_n$, where $R_n$ acts on the bottom. Then this is a downfree filtration over the base ring
$\Lambda \ot R_n$, with the same downfree character computed in \eqref{pwright}. This is because applying $\dup$ to a crossing diagram will not create any bubbles.

Using adjunction, every morphism in the space $\Hom_{\UC}(\EC^n \FC^n \1_{\lambda}, \1_{\lambda})$ is obtained by taking a morphism in $\Hom_{\UC}(\1_{\lambda} \EC^n, \1_{\lambda} \EC^n)$ and
placing caps on top. Thus crossing diagrams will form a basis for $\Hom_{\UC}(\EC^n \FC^n, \1)$ under the action of $\Lambda$ on the left and $R_n$ acting on the inwardly-pointing
boundary strands. This basis is again in bijection with $S_n$. However, the differential of a cap does introduce bubbles. We expect that this produces a downfree filtration over
$\Lambda \ot R_n$, but the formula for the downfree character is currently unknown, and will involve $e_1 \in \Lambda$.

For example, $\Hom_{\UC}(\EC \FC \1_{\lambda}, \1_{\lambda})$ has a basis with one diagram, given by the cap. By \cite[Definition 5.9]{EQpDGsmall}, $\dup$ sends this cup to itself
with a dot minus itself with a degree two bubble. Thus this Hom space, as a module over the $\sltwo$-algebra $R_1 \ot \Lambda$, is isomorphic to $(R_1 \ot \Lambda) \ang{x_1 - e_1}$.

Similar arguments to those used in \cite[\S 8]{LauSL2} will produce a basis for any morphism space in $\UC$. The basis will be a collection of reduced diagrams for oriented planar
matchings with $2n$ boundary points ($n$ oriented in, and $n$ oriented out). It is a basis over the left action of $\Lambda$ and the action of $R_n$ by placing dots on the
inwardly-pointing boundary strands. The basis is in bijection with $S_n$, though this bijection does not preserve the number of crossings in a diagram. We expect that this is a
downfree basis, with a partial order coming from crossing removal; in this case, the partial order does not coincide with the usual Bruhat order on $S_n$. The downfree character is
currently unknown.

\section{$\sltwo$ action on the thick calculus $\dot{\UC}$}

\begin{defn} Let $\dot{\UC}$ denote the $2$-category defined in \cite[\S 4]{KLMS}. \end{defn}

It can be hard to determine from \cite{KLMS} what precisely the generators and relations of $\dot{\UC}$ are, over and above the presentation of $\UC$, as it is not stated
explicitly there. We give the answer in \cite[Proposition 5.2]{EQpDGbig}. There are two new generators for each thickness $a$: the splitter $\EC^{(a)} \to \EC^{\ot a}$ and the merger $\EC^{\ot a} \to \EC^{(a)}$.
\begin{equation} \ig{1}{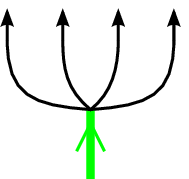} \qquad \ig{1}{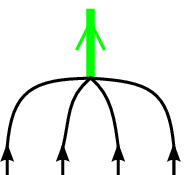} \end{equation} They compose in one direction to be the endomorphism $\psi_{w_0}$ inside $\NH_n$, this is \cite[(5.4a)]{EQpDGbig}.
\begin{equation} \label{eq:mergesplit} \ig{1}{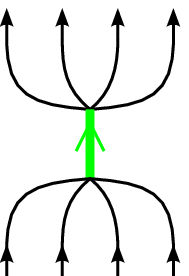} \quad = \quad \ig{1}{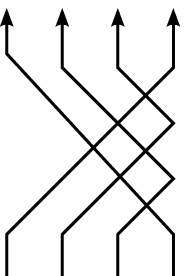} \end{equation} In the other direction, they compose to the identity if sufficiently many dots are placed in between, which is \cite[(5.4b)]{EQpDGbig}. A more general version of \cite[(5.4b)]{EQpDGbig} is
\begin{equation} \label{eq:splitfmerge} {
\labellist
\small\hair 2pt
 \pinlabel {$f$} [ ] at 31 60
\endlabellist
\centering
\ig{1}{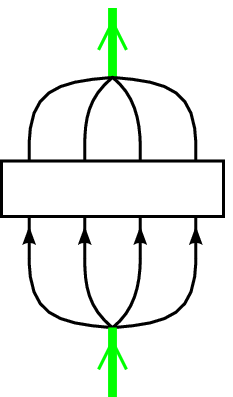}
} \quad = \quad {
\labellist
\small\hair 2pt
 \pinlabel {$\pa_{w_0} f$} [ ] at 15 29
\endlabellist
\centering
\ig{1}{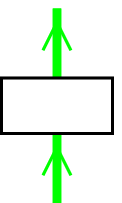}
}. \end{equation}

\begin{defn} We place an $\sltwo$ structure on $\dot{\UC}$, extending the operator $\dup$ from \cite[Definition 5.3]{EQpDGbig}, by asserting that $\ddown$ kills the new generating splitter and merger. \end{defn}

\begin{thm} The $\sltwo$ action on $\dot(\UC)$ is well-defined, and divided powers are defined integrally. \end{thm}

\begin{proof} We need to check that $\ddown$ preserves the relations. Clearly it preserves \eqref{eq:mergesplit} since both sides go to zero. Since $\ddown$ commutes with $\pa_{w_0}$, it is easy to verify that it preserves \eqref{eq:splitfmerge} as well.

We need to check that $[\ddown, \dup] = \dh$ on the new generators. This is straightforward, since $\ddown(-\delta) = \ddown(-\delta')$ agrees with the degree of the splitter. In fact, $\Hom(\EC^{(a)},\EC^{\ot a})$ is isomorphic to $R_a\ang{-\delta}$ as a left $(R_a,\sltwo)$-module, and $\Hom(\EC^{\ot a},\EC^{(a)})$ is isomorphic to $R_a\ang{-\delta'}$ as a right $(R_a,\sltwo)$-module.  Thus divided powers are well-defined by Proposition \ref{prop:Rnp}. For more on $\delta$ and $\delta'$ see the second proof of Theorem \ref{thm:NHmatrix}. \end{proof}

\addcontentsline{toc}{section}{References}


\bibliographystyle{alpha}
\bibliography{mastercopy}

\begin{thebibliography}{EMTW20}

\bibitem[BC18]{BeliakovaCooper}
Anna Beliakova and Benjamin Cooper.
\newblock Steenrod structures on categorified quantum groups.
\newblock {\em Fund. Math.}, 241(2):179--207, 2018.

\bibitem[BY13]{BezYunMonodromy}
Roman Bezrukavnikov and Zhiwei Yun.
\newblock On {K}oszul duality for {K}ac-{M}oody groups.
\newblock {\em Represent. Theory}, 17:1--98, 2013.

\bibitem[DGNO10]{DGNO}
Vladimir Drinfeld, Shlomo Gelaki, Dmitri Nikshych, and Victor Ostrik.
\newblock On braided fusion categories. {I}.
\newblock {\em Selecta Math. (N.S.)}, 16(1):1--119, 2010.

\bibitem[EK10]{EKho}
Ben Elias and Mikhail Khovanov.
\newblock Diagrammatics for {S}oergel categories.
\newblock {\em Int. J. Math. Math. Sci.}, pages Art. ID 978635, 58, 2010.

\bibitem[Eli]{EMSCI}
Ben Elias.
\newblock On diagrammatics for semisimple monoidal categories and their
  integral forms.
\newblock In preparation.

\bibitem[Eli19]{EDiamond}
Ben Elias.
\newblock A diamond lemma for {H}ecke-type algebras.
\newblock Preprint, 2019.
\newblock arXiv 1907:10571.

\bibitem[EMTW20]{EMTW}
Ben Elias, Shotaro Makisumi, Ulrich Thiel, and Geordie Williamson.
\newblock {\em Introduction to {S}oergel bimodules}, volume~5 of {\em RSME
  Springer Series}.
\newblock Springer, 2020.

\bibitem[EQ16a]{EQpDGsmall}
Ben Elias and You Qi.
\newblock An approach to categorification of some small quantum groups {II}.
\newblock {\em Adv. Math.}, 288:81--151, 2016.

\bibitem[EQ16b]{EQpDGbig}
Ben Elias and You Qi.
\newblock A categorification of quantum {$\mathfrak{sl}(2)$} at prime roots of
  unity.
\newblock {\em Adv. Math.}, 299:863--930, 2016.

\bibitem[EQ20]{EQHecke}
Ben Elias and You Qi.
\newblock Categorifying {H}ecke algebras at prime roots of unity, part {I}.
\newblock Preprint, 2020.
\newblock arXiv 2005.03128.

\bibitem[EW14]{EWHodge}
Ben Elias and Geordie Williamson.
\newblock The {H}odge theory of {S}oergel bimodules.
\newblock {\em Ann. of Math. (2)}, 180(3):1089--1136, 2014.

\bibitem[EW16a]{EWRel}
Ben Elias and Geordie Williamson.
\newblock Relative hard {L}efschetz for {S}oergel bimodules.
\newblock preprint, 2016.
\newblock arXiv 1607.03271.

\bibitem[EW16b]{EWGr4sb}
Ben Elias and Geordie Williamson.
\newblock Soergel calculus.
\newblock {\em Represent. Theory}, 20:295--374, 2016.
\newblock arXiv:1309.0865.

\bibitem[Kho16]{KhoHopf}
Mikhail Khovanov.
\newblock Hopfological algebra and categorification at a root of unity: the
  first steps.
\newblock {\em J. Knot Theory Ramifications}, 25(3):1640006, 26, 2016.

\bibitem[Kit18]{KitchlooSteenrod}
Nitu Kitchloo.
\newblock Soergel bimodules, the {S}teenrod algebra, and triply graded link
  homology.
\newblock Preprint, 2018.
\newblock arXiv 1305.4725.

\bibitem[KL09]{KhoLau09}
Mikhail Khovanov and Aaron~D. Lauda.
\newblock A diagrammatic approach to categorification of quantum groups. {I}.
\newblock {\em Represent. Theory}, 13:309--347, 2009.

\bibitem[KL11]{KhoLau11}
Mikhail Khovanov and Aaron~D. Lauda.
\newblock A diagrammatic approach to categorification of quantum groups {II}.
\newblock {\em Trans. Amer. Math. Soc.}, 363(5):2685--2700, 2011.

\bibitem[KLMS12]{KLMS}
Mikhail Khovanov, Aaron~D. Lauda, Marco Mackaay, and Marko Sto{\v{s}}i{\'c}.
\newblock Extended graphical calculus for categorified quantum {${\rm sl}(2)$}.
\newblock {\em Mem. Amer. Math. Soc.}, 219(1029):vi+87, 2012.

\bibitem[KQ15]{KQ}
M.~Khovanov and Y.~Qi.
\newblock An approach to categorification of some small quantum groups.
\newblock {\em Quantum Topol.}, 6(2):185--311, 2015.
\newblock \href{http://arxiv.org/abs/1208.0616}{arXiv:1208.0616}.

\bibitem[KR16]{KRWitt}
M.~Khovanov and L.~Rozansky.
\newblock Positive half of the {W}itt algebra acts on triply graded link
  homology.
\newblock {\em Quantum Topol.}, 7(4):737--795, 2016.
\newblock \href{https://arxiv.org/abs/1305.1642}{arXiv:1305.1642}.

\bibitem[KSQ17]{KQS}
Mikhail Khovanov, Joshua Sussan, and You Qi.
\newblock p-{DG} cyclotomic nil{H}ecke algebras.
\newblock Preprint, 2017.
\newblock arXiv 1711.07159.

\bibitem[Lau10]{LauSL2}
Aaron~D. Lauda.
\newblock A categorification of quantum {${\rm sl}(2)$}.
\newblock {\em Adv. Math.}, 225(6):3327--3424, 2010.

\bibitem[Lau20]{LaudaParameters}
Aaron~D. Lauda.
\newblock Parameters in categorified quantum groups.
\newblock {\em Algebr. Represent. Theory}, 23(4):1265--1284, 2020.

\bibitem[Lib08]{LibLL}
Nicolas Libedinsky.
\newblock Sur la cat\'egorie des bimodules de {S}oergel.
\newblock {\em J. Algebra}, 320(7):2675--2694, 2008.

\bibitem[LS12]{LicataSavageHeisenberg}
Anthony Licata and Alistair Savage.
\newblock A survey of {H}eisenberg categorification via graphical calculus.
\newblock {\em Bull. Inst. Math. Acad. Sin. (N.S.)}, 7(2):291--321, 2012.

\bibitem[Mon93]{Montgomerybook}
Susan Montgomery.
\newblock {\em Hopf algebras and their actions on rings}, volume~82 of {\em
  CBMS Regional Conference Series in Mathematics}.
\newblock Published for the Conference Board of the Mathematical Sciences,
  Washington, DC, 1993.

\bibitem[MS89]{ManSch}
Yu.~I. Manin and V.~V. Schechtman.
\newblock Arrangements of hyperplanes, higher braid groups and higher {B}ruhat
  orders.
\newblock In {\em Algebraic number theory}, volume~17 of {\em Adv. Stud. Pure
  Math.}, pages 289--308. Academic Press, Boston, MA, 1989.

\bibitem[MSV13]{MSV}
Marco Mackaay, Marko Sto{\v{s}}i{\'c}, and Pedro Vaz.
\newblock A diagrammatic categorification of the {$q$}-{S}chur algebra.
\newblock {\em Quantum Topol.}, 4(1):1--75, 2013.

\bibitem[Qi14]{QiHopf}
You Qi.
\newblock Hopfological algebra.
\newblock {\em Compositio Mathematica}, 150(01):1--45, 2014.

\bibitem[QS18]{QiSussan3}
You Qi and Joshua Sussan.
\newblock p-{DG} cyclotomic nil{H}ecke algebras {II}.
\newblock Preprint, 2018.
\newblock arXiv 1811.04372.

\bibitem[Rou08]{Rouq2KM-pp}
Rapha{\"e}l Rouquier.
\newblock 2-{K}ac-{M}oody algebras.
\newblock Preprint, 2008.
\newblock arXiv:0812.5023.

\end{thebibliography}

%

\vspace{0.1in}

\noindent B.~E.:  { \sl \small Department of Mathematics, University of Oregon, Eugene, OR 97403, USA}\newline \noindent {\tt \small email: belias@uoregon.edu}

\vspace{0.1in}

\noindent Y.~Q.: { \sl \small Department of Mathematics, University of Virginia, Charlottesville, VA 22904, USA} \newline \noindent {\tt \small email: yq2dw@virginia.edu}

%
\end{document}